\documentclass[11pt]{amsart}
\usepackage{amsthm,amsfonts,amssymb,amsmath,oldgerm}
\numberwithin{equation}{section}
\usepackage{fullpage}
\usepackage{amsmath}
\usepackage{setspace}
\usepackage{stmaryrd}
\usepackage{mathrsfs}
\usepackage{fancyhdr}
\usepackage{graphicx}
\usepackage{enumerate}
\usepackage{bm}
\usepackage{bbm}
\usepackage{dsfont}
\usepackage{multirow}
\usepackage{psfrag}
\usepackage[font=scriptsize]{caption}
\usepackage[font=scriptsize]{subcaption}
\usepackage{listings}
\usepackage{tikz}
\usepackage{empheq}
\usepackage{cases}
\usetikzlibrary{matrix} 
\usepackage[framemethod=TikZ]{mdframed}
\mdfdefinestyle{MyFrame}{backgroundcolor=gray!20!white}
\usepackage[makeroom]{cancel}

\usepackage[top=20mm,bottom=20mm,left=16mm,right=16mm,a4paper]{geometry}

\usepackage{hyperref,bookmark}

\usepackage[normalem]{ulem}
\definecolor{violet}{rgb}{0.580,0.,0.827}

\usepackage[textwidth= \marginparwidth,textsize=tiny]{todonotes}

\hypersetup{
    colorlinks,          
    citecolor=blue,        
    filecolor=blue,      
    urlcolor=blue,           
    linkcolor=blue,
}

\newcommand\dD{\mathrm{d}}

\def\eps{\varepsilon }

\newcommand{\J}{{\mathcal J}}

\newcommand{\beq}{\begin{equation}}
\newcommand{\eeq}{\end{equation}}
\newcommand{\beqa}{\begin{eqnarray}}
\newcommand{\eeqa}{\end{eqnarray}}
\newcommand\br{\begin{remark}}
\newcommand\er{\end{remark}}
\newcommand\bp{\begin{pmatrix}}
\newcommand\ep{\end{pmatrix}}
\newcommand{\be}{\begin{equation}}
\newcommand{\ee}{\end{equation}}
\newcommand\ba{\begin{equation}\begin{aligned}}
\newcommand\ea{\end{aligned}\end{equation}}
\newcommand\ds{\displaystyle}

\newcommand{\beg}{\begin{example}}
\newcommand{\eeg}{\end{exaplem}}
\newcommand{\bpr}{\begin{proposition}}
\newcommand{\epr}{\end{proposition}}
\newcommand{\bt}{\begin{theorem}}
\newcommand{\et}{\end{theorem}}
\newcommand{\bc}{\begin{corollary}}
\newcommand{\ec}{\end{corollary}}
\newcommand{\bl}{\begin{lemma}}
\newcommand{\el}{\end{lemma}}
\newcommand{\bd}{\begin{definition}}
\newcommand{\ed}{\end{definition}}
\newcommand{\brs}{\begin{remarks}}
\newcommand{\ers}{\end{remarks}}

\newtheorem{theorem}{Theorem}[section]
\newtheorem{proposition}[theorem]{Proposition}
\newtheorem{corollary}[theorem]{Corollary}
\newtheorem{lemma}[theorem]{Lemma}
\newtheorem{remark}[theorem]{Remark}
\newtheorem{definition}[theorem]{Definition}

\newtheorem{example}[theorem]{Example}

\newcommand{\N}{{\mathbb N}}

\newcommand{\R}{{\mathbb R}}
\newcommand{\T}{{\mathbb T}}

\newcommand\bx{{\bm x}}

\newcommand\bv{{\bm v}}
\newcommand\Div{{\rm div}}

\newcommand\bE{{\mathbf E}}

\newcommand\cA{{\mathcal A}}
\newcommand\cB{{\mathcal B}}

\newcommand\cD{{\mathcal D}}
\newcommand\cE{{\mathcal E}}

\newcommand\cH{{\mathcal H}}
\newcommand\cI{{\mathcal I}}
\newcommand\cJ{{\mathcal J}}

\newcommand\cM{{\mathcal M}}

\newcommand\cR{{\mathcal R}}

\usetikzlibrary{fit}
\usetikzlibrary{arrows}
\numberwithin{equation}{section}
\hypersetup{linkbordercolor=red}

\email{francis.filbet@math.univ-toulouse.fr}
\email{alain.blaustein@math.univ-toulouse.fr}

\title[On a discrete framework of hypocoercivity for kinetic
equations]{On a discrete framework of hypocoercivity for kinetic
  equations} 

\author{Alain Blaustein and Francis Filbet}



\begin{document}

\maketitle

\centerline{\scshape Alain Blaustein}
\medskip
{\footnotesize
     \centerline{Institut de Math\'ematiques de Toulouse, Universit\'e Paul Sabatier}
    \centerline{Toulouse, France}
  }
  \medskip

\centerline{\scshape Francis Filbet}
\medskip
{\footnotesize
    \centerline{Institut de Math\'ematiques de Toulouse, Universit\'e Paul Sabatier}
    \centerline{Toulouse, France}
}

\bigskip

\begin{abstract}
We propose and study a fully discrete finite volume scheme  for the
Vlasov-Fokker-Planck equation written as an hyperbolic system using
Hermite polynomials in velocity.  This approach naturally preserves
the stationary solution and the weighted $L^2$ relative entropy. Then,
we adapt the arguments developed in \cite{DMS} based the
hypocoercivity method to get quantitative estimates on the convergence
to equilibrium of the discrete solution. Finally, we prove that in the diffusive limit, the
scheme is asymptotic preserving with respect to both the time variable and the scaling parameter at play. 
\end{abstract}
\vspace{0.1cm}
\tableofcontents

\section{Introduction}
\label{sec:1}
\setcounter{equation}{0}
\setcounter{figure}{0}
\setcounter{table}{0}
The Vlasov-Fokker-Planck  equation is the kinetic description of the
Brownian motion of a large system of charged particles under the
effect of an electric field. For example, in electrostatic plasma,
where the Coulomb force are taken
into account, the time evolution of the electron distribution function
$f$ solves the Vlasov-Poisson-Fokker-Planck system, under the action of a self-consistent
potential $\Phi$:
\begin{equation*}
  \left\{
    \begin{array}{l}
\ds\frac{\partial f}{\partial t}\,+\,\bv\cdot\nabla_\bx f
      \,+\, \frac{q_e}{m_e}\,\bE\cdot\nabla_\bv f \,=\,
      \frac{1}{\tau_e}\,\Div_\bv\left(\bv f \,+\, T_0 \,\nabla_\bv f \right)\,,
      \\[1.1em]
      \ds-\eps_0\Delta\Phi \,=\, q_e\int_{\R^3} f \dD \bv, 
      \end{array}\right.
  \end{equation*}
  where $\eps_0$ is the vacuum permittivity, $q_e$ and $m_e$ are
  elementary charge and mass of the electrons, whereas $\tau_e$ is the
  relaxation time due to the collisions of the particles with the
  surrounding bath.

 Considering $\eps > 0$ as the ratio between the mean free path of
 particles and the length scale of observation, it allows to identify
 different regimes and the Vlasov equation
 may be written in a adimensional form
 \be
  \label{vfp0}
\eps\,\frac{\partial f}{\partial t}\,+\,\bv\cdot\nabla_\bx f
      \,+\, \bE\cdot\nabla_\bv f \,=\,
      \frac{\eps}{\tau(\eps)}\,\Div_\bv\left(\bv f \,+\, T_0 \,\nabla_\bv f \right)\,,
 \ee
 Our main purpose here is to build and analyse a numerical scheme able
 to capture two regimes of interest for equation \eqref{vfp0}, in a
 linear framework: the long time behavior $t\rightarrow \infty$ and
 the diffusive regime $\eps\rightarrow 0$. In
 various situations, the scaling parameters at play may be non
 homogeneous across the system leading to intricate situations, where
 both processes may coexist. Thus, we aim at designing a scheme robust
 enough to capture simultaneously these different behaviors.

More precisely, we consider the one dimensional Vlasov-Fokker-Planck equation with periodic boundary conditions in space, which reads
\beq
\label{vlasov}
\partial_t f \,+\, \frac{1}{\eps}\left( v\,\partial_x f
\,+\,E\,\partial_v f\right)  \,=\, \frac{1}{\tau(\eps)}\,
\partial_v \left( v \, f \,+\, T_0\,\partial_v f \right)\,,
\eeq
with $t\geq 0$, position $x \in \T$ and velocity
$v\in\R$, whereas the electric field derives from a potential $\Phi$
such that $E=-\partial_x \Phi$, with the following regularity assumption
\begin{equation}\label{reg:Phi}
	\Phi\in W^{2,\infty}\left(\T\right)\,.
\end{equation}
We also define the density $\rho$ by integrating the distribution
function in velocity,
\be
\label{def:rho}
\rho(t,x) \,=\, \int_{\R} f(t,x,v) \,\dD v. 
\ee

It is worth to mention that there are already several works on preserving large-time
behaviors of solutions to the Fokker-Planck equation or related
kinetic models. On the one hand, a fully discrete finite difference
scheme for the homogeneous Fokker-Planck equation has been proposed in
the pioneering work of Chang and Cooper \cite{ChangCooper:1970}. This
scheme preserves the stationary solution and the entropy decay of the
numerical solution.  On the other hand, finite volume schemes preserving the
exponential trend to equilibrium have been studied for non-linear
convection-diffusion  equations (see for example
\cite{Bessemoulin-Chatard2012,Burger2010,Chainais-Hillairet2007,Gosse2006}). More
recently, in  \cite{Pareschi2017}, the authors investigate the
question of describing correctly the equilibrium state of  non-linear
diffusion and kinetic models for high order schemes. Let us also
mention some works on boundary value problems \cite{Filbet2017,chainais_2020} where non-homogeneous Dirichlet boundary conditions are dealt with. 
  
  In the case of space non homogeneous kinetic equations, the convergence to
  equilibrium becomes tricky because of the lack of coercivity since
  dissipation occurs only in the velocity variable whereas transport
  acts in the space variable. Therefore,  only few results are
  available and a better understanding of hypocoercive structures at
  the discrete level is challenging.  Let us mention a first rigorous
  work  in this direction on the  Kolmogorov equation \cite{porretta_2017_numerical,Foster2017,Georgoulis}.
  In \cite{Foster2017}, a time-splitting scheme is applied and it is
  shown that  solutions  decay polynomially in time. In
  \cite{porretta_2017_numerical,Georgoulis},  a different approach has
  been used, based on the  work of H\'erau \cite{herau2017} and Villani \cite{Villani:AMS},
  for finite difference and a finite element schemes. Later,  Dujardin, H\'erau and Lafitte \cite{Dujardin} studied a finite difference scheme for the kinetic Fokker-Planck equation. Finally, in a more recent work \cite{MMT}, the authors established a
     discrete hypocoercivity framework based on the continuous
     approach provided in \cite{DMS}. It is based on  a modified discrete
     entropy,  equivalent to a weighted $L^2$ norm involving
     macroscopic quantities and the authors show quantitative
     estimates on the numerical solution for large time and in the
     limit $\eps\rightarrow 0$.

     The present contribution can be
     considered as a continuation of this latter work in order to discretize  the
     kinetic Fokker-Planck equation with an applied force field.  On the one hand, we consider the case where the interactions associated to collisions and
electrostatic effects have the same magnitude, that is,
$\tau(\eps)\sim\eps$, hence the limit $t/\eps\,\rightarrow\,+\infty$ corresponds to the
long time behavior of equation
\eqref{vlasov}. In this regime, the distribution function $f$
relaxes towards the stationary solution to the Vlasov-Fokker-Planck
equation $\rho_\infty \,\cM$, where  the Maxwellian $\cM$  is given by
$$
\cM(v)\,=\, \frac{1}{\sqrt{2\pi \,T_0}} \, \exp\left( -\frac{|v|^2}{2\,T_0}\right)\,,
$$
whereas the density $\rho_\infty$ is determined by
\be
\label{def:rho_inf}
\rho_\infty  \,=\, c_0 \, \exp\left(-\frac{\Phi}{T_0} \right),
\ee
where the constant $c_0$ is fixed by the conservation of mass,
that is,
$$
\int_{\T} \rho_\infty \,\dD x \,=\,  \iint_{\T\times \R} f_0(x,v)
\,\dD v\dD x\,.
$$
Thus, we set $f_\infty$  the stationary state of \eqref{vlasov},  defined as
$$
f_\infty(x,v) \,=\, \rho_\infty(x) \,  \cM(v)
$$ 
and we expect  that $f\rightarrow f_{\infty}$ as $t/\eps \rightarrow +\infty$.\\

On the other hand, the diffusive regime corresponds to a frontier where collisions dominate but still not enough to cancel completely the electrostatic effects. This situation occurs as $\eps\,\rightarrow\,0$ in the case where $\tau(\eps)\sim\tau_0\,\eps^2$, for some $\tau_0 >0$. Due to collisions, the distribution of velocities also relaxes towards a Maxwellian equilibrium. However, in this case, the spatial distribution converges to a time dependent distribution $\rho$ whose dynamics are driven by a drift-diffusion equation depending on the force field $E$. Indeed, performing the change of variable $x \rightarrow x \,+\,\tau_0\,\eps\,v$ in \eqref{vlasov} and integrating with respect to $v$, we deduce that the quantity
\[
\pi\left(t,x
\right)
\,=\,
\int_{\R}
f\left(t,
x-\tau_0\,\eps\,v
,v
\right)\,\dD v\,,
\]
solves the following equation
\begin{equation*}
	\ds\partial_t \, \pi
	\,+\,
	\tau_0\,\partial_x \left(\int_{\R}
	E\,f
	\left(t,
	x-\tau_0\,\eps\,v
	,v
	\right)
	\,
	\dD v\,-\,T_0\,\partial_x\,\pi\right)
	\,=\,
	0
	\,.
\end{equation*}
According to its definition, $\pi$ verifies:
$
\ds
\rho\,\sim\,\pi
$ in the limit $\eps\,\rightarrow\,0$. Therefore, we may formally replace
$
\ds\pi
$ with $\rho$ and $\eps$ with $0$ in the latter equation. This yields
\[
f(t,x,v)\,\underset{\eps\rightarrow0}{\longrightarrow}\,\rho_{\tau_0}(t,x) \,  \cM(v)
\,,
\]
where $\rho_{\tau_0}$ solves
\beq
\label{drift:diffusion}
\partial_t \rho_{\tau_0} \,+\,\tau_0\,\partial_x \left(E\,\rho_{\tau_0}\,-\,T_0\,\partial_x\,\rho_{\tau_0}\right)\,=\,0\,.
\eeq
To be noted that this regime is an intermediate situation which
contains more information than the long time asymptotic since we have
$\rho \rightarrow \rho_{\infty}$ by taking either $t \rightarrow
+\infty$ or $\tau_0 \rightarrow +\infty$.

At the discrete level, Asymptotic-Preserving
schemes have been developed to capture in a discrete setting the
diffusion limit, so that in the limit $\eps\rightarrow 0$, the
numerical discretization converges to  the macroscopic model (see
for instance \cite{jin2000, liu, jin2010,lemou2008} on finite
difference and finite volume schemes and \cite{dimarco2018, crestetto2017} on particle methods). 
 
In the present article, our aim is to design a numerical scheme which
is able to capture these two regimes but also all the intermediate
situations where $\ds\eps^2\,\lesssim\,
\tau(\eps)\,\lesssim\,\eps$. More precisely, we  suppose that
\begin{equation}
\label{cond:tau}
\sup_{\eps>0}\frac{\tau(\eps)}{\eps} \,\leq \, \overline{\tau}_0\in \left(0\,,\,+\infty\right).
\end{equation}
and distinguish two cases on $\tau(\eps)$ : 
\begin{itemize}
\item[$(i)$] either the diffusive regime assumption
\begin{equation}
\label{cond2a:tau}
\frac{\tau(\eps)}{\eps^2}\,\underset{\eps\rightarrow0}{\longrightarrow}\, \tau_0 < \,+\infty\,,
\end{equation}
where collisional effects strongly dominate; 
\item[$(ii)$] or the intermediate regime assumption
\begin{equation}
 \label{cond2b:tau}
	\frac{\tau(\eps)}{\eps^2}\,
	\underset{\eps\rightarrow 0}{\longrightarrow}\,+\infty\,,
      \end{equation}
which may for instance correspond to $\tau(\eps)=\eps^{\beta}$, with
$1\leq\beta<2$. It describes all the intermediate situations between long time and
diffusive regime.
\end{itemize}

The starting point of our analysis is the following estimate, obtained multiplying equation \eqref{vlasov} by $f\,/\,f_{\infty}$, and balancing the transport term with
the source term corresponding to the electric field thanks to the weight $f_{\infty}^{-1}$
\be
\label{esti:L2}
\frac{1}{2} \frac{\dD}{\dD t} \int_{\T^d\times\R^d}|f-f_\infty|^2\,f_{\infty}^{-1}\dD v\,\dD x \,+\,\frac{T_0}{\tau(\eps)}\,  \int_{\T\times\R} \left|
  \partial_v\left( \frac{f}{f_\infty}\right) \right|^2 \, f_\infty
\,\dD v\,\dD x \,=\, 0\,.
\ee
This estimate is important since it yields a $L^2$ stability result
on the solution to the Vlasov-Fokker-Planck equation
\eqref{vlasov}.

Our  purpose is to design a numerical scheme for which such estimate occurs. To this aim, we split our approach in two
steps: we apply a spectral  decomposition in velocity  of $f$
based on Hermite decomposition and we apply a structure preserving
finite volume scheme for the space discretization.  In the next
section (Section \ref{sec:2}), we provide explicit convergence rates for the continuous
model written in the Hermite basis (see Theorems \ref{th:main1} and
\ref{th:main2}). This first step allows us to present the general
strategy and to highlight the main properties of the transport
operator in order to design suitable numerical scheme.  Therefore, in
Section \ref{sec:3} we  adapt these latter results without any
loss to the fully discrete setting using a structure preserving finite
volume scheme and an implicit Euler scheme for the time discretization
(see Theorems \ref{th:main1:h} and \ref{th:main2:h}).  The variety of situations that we aim to cover may lead to various and intricate behaviors. Therefore, we successfully put great efforts into providing results which are uniform with respect to all parameters at play: time $t$, scaling parameters $(\eps,\tau_0)$ and eventually the numerical discretization. The result is worth the pain, since we propose in the Section \ref{sec:4} various simulations, in which we are able to capture, at low computational cost, a rich variety of situations.

%
\section{Hermite's decomposition for the velocity variable}
\label{sec:2}
\setcounter{equation}{0}
\setcounter{figure}{0}
\setcounter{table}{0}

The purpose of this section is to present a formulation of the
Vlasov-Fokker-Planck equation \eqref{vlasov} based on Hermite
polynomial and to provide quantitative results on $f$ when
$\eps\rightarrow 0$ and $t \rightarrow +\infty$. These results are identical to the ones obtained
in the continuous case except that there are formulated on the
corresponding Hermite's coefficients solution to a linear hyperbolic
system. This formulation  is well  adapted to prepare the fully discrete setting in Section \ref{sec:3}.
\\
We first use Hermite polynomials in the velocity
variable and write the Vlasov-Fokker-Planck equation
\eqref{vlasov}  as an infinite hyperbolic system for the Hermite coefficients
depending only on time and space. The idea is to apply a Galerkin
method only keeping a small finite set of orthogonal
polynomials rather than discretizing the distribution function in
velocity \cite{Armstrong1967,Joyce1971}. The merit to use
orthogonal basis like the so-called scaled Hermite basis has been
shown in \cite{Holloway1996,holloway2,Schumer1998} or more recently
\cite{Filbet2020,ref:5} for the Vlasov-Poisson system.  In this context the family of Hermite's functions  $\left(\Psi_k\right)_{k\in\N}$ defined as
\[
\Psi_{k}(v)\,=\, H_{k}\left(\frac{v}{\sqrt{T_0}}\right)\,\cM(v)\,,
\]
constitutes an orthonormal system for the inverse Gaussian weight, that is,
\[
\int_{\R}\,
\Psi_{k}(v)\,\Psi_{l}(v)\,\cM^{-1}(v)\dD v
\,=\,
\delta_{k,l}\,.
\]
In the latter definition, $
\ds
\left(
H_{k}
\right)_{k\in\N} 
$ stands for the family of Hermite polynomials defined recursively as follows
$H_{-1}=0$, $H_{0}=1$ and
\[
\xi\,H_{k}(\xi)\,=\,
\sqrt{k}\,H_{k-1}(\xi)
        \,+\,
        \sqrt{k+1}\,H_{k+1}(\xi)
\,,\quad\forall\, k\,\geq\,0\,.
\]
Let us also point out that Hermite's polynomials verify the following relation
\[
H_k'(\xi)\,=\,\sqrt{k}\,H_{k-1}(\xi)\,,\quad\forall\, k\,\geq\,0\,.
\]
Taking advantage of the latter relations, one can see why Hermite's functions arise naturally when studying the Vlasov-Poisson-Fokker-Planck model, especially in the diffusive regime, as they constitute an orthonormal basis which diagonalizes the Fokker-Planck operator:
\[
\partial_{v}\left[\,
v\,\Psi_{k}
\,+\,
T_0\,\partial_{v}\,\Psi_{k}\,
\right]
\,=\,
\,-\,k
\,\Psi_{k}\,.
\]
Therefore, we consider the decomposition of $f$ into its components 
$
C\,=\,
\left(
C_{k}
\right)_{k\in\N}
$
in the Hermite basis
\begin{equation}\label{f:decomp}
f
\left(t,x,v\right)
\,=\,
\sum_{k\in\N}\,
C_{k}
\left(t,x
\right)\,\Psi_{k}(v)\,.
\end{equation}
It's worth to mention that we also may consider a truncated series
neglecting high order coefficient in order to construct a
spectrally accurate approximation of $f$ in the velocity
variable.

As we have shown before, Hermite's decomposition with respect to the velocity variable is a suitable choice in our setting. When it comes to the space variable, we see from estimate \eqref{esti:L2} that the natural functional framework here is the $L^2$ space with weight $\rho_{\infty}^{-1}$. Unfortunately, it is not very  well adapted to  the space
discretization since it may generate additional spurious terms
difficult to control when dealing with
discrete integration by part. We bypass this difficulty by integrating the weight in the quantity of interest: instead of working directly with $\ds f$, we consider the quantity $\ds f\,/\,\sqrt{\rho}_{\infty}$ in order to get a well-balanced scheme
in the same spirit to what has been already done in \cite{chainais_2020,Filbet2017} for well-balanced finite volume schemes. More precisely, we set
$$
D_{k} \,:=\, \frac{C_{k}}{\sqrt\rho_\infty} 
$$
in \eqref{f:decomp}, and inject this ansatz in \eqref{vlasov}. Using that $\rho_\infty \, E \,=\, T_0\, \partial_x\rho_\infty $, we get
that
$D\,=\,(D_{k})_{k \in \N}$ satisfies the following system
  \begin{equation}
  \label{Hermite:D}
  \left\{
    \begin{array}{l}
  \ds\partial_t  D_{k} \,+\,
      \frac{1}{\eps}\,\left( \sqrt{k}\,
      \cA\,D_{k-1}\,- \, \sqrt{k +1}\,
        \cA^\star D_{k+1}\right) \,  =\,  - \,  \frac{k}{\tau(\eps)}\,  D_{k}\,,
        \\[1.2em]
         \ds D_{k}(t=0) =  D^{0,\eps}_{k}\,,
\end{array}\right.
\end{equation}
where operators $\cA$ and $\cA^{\star}$ are given by
\begin{equation*}
\left\{
\begin{array}{l}
\ds  \cA \,u  \,=\, +\sqrt{T_0}\,\partial_{x}  u \,-\,
\frac{E}{2\sqrt{T_0}}\, u\,, \\[1.1em]
\ds\cA^\star \,u  \,=\, -\sqrt{T_0}\,\partial_{x}  u \,-\,
\frac{E}{2\sqrt{T_0}}\, u\,.
\end{array}\right.
\end{equation*}
In this framework, the equilibrium $D_\infty$ to  \eqref{Hermite:D} is given by
\begin{equation}
  \label{Dinf}
D_{\infty,k} \,=\,
\left\{
  \begin{array}{l}
    \sqrt\rho_\infty, \,\, {\rm if } \,\, k=0\,,
    \\[0.9em]
    0, \, \, {\rm else\,,}
    \end{array}
  \right.
\end{equation}
and estimate \eqref{esti:L2} simply rewrites
\be
\label{estim:L2}
\frac{1}{2}\,\frac{\dD}{\dD t}\,\| D(t) - D_\infty \|^2_{L^2}
\,+\,
\frac{1}{\tau(\eps)}\,\sum_{k\in\N^\star} k\,\left\|D_k(t)\right\|^2_{L^2\left(\T\right)}
\,=\, 0\,,
\ee
where $\|\cdot\|_{L^2}$ stands for the overall $L^2$-norm \textbf{with no weight}
\[
\|D\|^2_{L^2}
\,=\,
\sum_{
	k\in\N
}
\| D_{k}\|_{L^2\left(\T\right)}^2\,.
\]
On top of that, the limit of the diffusive regime is given by
$D_{\tau_0}\,=\,\left(D_{\tau_0,k}\right)_{k\in \N}$ defined as follows
\be
\label{eq:lim1}
D_{\tau_0,k} = \left\{
\begin{array}{l}
	D_{\tau_0,0}, \,\,{\rm if}\,\,k=0\,,
	\\[1.em]
	0, \,\, {\rm else\,,} 
\end{array}\right.
\ee
where the first Hermite coefficient $D_{\tau_0,0}$ solves the following drift-diffusion equation
\be
\label{eq:lim0}
\partial_t  D_{\tau_0,0} \,+\, \tau_0\,\cA^\star \cA D_{\tau_0,0}   \, =\, 0\,,
\ee
which is obtained substituting $\rho_{\tau_0}$ with $D_{\tau_0,0}\,\sqrt{\rho}_{\infty}$ in equation \eqref{drift:diffusion}.\\

To conclude this section, we introduce some additional norms which arise naturally along our analysis. In Section \ref{sec:23}, we consider the following $H^{-1}$ norm defined on the $L^2$ subspace orthogonal to $\sqrt{\rho}_{\infty}$: for all $g \in L^2 \left(\T\right)$ which meets the condition
\begin{equation}\label{compatibility:g}
	\int_\T g\,\sqrt\rho_\infty\,\dD x\,=\,0\,,
\end{equation} 
we set 
\[
\left\|g\right\|_{H^{-1}}\,=\,\left\|\cA\,u\right\|_{L^2\left(\T\right)}\,,
\]
where $u$ solves the following elliptic equation
\be
\label{eq:elliptic}
\left\{
\begin{array}{l}
	\ds \cA^\star\cA\,u\ =\ g\,,
	\\[1.1em]
	\ds\int_\T u\,\sqrt\rho_\infty\,\dD x\,=\,0\,.
\end{array}\right.
\ee
	The latter equation admits a unique solution in $H^2\left(\T\right)$ for any data $\ds g\,\in\,L^2\left(\T\right)$ that meets the compatibility condition \eqref{compatibility:g}. This well-posedness result crucially relies on the Poincar\'e inequality \eqref{ineg}.

In Section \ref{sec:23}, we use the following $H^1$ norm, defined for all $\ds D\,=\,\left(D_k\right)_{k \in \N}$ as follows
\[
\|\cB\,D\|_{L^2}^2
\,=\,
\sum_{k\in\N} \|\cB_{k}\,D_{k}\|_{L^2\left(\T\right)}^2\,,
\]
where the family of differential operator 
$
\cB\,=\,\left(\cB_k\right)_{k\,\geq\,0}$ is defined as follows
\be
\label{def:B}
\cB_k = \left\{
\begin{array}{l}
	\cA\,,   \,{\rm if} \, \,k\,=\, 0\,,
	\\[0.9em]
	\cA^\star, \,\,{\rm else\,.}
\end{array}\right.
\ee
To end with, we introduce the notation $D_{\perp}\,=\,(D_{\perp,k})_{k
  \in \N}$, which corresponds to the Hermite coefficients of $f\,-\,\rho\,\cM$, that is

\be
\label{def:Dperp}
D_{\perp,k} = \left\{
\begin{array}{l}
	\,0, \,\,{\rm if}\,\,k=0\,,
	\\[1.em]
	\,D_k, \,\, {\rm else\,,} 
\end{array}\right.
\ee
so that
$$
\| D_\perp \|_{L^2} \,=\,  \|f\,-\,\rho\,\cM\|_{L^2(f_\infty^{-1})}.
$$
      \subsection{Main results}
      \label{sec:21}
In this section, we present two results which aim at describing the dynamics of \eqref{vlasov} in various regimes ranging from long time behavior to diffusive limit. We aim for result which capture simultaneously the limits $t\rightarrow +\infty$ and $\eps \rightarrow 0$, in order to lay the groundworks for our upcoming numerical analysis, in which we will build a scheme robust enough so that it captures all these situations.\\
Our first main result tackles the long time behavior of the solution $\ds D=\left(D_k\right)_{k \in \N}$ to \eqref{vlasov}. It is uniform with respect $\eps$ and covers all the regimes of interests since we only impose assumption \eqref{cond:tau} on the scaling parameter $\tau(\eps)$. This result is the first step towards its discrete analog, Theorem \ref{th:main1:h}
\begin{theorem}
	\label{th:main1}
	Suppose that condition \eqref{cond:tau}  on $\tau(\eps)$ is satisfied and let $D\,=\,(
	D_{k} )_{k\in\N}$ be the solution to \eqref{Hermite:D}
	with an initial datum $D^{0,\eps}$. There exists some positive constant $C$ depending only on $\Phi$ and $T_0$ such that 
\begin{itemize}
	\item[$(i)$]  under the condition $\ds\left\|D(0)\right\|_{L^2}\,<\,+\infty$, it holds for all times $t\,\geq\,0$
$$
	\left\|D(t)\,-\,D_\infty\right\|_{L^2} \,\leq \,
	\sqrt{3}\, \left\|D(0)\,-\,D_\infty\right\|_{L^2} \,\exp{
	\left(
	-
	\frac{\tau(\eps)}{\eps^2}\,
	\kappa\, t
	\right)
}\,;
$$
\item[$(ii)$]  under the condition $\ds\left\|\cB\,D(0)\right\|_{L^2}\,+\,\left\|D(0)\right\|_{L^2}\,<\,+\infty$, it holds for all times $t\,\geq\,0$
	$$
	\left\|\cB D(t)\right\|_{L^2} \,\leq \, \sqrt{3}\,\left( C
	\left(
	\overline{\tau}_0+1
	\right)\left\|D(0)\,-\,D_\infty\right\|_{L^2}
	  \,+\, \left\|\cB D(0)\right\|_{L^2} \right)\exp{
		\left(
		-
		\frac{\tau(\eps)}{\eps^2}\,
		\kappa\, t
		\right)
	}\,;
	$$
\end{itemize}
	where $\kappa >0$ is given by
	$$
	\kappa \,=\, \frac{1}{
		C\,(\overline{\tau}_0^2+1)
	} \,.
	$$
\end{theorem}
The proof  of this result is provided in  Section \ref{sec:23}. The main difficulty here consists in proving the convergence of the
first coefficient $D_0$ in the Hermite decomposition of $f$
towards the equilibrium $\sqrt{\rho}_{\infty}$. We adapt
hypocoercivity methods developed in \cite{Villani:AMS,DMS} to the
framework of Hermite decomposition. Instead of estimating directly the
quantities of interest, we introduce modified entropy functionnals
(see \eqref{eq:H0} and \eqref{eq:H1}), in order to recover dissipation
and thus a convergence rate on $D_0$. Then, the second item
tackles the convergence in a $H^1$ setting. Though a bit more
technical, this second convergence result contains no main additional
difficulty in comparison to the $L^2$ convergence result. Actually
this latter result is essentially motivated by the analysis of the
regime $\eps\rightarrow 0$ presented below. 
\\

This leads us to our second main result, which describes the behavior
of the system as $\eps$ vanishes. We distinguish the diffusive regime, which corresponds
to the case where $\tau(\eps)$ satisfies \eqref{cond2a:tau}  and the intermediate situations between long time and
diffusive regime where $\tau(\eps)$ satisfies \eqref{cond2b:tau}.  We will adapt this result into the fully discrete setting in Theorem \ref{th:main2:h}

\begin{theorem}
  \label{th:main2}
	Suppose that $\tau(\eps)$ meets assumption \eqref{cond:tau}. For all positive $\eps$, consider $D=(
	D_{k} )_{k\in\N}$ the solution to \eqref{Hermite:D}
	with an initial datum $D(0)$ such that
	\[
	\left\|D(0)\right\|_{H^{1}}^2\,:=\,
	\left\|\cB D(0)\right\|_{L^2}^2\,+\,
	\left\|D(0)\right\|_{L^2}^2\,<\,+\infty\,.
	\]
	The following statements hold true uniformly with respect to $\eps$
	\begin{itemize}
\item[$(i)$]  suppose that $\tau(\eps)$ satisfies
          \eqref{cond2a:tau}, that is $\tau(\eps)\sim \tau_0\,\eps^2$  and for simplicity, suppose
\begin{equation}\label{cond3:tau}
	\left|\frac{\tau(\eps)}{\tau_0\,\eps^2}\,-\,1\right|\,\leq\,\frac{1}{2}\,,\quad \forall\,\eps\,>\,0\,
\end{equation}
and consider $D_{\tau_0}=(D_{\tau_0,k})_{k\in\N}$
given by \eqref{eq:lim1}. On the one hand, it holds for all time $t\in\R^+$
 \begin{equation*}
\left\|D_{\perp}(t)\right\|_{L^2} \,\leq \,
\left\|D_{\perp}(0)\right\|_{L^2}\,e^{-t/(4\tau_0 \eps^2)}\,+\,\tau_0\,\eps\, C(\overline{\tau}_0+1)\,\left\|D(0)-D_\infty\right\|_{H^1}\,
e^{-
	\tau_0\,
	\kappa\, t
}\,,
\end{equation*}
where $D_{\perp}$ is given in \eqref{def:Dperp}; on the other hand, it holds
\begin{align*}
\left\| D_0(t)-D_{\tau_0,0}(t)\right\|_{H^{-1}} \,\leq \, 
&C\left(
\left\|D_0(0)-D_{\tau_0,0}(0)\right\|_{H^{-1}}
\,+\,\eps\,\tau_0\,(\overline{\tau}_0^3+1)\left\|D(0)-D_{\infty}\right\|_{H^{1}}\right)
e^{
	-
	\tau_0\,
	\kappa\, t}\\
+\,
&C\left|
\frac{\tau_0\eps^2}{\tau(\eps)}
-
1
\right|\left\|D_{\tau_0}(0)-D_{\infty}\right\|_{L^2}\,
e^{
	-
	\tau_0
	\kappa\, t}
\,;
\end{align*}
\item[$(ii)$] suppose that $\tau(\eps)$ satisfies
\eqref{cond2b:tau}, that is $\tau(\eps)/\eps^{2} \rightarrow +\infty  $ as $\eps$ vanishes. Then it holds for all time $t\in\R^+$
\begin{equation*}
	\left\|D_{\perp}(t)\right\|_{L^2} \,\leq \,
	\left\|D_{\perp}(0)\right\|_{L^2}\,e^{-t/(2\tau(\eps))}\,+\,\frac{\tau(\eps)}{\eps} C(\overline{\tau}_0+1)\,\left\|D(0)-D_\infty\right\|_{H^1}
	\,
	e^{-
		\frac{\tau(\eps)}{\eps^2}\,
		\kappa\, t
	}\,,
\end{equation*}
as well as
\begin{equation*}
\left\| D_0(t)-D_{\infty,0}\right\|_{H^{-1}} \,\leq \, 
C\left(
\left\|D_0(0) -D_{\infty,0}\right\|_{H^{-1}}
\,+\,\frac{\tau(\eps)}{\eps}\,(\overline{\tau}_0^3+1)\left\|D(0)-D_{\infty}\right\|_{H^{1}}\right)\,
e^{
-	
\frac{\tau(\eps)}{\eps^2}\,\kappa\, t}
\,.
\end{equation*}
\end{itemize}
In the latter estimate, constant $C$ only depends on $\Phi$ and $T_0$ and exponent $\kappa$ is given by
$$
\kappa \,=\, \frac{1}{
C\,(\overline{\tau}_0^2+1)
} \,.
$$
\end{theorem}
The proof  of this result is provided in  Section \ref{sec:24}, it
showcases two major difficulties. The first one is similar to the one
encountered in Theorem \ref{th:main1}; instead of estimating directly
the $H^{-1}$ norm between the first Hermite coefficient $D_0$ and its limit, we find the
right intermediate quantity in order to recover dissipation (see \eqref{eq:E2}). However, unlike in the case of Theorem \ref{th:main1},
we crucially need to incorporate derivatives of the solution $D$ to
\eqref{vlasov} in this quantity in order to obtain some convergence
rates. This leads us to the second
difficulty, which is that we propagate some regularity. Furthermore,
since Theorem \ref{th:main2} describes simultaneously the large time
behavior and the asymptotic $\eps \rightarrow 0$, it is not sufficient
to propagate derivative globally nor uniformly with respect to time,
we need instead to prove a convergence result in regular norms. This
motivates item $(ii)$ in Theorem \ref{th:main1}, which will play a key role in our proof. This regularity issue explains why we prove $H^{-1}$ convergence with respect to the first Hermite coefficient whereas we achieve strong $L^2$ convergence with respect to other coefficients. To be noted that strong $L^2$ convergence for the first coefficient may be achieved with our method at the price of loosing pointwise estimate with respect to time and thus considering integrated norms with respect to the time variable.\\

Theorems \ref{th:main1} and \ref{th:main2} fully answer their purpose,
which is to describe the dynamics of \eqref{vlasov} in the regime of
interests, uniformly with respect to all parameters at play here.

\subsection{Preliminary results}
\label{sec:22}
Let us first emphasize the important properties satisfied by $\cA$, which we will need to recover later on, in the discrete setting. First, $\cA^\star$ is its dual operator in $L^2(\T)$, indeed for all $u$, $v\in H^1(\T)$ it holds
 \be 
 \label{prop1:A}
\left\langle\cA^\star u ,\, v\right\rangle = \left\langle\cA v,\, u\right\rangle,
\ee
 where $\langle.,\,.\rangle$ denotes the classical scalar product in
 $L^2(\T)$.
Furthermore, we have $D_{\infty,0}$ lies in the kernel of $\cA$, indeed
\be
\label{prop2:A}
\cA \,D_{\infty,0} \,=\,  0\,;
\ee
in this setting, conservation of mass is ensured by the following property
\begin{equation}\label{cons:mass1}
	\int_{\T}\cA^{\star}\,u\,\sqrt{\rho}_{\infty}\,\dD x\,=\,0\,,
\end{equation}
indeed, considering equation \eqref{Hermite:D} with index $k=0$ integrated over $\T$ and applying the latter relation with $u\,=\,D_1$, we obtain 
\[
\frac{\dD}{\dD t}\int_{\T}D_0(t)\,\sqrt{\rho}_{\infty}\,\dD x\,=\,0\,,
\]
and therefore
\begin{equation}\label{cons:mass}
\int_{\T}D_0(t)\,\sqrt{\rho}_{\infty}\,\dD x\,=\,\int_{\T}D_{\infty,0}\,\sqrt{\rho}_{\infty}\,\dD x\,;
\end{equation}
we also point out that since 
\[
\sqrt{T_0}\,(\cA+\cA^\star)\,=\,\partial_x\Phi\,,
\]
it holds
\be
\label{prop:A+Astar}
\| \, \left(\cA + \cA^\star\right) u\|_{L^2} \,\leq \, \frac{1}{\sqrt{T_0}}\,\|\Phi\|_{W^{1,\infty}} \|u\|_{L^2}\,,
\ee
on top of that, operators $\cA$ and $\cA^\star$ do not commute and we have
$$
 [\cA, \, \cA^\star] \,=\,  \cA \, \cA^\star - \cA^\star \, \cA \,=\, \partial_{xx}\Phi\,,
 $$
which yields
\be
\label{prop3:A}
\| \, [\cA, \, \cA^\star] \,u\|_{L^2} \,\leq \, \|\Phi\|_{W^{2,\infty}}\, \|u\|_{L^2}\,;
\ee
the last key property verified by operator $\cA$ is the following Poincar\'e-Wirtinger inequality: under the compatibility condition \eqref{compatibility:g} on $u \in H^{1}\left(\T\right)$ it holds 
\be
\label{ineg}
\|u\|_{L^2}
\,\leq\,C_P\sqrt{T_0}\,\left(\int_{\T}\left|\partial_x\left(\frac{u}{\sqrt\rho_\infty}\right)\right|^2\,\rho_\infty\,\dD
x\right)^{1/2} \,=\,C_P\| \cA\,u\|_{L^2}\,,
\ee
for some positive constant $C_P$ depending only on the potential
$\Phi$ and $T_0$. A proof of this result will be given in the discrete setting (see Lemma \ref{lem:3.1}), we do not detail it in the continuous case since it is not our main interest here.\\

\subsection{Proof of Theorem \ref{th:main1}}
\label{sec:23}
It is worth to mention that estimate \eqref{estim:L2} itself is not sufficient to conclude
on the rate of convergence of $D$ to the equilibrium $D_\infty$, since there is no dissipation with respect to the zero-th Hermite coefficient $D_{0}$. Therefore, it does not provide quantitative estimates when it comes to its convergence towards $D_{\infty,0}$. Recovering this dissipation is
the key feature of hypocoercivity \cite{Villani:AMS,DMS}. In our setting it is done by
combining the equations on $D_0$ and $D_{1}$,  to remove stiff terms
\be
\label{eq:D0}
\partial_t \left(  D_{0} + \frac{\tau(\eps)}{\eps} \cA^\star
D_{1}  \right) \,+\, \frac{\tau(\eps)}{\eps^2}  \,\left( \cA^\star \cA\,D_{0}  \,- \, \sqrt{2}\,
\left(\cA^\star\right)^2 D_{2}\right) 
\, =\, 0\,.
\ee
To  prove quantitative estimates on the solution to
\eqref{Hermite:D}, we therefore introduce the ''modified entropy
functional''  \cite{DMS, Villani:AMS}: for any $\alpha_0>0$, which will be specified later, we define $\cH_0$ as 
\begin{equation}
\cH_0[D| D_\infty] \,=\,
\frac12\,\|D(t)-D_\infty\|_{L^2}^2\,+\, \alpha_0\,\left\langle
  \frac{\tau(\eps)}{\eps}\,\cA^\star D_1, u^\eps\right\rangle\,,
 \label{eq:H0}
\end{equation}
where $u^\eps$ is the particular solution to equation \eqref{eq:elliptic} with source term is $g\,=\,D_0\,-\,D_{\infty,0}$. To be noted that $g\,=\,D_0\,-\,D_{\infty,0}$ fullfils the compatibility condition \eqref{compatibility:g},
thanks to the conservation of mass property \eqref{cons:mass1}.

The first step consists in proving some intermediate results on the solutions $u^\eps$ to \eqref{eq:elliptic} 
 
\bl
\label{lem:1}
 Consider any  $g \in L^2(\T)$ which meets condition \eqref{compatibility:g} and  $u$ the
 corresponding solution to \eqref{eq:elliptic}. Then,
 $u$ satisfies the following estimate
 \be
 \label{lem:110}
    \|\cA\, u \|_{L^2} \,\leq\, C_P\,\|g\|_{L^2}\,,
    \ee
    and
    \be
    \label{lem:111}
  \|\cA^2\, u\|_{L^2}  \,\leq\, \left(1\,+\,\frac{C_P}{\sqrt{T_0}}\,\|\Phi\|_{W^{1,\infty}}\right)\|g\|_{L^2}\,,
\ee
where $C_P$ is the Poincar\'e constant in \eqref{ineg}.
\\
Moreover, considering now the solution $D$ to
\eqref{Hermite:D} and $u^\eps$ the solution to \eqref{eq:elliptic} with source term $g\,=\,D_0\,-\,D_{\infty,0}$, it holds for all time $t\,\geq\,0$
\be
\label{lem:12}
\eps\,\|\cA \,\partial_t u^\eps(t)\|_{L^2}   \,\leq\, \|D_1(t) \|_{L^2}\,.
\ee
\el
\begin{proof}
 The first estimate is obtained by testing the elliptic equation \eqref{eq:elliptic} against $u$ and applying \eqref{prop1:A} 
 
 \[
  \|\cA\, u \|_{L^2}^2 \,\leq\,\| g\|_{L^2}\,\|u\|_{L^2}\,,
\]
hence the Wirtinger-Poincar\'e inequality \eqref{ineg} yields,
$$
\|\cA \,u\|_{L^2} \,\leq\,C_P\,\| g\|_{L^2}\,.
  $$
  For the second estimate, we rewrite $\cA^2 u$ as follows
  
  $$
  \cA^2 \,u \,=\, -\cA^\star\cA \,u \,+
  \left(\cA+\cA^\star\right) \cA \,u\,,  
  $$
  then we replace $\cA^\star\cA \,u $ according to equation \eqref{eq:elliptic}, take the $L^2$ norm on both sides of the relation and apply in turn \eqref{prop:A+Astar} to estimate operator $\ds\cA+\cA^\star$ and item \eqref{lem:110} to estimate the norm of $\cA\,u$, it yields
$$\|\cA^2\,u \|_{L^2} \,\leq \, \left( 1 \,+\, \frac{C_P}{\sqrt{T_0}}\,\,\|\Phi\|_{W^{1,\infty}}\right)\|g\|_{L^2}\,.
$$
 For the third estimate we consider now that $D$ is solution to
 \eqref{Hermite:D} and  first take the time derivative of the
 elliptic equation \eqref{eq:elliptic} and use the equation
 \eqref{Hermite:D} on  $D_0$  to get
$$
  \eps\,\partial_t(\cA^\star\cA\, u^\eps)\, =\, \eps\,\partial_t(D_0 -
  D_{\infty,0}) \,=\, -\cA^\star D_1\,.
  $$ 
Then multiply by $\partial_t u^\eps$ and use \eqref{prop1:A} to get 
\[
 \|\partial_t\cA \,u^\eps\|_{L^2}^2\, =\, -\frac{1}{\eps}\left\langle
   D_1, \partial_t\cA \,u^\eps\right\rangle
 \,\leq\,\frac{1}{\eps}\,\|D_1\|_{L^2}\,\|\partial_t\cA \,u^\eps\|_{L^2}\,.
\]
\end{proof}

Thanks to the latter result  we now prove that for small enough $\alpha_0>0$, the square root of the modified entropy is
equivalent to the $L^2$ norm of  $D-D_\infty$

\bl
Suppose that condition \eqref{cond:tau}  on $\tau(\eps)$ is satisfied.  Then for all
$\alpha_0\in(0,\overline{\alpha}_0)$,  with $\overline{\alpha}_0=1/(4\,\overline{\tau}_0\,C_P)$ 
and $D\in L^2(\T)$ such that $D_0-D_\infty$ satisfies the
compatibility condition \eqref{compatibility:g}, one has
 \begin{equation}
   \label{equiv}
 \|D -D_\infty\|_{L^2}^2 \leq\,4\,\cH_0[D|D_\infty]\,\leq \,3\,\|D -D_\infty\|_{L^2}^2\,.
 \end{equation}
 \label{lem:2}
\el
\begin{proof}
We estimate the additional term in the expression of $\cH_0$ by applying the duality formula \eqref{prop1:A} and then Cauchy-Schwarz inequality
\[
|\left\langle \cA^\star D_1, u^\eps \right\rangle |\,=\,
|\left\langle D_1, \cA\,u^\eps
\right\rangle_{L^2}|\,\leq\,\|D_1\|_{L^2}\,\|\cA\,u^\eps\|_{L^2}\,.
\]
Then, we apply item \eqref{lem:110} of Lemma \ref{lem:1} with $u^\eps$ and $g\,=\,D_0\,-\,D_{\infty,0}$ and upper bound the norm of $\cA\,u^\eps$ accordingly
 \[
  \|D_1\|_{L^2}\,\|\cA\,u^\eps\|_{L^2}\,\leq\,C_P\,\|D-D_{\infty}\|^2_{L^2}\,,
\]
hence, applying assumption \eqref{cond:tau}, we deduce
$$
 \alpha_0\,\frac{\tau(\eps)}{\eps} \,  |\left\langle \cA^\star
   D_1, u^\eps \right\rangle |  \, \leq \, \alpha_0  \, \overline{\tau}_0\,C_P\, \|D-D_{\infty}\|^2_{L^2}\,.
$$
Choosing $\overline{\alpha}_0=1/(4\,\overline{\tau}_0\,C_P)$,  the result follows for $\alpha_0 \in (0,\overline\alpha_0)$.
\end{proof}

Relying on the previous lemmas, we are now able to carry out the proof
of the first item $(i)$ of Theorem \ref{th:main1}.
We compute the time derivative of the modified relative entropy and split into three terms
\[
\frac{\dD}{\dD t}\cH_0[D(t)|D_\infty]\, =\, \cI_1(t) \,+\, \alpha_0\,\cI_2(t) \,+\, \alpha_0\,\cI_3(t)\,,
\]
where the first one corresponds to the dissipation of the $L^2$ norm \eqref{estim:L2},
$$
 \cI_1\,=\, -\frac{1}{\tau(\eps)}\,\sum_{k\in\N}k\,\left\|D_k\right\|_{L^2}^2\,,
 $$
 whereas the other ones correspond to the additional term of the
 modified relative entropy,  
 $$
 \left\{\begin{array}{l}
   \ds\cI_2\,:=\,
   -\frac{\tau(\eps)}{\eps^2}\,\left\langle\cA^\star\cA\,\left(D_0-D_{\infty,0}\right) \,-\,
   \sqrt{2}\,(\cA^\star)^2\, D_2 ,\, u^\eps\right\rangle  \,
   -\, \frac{1}{\eps} \left\langle\cA^\star D_1
   ,\, u^\eps\right\rangle\,,
   \\[1.em]
  \ds \cI_3\,:=\, +\frac{\tau(\eps)}{\eps}\,\left\langle\cA^\star
  D_1,\,\partial_t u^\eps\right\rangle\,.
\end{array}\right.
$$
On the one hand, the term $\cI_2$ gives the expected dissipation on $(D_0 -
D_{\infty,0})$ since $u^\eps$ solves \eqref{eq:elliptic} with source term $(D_0 -
D_{\infty,0})$. On the other hand we get
some additional terms which can be estimated thanks to \eqref{lem:110}
and \eqref{lem:111} in Lemma
\ref{lem:1}, it yields,
\begin{eqnarray*}
\cI_2 &\leq& -\frac{\tau(\eps)}{\eps^2}\, \| D_0-D_{\infty,0}
\|_{L^2}^2 \,+\, \frac{\tau(\eps)}{\eps^2} \,  \sqrt 2\,\left( 1
  \,+\,\frac{C_P}{\sqrt{T_0}}\,\|\Phi\|_{W^{1,\infty}}\right)\,\|D_0
                  -D_{\infty,0}\|_{L^2}\| D_2\|_{L^2} \\[1.0em]
  && +\, \frac{C_P}{\eps}\, \| D_0 -
     D_{\infty,0}\|_{L^2}\, \| D_1\|_{L^2}\,,
     \\[1.0em]
&\leq& -\frac{\tau(\eps)}{\eps^2}\,\left(1\,-\,
       C\,\eta\right) \, \| D_0-D_{\infty,0}
\|_{L^2}^2 \,+\,  \frac{C}{2\,\eta}\, \left( \frac{\tau(\eps)}{\eps^2}\| D_2\|_{L^2} ^2 \,+\,  \frac{1}{\tau(\eps)}\, \| D_1\|_{L^2}^2\right)\,,
\end{eqnarray*}
for any positive $\eta$ and for some positive constant $C$ depending only on $T_0$ and $\Phi$. The term $\cI_3$ is estimated directly by
applying \eqref{lem:12} of Lemma \ref{lem:1},
$$
\cI_3 \,\leq\, \frac{\tau(\eps)}{\eps^2}\|D_1\|_{L^2}^2\,.
$$
From these latter estimates and taking $\eta=1/(2C)$, we get the following inequality
\begin{eqnarray*}
  &&\frac{\dD}{\dD t}\cH_0[D|D_\infty]
  \\
  &&\leq\, -
\frac{\tau(\eps)}{\eps^2}\,
\left(
\frac{\alpha_0}{2}\,
\| D_0-D_{\infty,0}
\|_{L^2}^2 +\left(\frac{\eps^2}{\tau(\eps)^2}\, -\,  C^2\,
\left(1\,+\frac{\eps^2}{\tau(\eps)^2}\right)\, \alpha_0\right) \,\sum_{k\in\N}k\,\left\|D_k\right\|_{L^2}^2
\right)\,.
  \end{eqnarray*}
Under the following condition
\[
\alpha_0
\,\leq\,
\textrm{argmax}_{\alpha>0}\,
\min\left(
\frac{\alpha}{2},\,
\frac{\eps^2}{\tau(\eps)^2}\, -\,  C^2\,
\left(1\,+\frac{\eps^2}{\tau(\eps)^2}\right)\, \alpha\right)\,,
\]
which, according to assumption \eqref{cond:tau} on $\tau(\eps)$, is fulfilled as long as 
\[
\alpha_0
\,\leq\,
\frac{1}{
C\,(\overline{\tau}_0^2+1)
}\,,
\]
for some constant $C$ depending only on $\Phi$ and $T_0$, and taking $\kappa_0$ such that $3\,\kappa_0/4\,=\,\alpha_0/2$,
we derive the following estimate
\[
 \frac{\dD}{\dD t}\cH_0[D|D_\infty]  +
 \frac{\tau(\eps)}{\eps^2}\,
 \frac{3\,\kappa_0}{4}\, \|D - D_\infty \|_{L^2}^2 \,\leq\,0\,.
\]
Then applying Lemma~\ref{lem:2} and taking $\alpha_0\leq\overline{\alpha}_0$, we deduce
\[
\frac{\dD}{\dD t}\cH_0[D|D_\infty] \,+\,
\frac{\tau(\eps)}{\eps^2}\,\kappa_0\,\cH_0[D |D_\infty] \,\leq\,0\,,
\]
which yields after applying Gronwall's lemma, for any $t\geq 0$,
\begin{equation*}
\cH_0[D(t)|D_\infty] \,\leq \, \cH_0[D(0)|D_\infty] \,\exp{
	\left(
	-
	\frac{\tau(\eps)}{\eps^2}\,
	\kappa_0\, t
	\right)
}\,.
\end{equation*}
We conclude this proof by applying Lemma \ref{lem:2} in order to substitute $\cH_0$ with the $L^2$ norm of $D\,-\,D_{\infty}$ in the latter estimate.
\\

We now turn to the proof of the second item $(ii)$ of Theorem
\ref{th:main1}. To estimate the norm of $\cB\,D$, we apply the operator $\cB_k$ to \eqref{Hermite:D}  and next  multiply by $\cB_k D_k$, integrate with
respect to $x\in\T$ and sum over $k\in\N$, it yields
\begin{equation*}
\frac{1}{2}\,\frac{\dD}{\dD t}\,\| \cB D(t) \|^2_{L^2}
  \,=\,\cJ_1(t)\,,
\end{equation*}
where $\cJ_1$ is defined as follows
\begin{equation*}
\cJ_1\,=\,\sum_{k\in\N^\star} -
\frac{k}{\tau(\eps)}\left\|\cB_k D_k\right\|^2_{L^2}\,+\, \frac{\sqrt{k}}{\eps}
     \left( \left\langle\cB_{k-1} \cA^\star
     D_{k},\,\cB_{k-1} D_{k-1}\right\rangle \,- \,  
  \left\langle\cB_k\cA\,D_{k-1},\, \cB_k
     D_{k}\right\rangle\right)\,,
\end{equation*}
where we use that $\cA\,D_{\infty,0}=0$ and $D_{\infty,k}=0$ for
$k>0$. Hence applying an integration by part and from the specific choice \eqref{def:B}
of $\cB$, we have
\begin{equation}\label{def:Jeps1}
\cJ_1
\,=\,-\frac{1}{\tau(\eps)}\,\sum_{k\in\N^\star} k\,\left\|\cB_k D_k\right\|^2_{L^2}
  \,-\, \frac{1}{\eps}\,\sum_{k\geq 2} \sqrt{k}\,
   \left\langle\, [\cA^\star, \cA]\,D_{k-1},\, \cA^\star D_{k}\right\rangle\,.
\end{equation}
Applying Young inequality and property \eqref{prop3:A} on the commutator $ [\cA^\star, \cA]$, we get that
 \begin{equation*}
 \cJ_1
\,\leq\,\frac{1}{\tau(\eps)}
\left(\frac{\eta}{2} \,\|\Phi\|_{W^{2,\infty}}^2 
-1
\right)
\sum_{k\in\N^\star} k\,\left\|\cB_k  D_k\right\|^2_{L^2}
\, +\,  \frac{1}{2 \,\eta}\,\frac{\tau(\eps)}{\eps^2}\,\sum_{k\geq 1} \left\| D_{k}\right\|^2_{L^2}\,.
 \end{equation*}
Therefore, choosing $\eta \leq 1/\|\Phi\|_{W^{2,\infty}}^2$, it yields
\begin{equation}
   \label{estim:H1}
\frac{1}{2}\,\frac{\dD}{\dD t}\,\| \cB D \|^2_{L^2}\,+\,\frac{1}{2\,\tau(\eps)}\,\sum_{k\in\N^\star} k\,\left\|\cB_k D_k\right\|^2_{L^2}\,\leq\, C\,\frac{\tau(\eps)}{\eps^2}\,\sum_{k\geq 1} \left\| D_{k}\right\|^2_{L^2}\,.
\end{equation}
Again since there is no dissipation on the zero-th Hermite coefficient
of $\cB_0\,D_0$, we proceed as for the $L^2$ estimate and  introduce a correction $\cH_1$ given by
\begin{equation}
\cH_1[D|D_\infty] \,=\,
\frac12\,\|\cB D \|_{L^2}^2\,+\, \alpha_1\,\left\langle
 \frac{\tau(\eps)}{\eps}\,\cA\,D_0, D_1\right\rangle\,,
 \label{eq:H1}
\end{equation}
where $\alpha_1$ has to be determined. First, we point out that for small enough $\alpha_1>0$, the modified entropy $\cH_1$ is
controlled by the squares of the $L^2$ norms of  $D-D_\infty$ and $\cB D$.

\bl
Suppose that condition \eqref{cond:tau}  on $\tau(\eps)$ is satisfied.  Then for all
$\alpha_1\in(0,\overline{\alpha}_1)$,  with $\overline{\alpha}_1=1/(2\,\overline{\tau}_0)$
and $D\in L^2(\T)$, one has
\begin{equation}
\label{equiv2}\|\cB D \|_{L^2}^2
-
\|D -D_\infty\|_{L^2}^2
\leq\,4\,\cH_1[D|D_\infty]\,\leq \,3\,\|\cB D \|_{L^2}^2
+\|D -D_\infty\|_{L^2}^2\,.
\end{equation}
\label{lem:3}
\el 
\begin{proof}
The result is obtained applying the Young inequality to the additional term in the definition \eqref{eq:H1} of $\cH_1$
\end{proof}

To complete the proof of the second item $(ii)$ in Theorem \ref{th:main1}, we compute the time derivative of the modified relative entropy and split into two terms
\[
\frac{\dD}{\dD t}\cH_1[D|D_\infty]\, =\, \cJ_1 + \alpha_1\,\cJ_2\,,
\]
where the first one corresponds to the dissipation of the $L^2$ norm
of $\cB\,(D-D_\infty)$ for which we already have
an estimate  \eqref{estim:H1},  that is,
$$
 \cJ_1 \,\leq\,
 -\frac{1}{2\,\tau(\eps)}\,\sum_{k\in\N^\star}k\,\left\|\cB_k D_k\right\|_{L^2}^2\,+\, C\, \frac{\tau(\eps)}{\eps^2}\, \sum_{k\geq 1} \left\| D_{k}\right\|^2_{L^2},
 $$
 whereas the other ones correspond to the additional term of the
 modified relative entropy,
 $$
 \cJ_2\,:=\,
   \frac{\tau(\eps)}{\eps^2}\,\left( \left\langle\cA\cA^\star D_1 ,\, D_1\right\rangle  \,-\,
   \left\|\cA\,D_0\right\|_{L^2}^2  
   \,+\,\sqrt{2} \left\langle\cA\,D_0,\,\cA^\star
     D_2\right\rangle\right)  \,-\, \frac{1}{\eps} \left\langle D_1
   ,\,\cA\,D_0\right\rangle\,.
 $$
From \eqref{prop1:A} and \eqref{prop2:A}
on the operators $(\cA,\,\cA^\star)$,  we have
$$
\frac{1}{\eps}\,\langle D_1, \,\cA\,D_0\rangle \,=\, \left\langle
\frac{1}{\tau(\eps)^{{1}/{2}}}\,\cA^\star D_1,\,
\frac{\tau(\eps)^{{1}/{2}}}{\eps} \left(D_0-D_{\infty,0}\right)\right\rangle\,,
$$
hence  applying twice the Young inequality on the third term of the right
hand side  and on the latter term, it yields
$$
\cJ_2 \,\leq\, -\frac{\tau(\eps)}{\eps^2}\,\left[ \frac{1}{2}\,\|\cA\,D_0\|^2_{L^2} \,-\,
                    \left( 1 + \frac{\eps^2}{\tau(\eps)^2}\right)\,\sum_{k\in\N^\star}
                  k\,\left\|\cB_k D_k\right\|^2_{L^2} \,-\, \|D_0-D_{\infty,0} \|_{L^2}^2\right]\,.
                  $$
Therefore, from these estimates, we get the following inequality
\begin{eqnarray*}                    
  \frac{\dD}{\dD t}\cH_1[D|D_\infty] &\leq& (C+\alpha_1)\,\frac{\tau(\eps)}{\eps^2}\, \|D_0-D_{\infty,0} \|_{L^2}^2
  \\
  &-&\frac{\tau(\eps)}{2\,\eps^2} \left[ {\alpha_1}\,\|\cA\,D_0\|^2_{L^2} \,+\,
                    \left(\frac{\eps^2}{\tau(\eps)^2}-2\,\alpha_1\,\left( 1 + \frac{\eps^2}{\tau(\eps)^2}\right)\right)\,\sum_{k\in\N^\star}
                  k\,\left\|\cB_k D_k\right\|^2_{L^2} \right]\,, 
\end{eqnarray*}
hence choosing  $\alpha_1$
\[
\alpha_1
\,\leq\,
\textrm{argmax}_{\alpha>0}\,
\min\left(\alpha,\,
\frac{\eps^2}{\tau(\eps)^2}\, -\,
2\,\alpha\left(1\,+\frac{\eps^2}{\tau(\eps)^2}\right)\,\right)
\,=\,
\frac{1}{2+3\,\frac{\tau(\eps)^2}{\eps^2}}\,,
\]
which is verified under the following condition
\[
\alpha_1\,\leq\,
\frac{1}{2+3\,\overline{\tau}_0^2}\,,
\]
we get that
         \[
\frac{\dD}{\dD t}\cH_1[D|D_\infty] \,+\, \frac{\tau(\eps)}{\eps^2}\,\frac{\alpha_1}{2}
\left\|\cB
  D \right\|^2_{L^2}   \, \leq\,  C\,\frac{\tau(\eps)}{\eps^2}\,
\|D - D_{\infty}\|_{L^2}^2\,.
\]         
Furthermore, taking $\alpha_1\leq1/(2\,\overline{\tau}_0)$ and applying Lemma \ref{lem:3}, we obtain
$$
\frac{\dD}{\dD t}\cH_1[D|D_\infty] \,+\, \frac{\tau(\eps)}{\eps^2}\,
\frac{2\,\alpha_1}{3}\,
\cH_1[D|D_\infty]  \, \leq\,
C\,\frac{\tau(\eps)}{\eps^2}\, 
\|D - D_{\infty}\|_{L^2}^2\,.        
$$
Then we set
$$
\kappa_1\,=\,\min\left(\frac{2\,\alpha_1}{3}\,,\,\kappa_0\right)
$$
and multiply the latter inequality by $\exp{\left(\frac{\tau(\eps)}{\eps^2}\,
	\frac{2\,\alpha_1}{3}\,t\right)}$, integrate in time and apply the first  item $(i)$ of Theorem \ref{th:main1} to estimate the
right hand side, this yields
$$
\cH_1[D(t)|D_\infty] \,\leq \, \left( C
\left(
\overline{\tau}_0^2+1
\right)
\left\|D(0)\,-\,D_\infty\right\|_{L^2}^2  \,+\,  \cH_1[D(0)|D_\infty] \right)\,\exp{
	\left(
	-
	\frac{\tau(\eps)}{\eps^2}\,
	\kappa_1\, t
	\right)
}\,.
$$
We conclude this proof by substituting $\cH_1$ with the norm of $\ds\cB D$ in the latter estimate according to Lemma \ref{lem:3}.

\subsection{Proof of Theorem \ref{th:main2}}
\label{sec:24}

Once again, instead of estimating directly the $H^{-1}$ norm of $D_0-D_{\tau_0}$, we introduce the following quantity, meant to recover dissipation on the zero-th Hermite coefficient
\begin{equation}\label{eq:E2}
	\cE(t)\,=\,\frac{1}{2}\, \|\cA\,v^\eps(t)\|_{L^2}^2\,,
\end{equation}
where $v^\eps(t)$ solves the elliptic equation \eqref{eq:elliptic} with source term given by 
\[\ds g(t)\,=\,D_0(t)\,
+\,\frac{\tau(\eps)}{\eps}\,\cA^\star D_1(t)
\,-\,D_{\tau_0,0}(t)\,,\]
where $D_0(t)$ and  $D_1(t)$ are the first two components of the solution
$D(t)$ of \eqref{Hermite:D} and $D_{\tau_0,0}(t)$ is either the unique solution to the
convection-diffusion equation \eqref{eq:lim0} when $\tau_0$ is finite
or the stationary solution  $D_{\infty,0}$ given by \eqref{Dinf} when $\tau_0=\infty$. The latter right hand side is motivated by equation \eqref{eq:D0} since it is given by the difference between $D_0\,
+\,\frac{\tau(\eps)}{\eps}\,\cA^\star D_1$ and $D_{\tau_0,0}$. We point out that the latter source term meets the compatibility condition \eqref{compatibility:g} thanks to property \eqref{cons:mass1}, which ensures that $\ds\cA^\star D_1(t)$ is orthogonal to $\ds\sqrt{\rho}_{\infty}$ in $\ds L^2\left(\T\right)$.\\

Before proving the first item of Theorem \ref{th:main2}, let us
present some preliminary results. On the one hand, the following Lemma ensures that $\cE(t)$ is
controlled by the squares of the $L^2$ norm of $\cB D(t)$ and the $H^{-1}$ norm of $D_0(t)-D_{\tau_0,0}(t)$
\bl
We consider $\cE(t)$ defined by \eqref{eq:E2}. It holds uniformly with respect to $\eps$
 \begin{equation}
   \label{equiv3}
   \cE(t)\,\leq\, \|D_0(t) -D_{\tau_0,0}(t)\|_{H^{-1}}^2
 \,+\,C_P^2\,\frac{\tau(\eps)^2}{\eps^2}\,
   \|\cB D(t)\|_{L^2}^2\,,
 \end{equation}
 and
 \begin{equation}
   \label{equiv4}
   \frac{1}{4}\,
   \|D_0(t) -D_{\tau_0,0}(t)\|_{H^{-1}}^2
 \,-\,C_P^2\,\frac{\tau(\eps)^2}{2\,\eps^2}\,
   \|\cB D(t)\|_{L^2}^2
   \,\leq\,
   \cE(t)
   \,.
 \end{equation}
 \label{lem:4}
\el 
\begin{proof}
Defining $w^\eps$ and $u_{\tau_0}$ as the respective solutions to \eqref{eq:elliptic} with source term $g\,=\,\cA^\star D_1$ and $D_{\tau_0,0}\,-\,D_{\infty,0}$, it holds
\[
v^\eps\,=\,u^\eps\,-\,u_{\tau_0}\,+\,\frac{\tau(\eps)}{\eps}\,w^\eps\,.
\]
We apply operator $\cA$ to the latter relation, take the $L^2$ norm, and apply the triangular inequality, it yields 
\[
\sqrt{2\,\cE}\,\leq\,
\left\|
\cA\left(u^\eps-u_{\tau_0}\right)
\right\|_{L^2}
\,+\,
\frac{\tau(\eps)}{\eps}\,
\left\|
\cA\,w^\eps
\right\|_{L^2}\,,
\]
and
\[
\left\|
\cA\left(u^\eps-u_{\tau_0}\right)
\right\|_{L^2}
\,-\,
\frac{\tau(\eps)}{\eps}\,
\left\|
\cA\,w^\eps
\right\|_{L^2}
\,\leq\,\sqrt{2\,\cE}\,.
\]
We estimate 
$
\left\|
\cA\,w^\eps
\right\|_{L^2}
$ applying \eqref{lem:110} in Lemma \ref{lem:1} with source term $g\,=\,\cA^\star D_1$, this yields
\[
\sqrt{2\,\cE}\,\leq\,
\|D_0 -D_{\tau_0,0}\|_{H^{-1}}
\,+\,
\frac{\tau(\eps)}{\eps}\,
C_P\,\|\cB D\|_{L^2}\,,
\]
and
\[
\|D_0 -D_{\tau_0,0}\|_{H^{-1}}
\,-\,
\frac{\tau(\eps)}{\eps}\,
C_P\,\|\cB D\|_{L^2}
\,\leq\,\sqrt{2\,\cE}
\,.
\]
We obtain the result taking the square of the latter inequalities and applying Young's inequality.\\
\end{proof}

On the other hand, when $\tau_0$ is finite, we observe that the long
time behavior of $\ds D_{\tau_0,0}$ may be easily  investigated. Indeed, since $\ds \cA\,D_{\infty,0}\,=\,0$, we have that $D_{\tau_0,0} - D_{\infty,0}$ also solves \eqref{eq:lim0}. Therefore, multiplying \eqref{eq:lim0} by $D_{\tau_0,0} - D_{\infty,0}$, integrating over $\T$ and applying the Poincar\'e inequality \eqref{ineg}, we obtain the following estimate after applying Gronwall lemma
\begin{equation}\label{cv:D0}
	\|D_{\tau_0}(t) -D_{\infty}\|_{L^2}
	\,\leq\,
	\|D_{\tau_0}(t) -D_{\infty}\|_{L^2}\,\exp{
		\left(-\frac{\tau_0}{C_P^2}\,t\right)}\,,\quad\forall\,t\in\R^+\,.
\end{equation}
We are now able to prove the first item  $(i)$ of Theorem
\ref{th:main2}, which treats the case where $\tau(\eps) \sim
\tau_0\,\eps^2$, when $\eps\rightarrow 0$ where $\tau_0 \in
\R^+_\star$.
To derive the first estimate in item $(i)$ of Theorem \ref{th:main2}, our starting point is the $L^2$ estimate \eqref{estim:L2} which ensures
\begin{eqnarray*}
\frac{1}{2}\frac{\dD}{\dD t}\left\|D_{\perp}(t)\right\|^2_{L^2}\, +
\frac{1}{\tau(\eps)}\, \left\|D_{\perp}(t)\right\|^2_{L^2}
&\leq& -\frac{1}{2}\,\frac{\dD}{\dD t}\,\| D_0(t) -D_{\infty,0}\|^2_{L^2} 
\\
&\leq& -\frac{1}{\eps} \left\langle \cA^\star D_1(t),\,
D_0(t)-D_{\infty,0}\right\rangle\,\\ &=&\,-\frac{1}{\eps} \left\langle  D_1(t),\,
\cA\left(D_0(t)-D_{\infty,0}\right)\right\rangle\,, 
\end{eqnarray*}
hence it gives from the Young inequality
$$
\frac{\dD}{\dD t}\left\|D_{\perp}(t)\right\|^2_{L^2}  \, +\, \frac{1}{\tau(\eps)}\,\left\|D_{\perp}(t)\right\|^2_{L^2}
\,\leq\,\frac{\tau(\eps)}{\eps^2}\,\|\cB D(t)\|_{L^2}^2\,. 
$$
We bound $\ds\|\cB D(t)\|_{L^2}^2$ applying item $(ii)$ of Theorem \ref{th:main1}. After multiplying the latter estimate by $\ds e^{t/\tau(\eps)}$ and integrating with respect to time, it yields
\begin{eqnarray*}
\left\|D_{\perp}(t)\right\|^2_{L^2} &\leq&
\ds\left\|D_{\perp}(0)\right\|^2_{L^2}\,\exp\left(-\frac{t}{\tau(\eps)}\right)
  \\
  &+&
\left(C(\overline{\tau}_0^2+1)\,\left\|D(0)-D_\infty\right\|^2_{L^2} + \left\|\cB D(0)\right\|_{L^2}^2 
\right)
\frac{3\,\tau(\eps)^2
}{\eps^2-\kappa\,\tau(\eps)^2}\,
\exp\left(-\frac{\tau(\eps)}{\eps^2}\,
\kappa\, t
\right)\,,
\end{eqnarray*}
where $C$ is a positive constant depending only on $\Phi$ and $T_0$ and $\kappa\,=\,\left(C(\overline{\tau}_0^2+1)\right)^{-1}$. Then we apply condition \eqref{cond:tau} on $\tau(\eps)$, which ensures that taking $C$ greater than $2$ in the definition of $\kappa$, it holds 
$
\ds
1/2\,\leq\, 1-\kappa\,\tau(\eps)^2/\eps^2
$ uniformly with respect to $\eps$.
Therefore, we deduce the following estimate, which yields the first result in $(i)$ of Theorem \eqref{th:main1}, after taking its square root and applying assumption \eqref{cond3:tau} in order to substitute $\tau(\eps)$ with $\tau_0\,\eps^2$
 \begin{equation*}
 \left\|D_{\perp}(t)\right\|^2_{L^2} \,\leq \,
 \left\|D_{\perp}(0)\right\|^2_{L^2}\,e^{-\frac{t}{\tau(\eps)}}\,+\, 6\left(C(\overline{\tau}_0^2+1)\,\left\|D(0)-D_\infty\right\|^2_{L^2} + \left\|\cB D(0)\right\|_{L^2}^2 
 \right)
\,\frac{\tau(\eps)^2}{\eps^2}\,
e^{-
\frac{\tau(\eps)}{\eps^2}\,
\kappa\, t
}\,. 
 \end{equation*}
We now prove the second result in item $(i)$ of Theorem \ref{th:main2}. To do so, we evaluate $\cE$ observing that
 $$
\frac{\dD \cE }{\dD t} \,=\, \left\langle \partial_t\left(D_0\,
  +\,\frac{\tau(\eps)}{\eps}\,\cA^\star D_1
  \,-\,D_{\tau_0,0}\right),\,v^\eps\right\rangle\,.
$$
Therefore, relying on equations \eqref{eq:D0} and \eqref{eq:lim0} we deduce 
\begin{equation*}
\frac{\dD \cE }{\dD t}  \,=\, -\,\frac{\tau(\eps)}{\eps^2}
\, \| D_0
+\frac{\tau(\eps)}{\eps}\,\cA^\star D_1
-D_{\tau_0,0}\|_{L^2}^2
\,+\,
\cE_{1}
\,+\,
\cE_{2}
\,+\,
\cE_{3}\,,
\end{equation*} 
where
\begin{equation*}
\left\{
  \begin{array}{l}
  \ds \cE_{1}
\,=\, \left(
\tau_0-\frac{\tau(\eps)}{\eps^2}
\right)\,
\left\langle 
\cA^*\cA\,D_{\tau_0,0}\,,\,
v^\eps\right\rangle\,,
  \\[1.1em]
  \ds\cE_{2}
\,=\,\frac{\tau(\eps)^2}{\eps^3}\,
\left\langle 
\cA^*\cA\, D_1\,,\,
v^\eps\right\rangle\,,
\\[1.1em]
  \ds\cE_{3}
\,=\,\sqrt{2}\,
\frac{\tau(\eps)}{\eps^2}\,
\left\langle 
\left(\cA^*\right)^2
D_2\,,\,
v^\eps\right\rangle\,.
  \end{array}\right.
\end{equation*} 
We rewrite $\cE_{1}$, $\cE_{2}$ and $\cE_{3}$ according to the following considerations: first, we notice that $D_{\infty,0}$ solves \eqref{prop2:A} and therefore add $D_{\infty,0}$ to the left hand side of the bracket in $\cE_{1}$, second we apply the duality formula \eqref{prop1:A} in $\cE_{1}$, $\cE_{2}$ and $\cE_{3}$ and then replace $v^\eps$ in  $\cE_{1}$ and $\cE_{2}$ according to the relation
\[
\cA^\star \cA\,v^\eps\,=\,D_0
+\frac{\tau(\eps)}{\eps}\,\cA^\star D_1
-D_{\tau_0,0}\,.
\]
Hence, we obtain
\begin{equation*}
	\left\{
	\begin{array}{l}
		\ds \cE_{1}
		\,=\,
		\left(
		\tau_0-\frac{\tau(\eps)}{\eps^2}
		\right)\,
		\left\langle 
		D_{\tau_0,0}-D_{\infty,0}\,,\,
		 D_0
		+\frac{\tau(\eps)}{\eps}\,\cA^\star D_1
		-D_{\tau_0,0}\right\rangle\,,
		\\[1.1em]
		\ds\cE_{2}
		\,=\,\frac{\tau(\eps)^2}{\eps^3}\,
		\left\langle 
		D_1\,,\,
		 D_0
		+\frac{\tau(\eps)}{\eps}\,\cA^\star D_1
		-D_{\tau_0,0}\right\rangle\,,
		\\[1.1em]
		\ds\cE_{3}
		\,=\,\sqrt{2}\,
		\frac{\tau(\eps)}{\eps^2}\,
		\left\langle 
		D_2\,,\,\cA^2\,v^\eps\right\rangle\,.
	\end{array}\right.
\end{equation*} 
To estimate $\cE_{1}$, we apply Young's inequality, which yields
\[
\cE_{1}\,\leq\,\frac{\eta}{2}\frac{\tau(\eps)}{\eps^2}
\, \|  D_0
+\frac{\tau(\eps)}{\eps}\,\cA^\star D_1
-D_{\tau_0,0}\|_{L^2}^2
\,+\,
\frac{1}{2\eta}
\frac{\eps^2}{\tau(\eps)}
\left|
\tau_0
-
\frac{\tau(\eps)}{\eps^2}
\right|^2\|D_{\tau_0}-D_{\infty}\|_{L^2}^2\,,
\]
for all positive $\eta$. To estimate $\cE_{2}$, we apply Young's inequality and then assumption \eqref{cond:tau} which ensures that $\ds \tau(\eps)^3/\eps^4\,\leq\,\left(\overline{\tau}_0^2\,\tau(\eps)\right)/\eps^2$, this gives
\[
\cE_{2}\,\leq\,\frac{\eta}{2}\frac{\tau(\eps)}{\eps^2}
\, \|  D_0
+\frac{\tau(\eps)}{\eps}\,\cA^\star D_1
-D_{\tau_0,0}\|_{L^2}^2
\,+\,
\frac{1}{\eta}\frac{\tau(\eps)}{\eps^2}\,\overline{\tau}_0^2\,\left\|D_{\perp}\right\|^2_{L^2}\,,
\]
for all positive $\eta$. To estimate $\cE_{3}$, we apply Young's inequality and then bound the norm of $\cA^2\,v^\eps$ by applying item \eqref{lem:111} in Lemma \ref{lem:1} with source term
\[g\,=\,D_0
+\frac{\tau(\eps)}{\eps}\,\cA^\star D_1
-D_{\tau_0,0}\,,
\]
it yields
\[
\cE_{3}\,\leq\,\eta\,\frac{\tau(\eps)}{\eps^2}
\, \|  D_0
+\frac{\tau(\eps)}{\eps}\,\cA^\star D_1
-D_{\tau_0,0}\|_{L^2}^2
\,+\,
\frac{C}{\eta}\frac{\tau(\eps)}{\eps^2}\,\left\|D_{\perp}\right\|^2_{L^2}\,,
\]
for some constant $C$ depending only on $\Phi$ and $T_0$. We gather the latter estimates, take $\eta=1/4$ and apply item \eqref{lem:110} in Lemma \ref{lem:1}, which ensures that
\[
\cE
\,\leq\,
\frac{C_P^2}{2}\,
\| D_0
+\frac{\tau(\eps)}{\eps}\,\cA^\star D_1
-D_{\tau_0,0}\|_{L^2}^2\,.
\]
Therefore, we obtain
\begin{align*}
\frac{\dD \cE }{\dD t}
+\,\frac{\tau(\eps)}{C_P^2\,\eps^2}
\, \cE
\,\leq
\,C\,
\frac{\tau(\eps)}{\,\eps^2}
\,\left(1+\overline{\tau}_0^2\right)\left\|D_{\perp}\right\|^2_{L^2}
\,+\,
C\,
\frac{\eps^2}{\tau(\eps)}
\left|
\tau_0
-
\frac{\tau(\eps)}{\eps^2}
\right|^2\,\|D_{\tau_0}-D_{\infty}\|_{L^2}^2\,,
\end{align*} 
for some constant $C$ depending only on $\Phi$ and $T_0$. Then we multiply the latter estimate by $\ds \exp{
\left(
\frac{\tau(\eps)}{C_P^2\,\eps^2}\,t
\right)
}$ and integrate with respect to time. After applying \eqref{cv:D0} to estimate $\ds \|D_{\tau_0} -D_{\infty}\|_{L^2}$ and the first result in item $(i)$ of Theorem \ref{th:main2} to estimate the norm of $D_{\perp}$, it yields
\begin{align*}
 \cE(t) \,&\leq \, 
\left(
\cE(0)
\,+\,C\,\frac{\tau(\eps)^2}{\eps^2}\,(\overline{\tau}_0^6+1)\left\|D(0)-D_{\infty}\right\|_{H^{1}}^2\right)\,
\exp{\left(
-
\frac{\tau(\eps)}{\eps^2}\,
\kappa\, t\right)}\\
&+
C\,\left|
\frac{\tau_0\,\eps^2}{\tau(\eps)}
-
1
\right|^2\,\|D_{\tau_0}(0)-D_{\infty}\|_{L^2}^2\,
\left(
\frac{2\,\tau_0\,\eps^2}{\tau(\eps)}
-
1
\right)^{-1}\,
\exp{\left(
-
\frac{\tau(\eps)}{\eps^2}\,
\kappa\, t\right)
}
\,. 
\end{align*}
To conclude, we substitute $\cE(t)$ (resp. $\cE(0)$) in the latter estimate according to \eqref{equiv4} (resp. \eqref{equiv3}) in Lemma \ref{lem:4} and then apply assumption \eqref{cond3:tau} on $\tau(\eps)$, which ensures $\ds\left(
\frac{2\,\tau_0\,\eps^2}{\tau(\eps)}
-
1
\right)^{-1}\,\leq\, 3$, this yields 
\begin{align*}
& \|D_0(t) -D_{\tau_0,0}(t)\|_{H^{-1}}^2 \,\leq \, 
  \\
  &\quad \ds C\left(
\|D_0(0) -D_{\tau_0,0}(0)\|_{H^{-1}}^2
\,+\,\frac{\tau(\eps)^2}{\eps^2}\,(\overline{\tau}_0^6+1)\left\|D(0)-D_{\infty}\right\|_{H^{1}}^2\right)
e^{
-
\frac{\tau(\eps)}{\eps^2}\,
    \kappa\, t}\, +
  \\
&\quad \ds C\left|
\frac{\tau_0\,\eps^2}{\tau(\eps)}
-
1
\right|^2\|D_{\tau_0}(0)-D_{\infty}\|_{L^2}^2\,
e^{
-
\frac{\tau(\eps)}{\eps^2}\,
\kappa\, t}
\,. 
\end{align*}
We obtain the second estimate provided in $(i)$ of Theorem \ref{th:main2} taking the square root in the latter estimate and applying assumption \eqref{cond3:tau} in order to substitute $\tau(\eps)$ with $\tau_0\,\eps^2$.

To prove the second item $(ii)$ of Theorem \ref{th:main2}, we  follow the same
lines as the ones for item $(i)$ replacing  $D_{\tau_0}$ by $D_\infty$ and observing that $D_\infty$ also
solves the equation \eqref{eq:lim0} since it is a stationary
solution. Therefore, computations are even simpler since the term $\cE_{1}$
vanishes in this case. As a consequence the estimate provided in item $(ii)$
follows. 

%
\section{Finite volume discretization for the space variable}
\label{sec:3}
\setcounter{equation}{0}
\setcounter{figure}{0}
\setcounter{table}{0}

In this section we present a finite volume scheme for \eqref{Hermite:D}. Then we prove discrete hypocoercive estimates on
the discrete solution to investigate the long time behavior and the
speed of convergence to the steady state. Finally, we prove an
asymptotic preserving property for the diffusive limit taking
$\tau(\eps)\sim\tau_0\, \eps^2$ with error estimates with respect to $\eps$. Thanks to the groundworks laid in the previous Section, we are able to propose a scheme which describes all the variety of regimes that we aim to capture in this article.

\subsection{Numerical scheme}
For simplicity purposes, we consider the problem
in one space dimension. It will be straightforward to generalize this
construction for Cartesian meshes in multidimensional case. In a
one-dimensional setting, we consider an interval $(a,b)$ of $
\mathbb{R}$ and  for  $N_{x}\in\N^\star$, we introduce the set
$\J=\{1,\ldots, N_x\}$ and  a family of control volumes
$\left(K_{j}\right)_{j\in\J}$ such that
$K_{j}=\left]x_{j-{1}/{2}},x_{j+{1}/{2}}\right[$ with $x_{j}$ the
middle of the intervall $K_j$ and 
\begin{equation*}
a=x_{{1}/{2}}<x_{1}<x_{{3}/{2}}<...<x_{j-{1}/{2}}<x_{j}<x_{j+{1}/{2}}<...<x_{N_{x}}<x_{N_{x}+{1}/{2}}=b\,.
\end{equation*}
Let us set
$$
\left\{
  \begin{array}{l}
\ds\Delta x_{j}=x_{j+{1}/{2}}-x_{j-{1}/{2}}, \,\text{ for } j \in\J\,,
\\[0.9em]
\ds\Delta x_{i+{1}/{2}} = x_{j+1}-x_{j}, \,\text{ for } 1 \leq j \leq N_{x}-1\,.
  \end{array}\right.
$$
We also introduce the parameter $h$ such that
$$
h \,=\, \max_{j \in\J} \Delta x_j\,.
$$
Let $ \Delta t$ be the time step. We set $t^{n}=n \Delta t$ with
$n\in\N$. A time discretization of $\R^+$ is then given  by the
increasing sequence of $(t^{n})_{n \in\N}$. In the sequel, we will
denote by $D_k^n$ the approximation of $D_k(t^n)$, where the index $k$
represents the $k$-th mode of the Hermite decomposition, whereas
$\cD_{k,j}^n$ is an approximation of the mean value of $D_k$
over the cell $K_{j}$ at time $t^n$.

First of all, the initial condition is discretized on each cell $K_{j}$ by:
\begin{equation*}
\cD_{k,j}^{0}=\frac{1}{\Delta x_{j}}\int_{K_{j}}D_{k}(t=0,x)\, dx, \quad j\in\J\,.
\end{equation*}
The finite volume scheme is obtained by integrating the equation
\eqref{Hermite:D} over each control volume $K_{j}$ and over each time
step. Concerning the time discretization, we can choose any implicit
method (backward Euler, Implicit Runge-Kutta,...). Since in this paper we are
interested in the spatial discretization, we will only consider a backward Euler method afterwards. Let us now focus on the spatial
discretization.

By integrating equation
\eqref{Hermite:D} on $K_{j}$ for $j\in\J$, we obtain the numerical scheme: for $D_k^n=(\cD_{k,j}^n)_{j\in\J}$
\begin{equation}
	\frac{D_{k}^{n+1} -D_{k}^{n} }{\Delta t}  \,+\,
	\frac{1}{\eps}\left(\sqrt{k}\,\mathcal{A}_{h}\,D_{k-1}^{n+1} \,-\,\sqrt{k+1}\,\mathcal{A}_{h}^\star\,D_{k+1}^{n+1}\right)  \,=\,
	-\frac{k}{\tau(\eps)}\, D_{k}^{n+1}\,, 
	\label{discrete}
\end{equation}
where $\mathcal{A}_{h}$ (resp. $\mathcal{A}_{h}^\star$)  is an
approximation of the operator $\cA$ (resp. $\cA^\star$) given by
\be
\label{eq:Ah}
\cA_h = (\cA_j)_{j \in\J} \quad{\rm and}\quad \cA_h^\star = (\cA_j^\star)_{j \in\J}\,. 
\ee
and where for
$D=(\cD_{j})_{j\in\J}$ it holds
\be
\label{def:Ah}
\left\{
  \begin{array}{l}
\ds\cA_{j} D \,=\, +\sqrt{T_0}\, \left( \frac{\cD_{j+1} - \cD_{j-1}}{2\Delta
    x_j}  \,-\, \frac{E_j}{2\,T_0}\, \cD_{j}\right)\,, \quad j\in\J\,,
    \\[1.1em]
 \ds\cA_{j}^\star D \,=\, -\sqrt{T_0}\, \left( \frac{\cD_{j+1} - \cD_{j-1}}{2\Delta
    x_j}  \,+\, \frac{E_j}{2\,T_0}\, \cD_{j}\right)\,, \quad j \in\J\,,
  \end{array}\right.
\ee
whereas the discrete electric field $E_j$ is given by
\be
\label{def:Ei}
E_j \,=\, -\frac{\Phi_{j+1}-\Phi_{j-1}}{2\Delta x_j} \,=\, \frac{2\,T_0}{\sqrt\rho_{\infty,j}}\,\frac{\sqrt\rho_{\infty,j+1}- \sqrt\rho_{\infty,j-1}}{2\,\Delta x_j}\,,
\ee
where $\rho_{\infty,j}$ is an approximation of the stationary density $\rho_{\infty}$ on the cell
$K_j$. This latter formula is consistent with the definition of
$\sqrt\rho_{\infty}= c_0\, e^{-\Phi/(2T_0)}$ and the fact that
$$
\frac{1}{2\,T_0}\,\partial_x\Phi \,=\, -\frac{1}{\sqrt\rho_{\infty}} \, \partial_x\sqrt\rho_{\infty}\,.
$$

This choice of discretization is motivated by preserving at the discrete level the key
properties \eqref{prop1:A}-\eqref{ineg}. In the end, we propose the following approximation of the continuous solution $f$ to \eqref{vlasov}
\[
f^n(x,v)\,=\,
\sum_{k\in \N} \sqrt{\rho}_{\infty}(x)\,D_{k}^n(x)\,\Psi_k(v)\,,
\]
where for each $k\geq 0$ and $n\geq 0$, we define a piecewise constant
function $D_k^n$ from the numerical values $(\cD_{k,j}^n)_{j\in\J}$ as
$$
D_{k}^n(x) = \cD_{k,j}^n, \qquad x\in K_j\,.
$$
In this context the equilibrium $D_{\infty}$ is given by  
\begin{equation}\label{D:infty:h}
	D_{\infty,k} \,=\,
	\left\{
	\begin{array}{l}
		\sqrt\rho_\infty, \,\, {\rm if } \,\, k=0\,,
		\\[0.9em]
		0, \, \, {\rm else\,;}
	\end{array}
	\right.
\end{equation}
as for the limit in the diffusive regime $D_{\tau_0}^n=(D_{\tau_0,k}^n)_{k\in\N}$, it is given by
\be
\label{disc:eq:lim1}
D_{\tau_0,k}^n = \left\{
\begin{array}{l}
	D_{\tau_0,0}^n, \,\,{\rm if}\,\,k=0\,,
	\\[1.em]
	0, \,\, {\rm else\,,} 
\end{array}\right.
\ee
where $D_{\tau_0,0}^n$ solves the following discrete version of equation \eqref{eq:lim0}
\be
\label{disc:eq:lim0}
\frac{D_{\tau_0,0}^{n+1} -D_{\tau_0,0}^{n} }{\Delta t}   \,+\, \tau_0\,\mathcal{A}_{h}^\star \mathcal{A}_{h} D_{\tau_0,0}^{n+1}   \, =\, 0\,.
\ee
We now introduce the norms we will work with in this section. We denote by $\langle.,.\rangle$ the $L^2$ scalar product for
any  $u=(u_{j})_{j\in\J}$ and $v=(v_{j})_{j\in\J}$,
$$
\langle u, \, v \rangle \,=\, \sum_{j\in\J}  \Delta x_j\, u_{j}\, v_j
$$
and
$$
\|u\|_{L^2} \,=\,   \left(\sum_{j\in\J}  \Delta x_j\, u_{j}^2\right)^{1/2}.
$$
As in the \eqref{compatibility:g}, we consider the following $H^{-1}$ norm defined on the $L^2$ subspace orthogonal to $\sqrt{\rho}_{\infty}$: for all $g_h =\left(g_j\right)_{j \in \cJ}$ which meets the condition
\be
\sum_{j\in\J} \Delta x_j \,g_j \, \sqrt\rho_{\infty,j}\,=\,0\,,
\label{brac}
\ee
we set 
\[
\left\|g_h\right\|_{H^{-1}}\,=\,\left\|\cA\,u_h\right\|_{L^2\left(\T\right)}\,,
\]
where $\ds u_h=(u_j)_{j\in\J}$ is the solution to the discrete equivalent of equation \eqref{eq:elliptic}
\be
\label{eq:elliptic:h}
\left\{
\begin{array}{l}
	\ds (\cA_h^\star\,\cA_h) \,u_h\ =\ g\,,
	\\[1.1em]
	\ds\sum_{j\in\J}\Delta x_j \, u_j \,\sqrt\rho_{\infty,j}\,=\,0\,.
\end{array}\right.
\ee
We also use the $H^1$ norm, analog to the one given in \eqref{def:B}, defined for all $\ds D\,=\,\left(D_k\right)_{k \in \N}$ as follows
\[
\|\cB_h\,D\|_{L^2}^2
\,=\,
\sum_{k\in\N} \|\cB_{k}\,D_{k}\|_{L^2}^2\,,
\]
where the family of discrete operator 
$
\cB_h\,=\,\left(\cB_{h,k}\right)_{k\,\geq\,0}$ is given as follows
\be
\label{def:Bh}
\cB_{h,k} = \left\{
\begin{array}{l}
	\cA_h\,,   \,{\rm if} \, \,k\,=\, 0\,,
	\\[0.9em]
	\cA_h^\star\,, \,\,{\rm else\,.}
\end{array}\right.
\ee
To conclude with this section, we take the same definition of $D_{\perp}$ as in the continuous setting.

\subsection{Main results}
We can now release the two results that constitute the core of this article. Thanks to our choice of discretization, they are an exact translation of their continuous analogs, Theorems \ref{th:main1} and \ref{th:main2}, into the discrete setting, without any loss of accuracy nor uniformity with respect to the parameters at play in our analysis. On top of that, the results are also uniform with respect to the discretization parameters.\\

This first result is the continuous analog of Theorem \ref{th:main1}, it ensures that our scheme has the same long time behavior as the continuous model
\begin{theorem}
  \label{th:main1:h}
  	Suppose that condition \eqref{cond:tau}  on $\tau(\eps)$ is satisfied and Let $D^n=(D_{k}^n)_{k\in\N}$ be the solution to \eqref{discrete}. The following statements hold true
  	\begin{enumerate}
        \item[$(i)$]
          there exists some positive constant $C_0$ depending only on $\Phi$ and $T_0$ such that for all $\eps\,>\,0$ and all $n\,\geq\,0$, we have
  		$$
  		\left\|D^n\,-\,D_\infty\right\|_{L^2} \,\leq \,
  		\sqrt{3}\, \left\|D^{0}\,-\,D_\infty\right\|_{L^2}\left(
  		1\,+\,
  		\frac{\tau(\eps)}{\eps^2}\,
  		\kappa_0\, \Delta t
  		\right)^{-n/2}\,;
  		$$
              \item[$(ii)$]
                suppose in addition
  		that the mesh is regular enough so that the quantity 
  		\begin{equation}\label{hyp:mesh}
  			R_h\,=\,\sup_{(i,j) \in \cJ^2} \left|\Delta x_j \Delta x_i^{-1} -1\right|
  		\end{equation}
  		stays uniformly bounded with respect to the discretization parameter $h$. Then  there exists a positive constant
  		$C_1$ (depending only on $\Phi$, $T_0$ and $R_h$) such that that for all $\eps\,>\,0$ and all $n\,\geq\,0$, we have
  		$$
  		\left\|\cB_h D^n\right\|_{L^2} \,\leq \, \sqrt{3}\,\left( C_1
  		\left(
  		\overline{\tau}_0+1
  		\right)
  		\left\|\cB_h D^{0}\right\|_{L^2}  \,+\,  \left\|D^{0}\,-\,D_\infty\right\|_{L^2}\right)\left(
  		1\,+\,
  		\frac{\tau(\eps)}{\eps^2}\,
  		\kappa_1\, \Delta t
  		\right)^{-\frac{n}{2}}\,,
  		$$
  	\end{enumerate}
 In the previous estimates $\kappa_i >0$ is
 given by
 $$
 \kappa_i \,=\, \frac{1}{
 	C_i\,(\overline{\tau}_0^2+1)
 } \,.
 $$
\end{theorem}

Our second result deals with the asymptotic $\eps \rightarrow 0$, it is the discrete analog of Theorem \ref{th:main2}
\begin{theorem}
 \label{th:main2:h}
Suppose that $\tau(\eps)$ meets assumption \eqref{cond:tau} and that the mesh meets assumption \eqref{hyp:mesh}. Consider the solution $D^n=(
D^n_{k} )_{k\in\N}$ to \eqref{discrete}. 
The following statements hold true uniformly with respect to $\eps$
\begin{itemize}
\item[$(i)$]  suppose that $\tau(\eps)$ satisfies
\eqref{cond2a:tau}  and \eqref{cond3:tau} and
consider $D^n_{\tau_0}=(D^n_{\tau_0,k})_{k\in\N}$ given by
\eqref{disc:eq:lim1}. Then it holds for all $n\geq0$,
\begin{equation*}
\left\|D_{\perp}^n\right\|_{L^2} \,\leq \,
\left\|D_{\perp}^0\right\|_{L^2}\,\left(1+\frac{\Delta t}{2\,\tau_0\, \eps^2}\right)^{-\frac{n}{2}}+\,\tau_0\,\eps\,C(\overline{\tau}_0+1)\,\left\|D^0-D_\infty\right\|_{H^1}\,
\left(1+\tau_0\,\kappa\,\Delta t\right)^{-\frac{n}{2}}\,,
\end{equation*}
and
\begin{align*}
	\ds
\left\| D^n_0-D^n_{\tau_0,0}\right\|_{H^{-1}}\leq\;
&\ds C\left(
\left\|D^{0}_0 -D_{\tau_0,0}^0\right\|_{H^{-1}}
+\,\eps\,\tau_0\,(\overline{\tau}_0^3+1)\left\|D^{0}-D_{\infty}\right\|_{H^{1}}
\right)
\left(1+\tau_0\kappa\Delta t\right)^{-\frac{n}{2}}\,,\\
&\ds C
\left|
\frac{\tau_0\eps^2}{\tau(\eps)}
-
1
\right|\left\|D_{\tau_0}^0-D_{\infty}\right\|_{L^2}
\left(1+\tau_0\kappa\Delta t\right)^{-\frac{n}{2}}
\,;
\end{align*}
\item[$(ii)$]  suppose that $\tau(\eps)$ satisfies
\eqref{cond2b:tau}. Then it holds for any $n\geq0$
\begin{equation*}
	\left\|D_{\perp}\right\|^2_{L^2} \,\leq \,
	\left\|D_{\perp}\right\|^2_{L^2}\,\left(1+\frac{\Delta t}{\tau(\eps)}\right)^{-\frac{n}{2}}+\,\frac{\tau(\eps)}{\eps}\,C(\overline{\tau}_0+1)\,\left\|D^0-D_\infty\right\|_{H^1}\,
	\left(1+\frac{\tau(\eps)}{\eps^2}\kappa\Delta t\right)^{-\frac{n}{2}}\,,
\end{equation*}
and
\begin{equation*}
\left\| D^n_0-D^n_{\infty,0}\right\|_{H^{-1}}\leq C\left(
\left\|D^{0}_0 -D_{\infty,0}^0\right\|_{H^{-1}}
+\frac{\tau(\eps)}{\eps}\,(\overline{\tau}_0^3+1)\left\|D^{0}-D_{\infty}\right\|_{H^{1}}\right)
\left(1+\frac{\tau(\eps)}{\eps^2}\kappa\Delta t\right)^{-\frac{n}{2}}.
\end{equation*}
\end{itemize}
In the latter estimate, constant $C$ only depends on $\Phi$, $T_0$ and $R_h$ and exponent $\kappa$ is given by
$$
\kappa \,=\, \frac{1}{
	C\,(\overline{\tau}_0^2+1)
} \,.
$$
Furthermore the shorthand notation $\left\|\cdot\right\|_{H^{1}}$ stands for
\[
\left\|D\right\|_{H^{1}}^2\,:=\,
\left\|\cB D\right\|_{L^2}^2\,+\,
\left\|D\right\|_{L^2}^2\,.
\]
\end{theorem}

The proof of these results follows almost exactly the same lines as the proof of Theorems \ref{th:main1} and \ref{th:main2} thanks to the Lemma \ref{lem:3.1}, which constitutes the keystone of our analysis and which ensures that our discretization $\cA_h$ of operator $\cA$ shares all the important properties \eqref{prop1:A}-\eqref{ineg} of its continuous analog. The only difference comes down to some numerical remainder terms that we easily control applying methods already developed in the continuous section.\\

\subsection{Preliminary  properties}
This section is dedicated to the following fundamental Lemma, which ensures that the key properties \eqref{prop1:A}-\eqref{ineg} of the continuous operator $\cA$ are preserved by its discrete analog $\cA_h$. Thanks to this Lemma, all the computations carried in Section \ref{sec:2} directly translate into the discrete framework.
\bl
\label{lem:3.1}
Consider the discrete operators $\cA_h$ and  $\cA_h^\star$ given in
\eqref{eq:Ah}. Then we have for any $u=(u_j)_{j \in\J}$ and $v=(v_j)_{j \in\J}$
\begin{enumerate}
\item\label{prop1:A:h} preservation of the duality formula
  $$
\left\langle \cA_h u, \, v \right\rangle \,=\,  \left\langle u, \,
  \cA_h^\star v \right\rangle\,;
$$
\item\label{prop2:A:h} preservation of the kernel of operator $\cA_h$
$$
\cA_h D_{\infty,0} \,=\, 0\,,
$$ 
where the equilibrium $D_{\infty}$ is given by \eqref{D:infty:h};
\item\label{cons:mass0:h} preservation of the mass conservation properties
\begin{equation}\label{cons:mass1:h}
	\sum_{j\in\J} \Delta x_j \,\cA^{\star}_j u \, \sqrt\rho_{\infty,j}\,=\,0\,,
\end{equation}
and for all $n\,\geq\,0$, the solution $D_0^n=(\cD_{0,j}^n)_{j\in\J}$ to \eqref{discrete} with index $k=0$ verifies
\begin{equation}\label{cons:mass:h}
	\sum_{j\in\J} \Delta x_j \,\cD^n_{0,j} \, \sqrt\rho_{\infty,j}\,=\,	\sum_{j\in\J} \Delta x_j \,\rho_{\infty,j}\,;
\end{equation}
\item\label{prop:A+Astar:h} preservation of the sum property
\begin{equation*}
\| \left(\cA_h + \cA^\star_h\right) u\|_{L^2} \,\leq \, \frac{1}{\sqrt{T_0}}\,\|\Phi\|_{W^{1,\infty}} \|u\|_{L^2}\,;
\end{equation*}
\item\label{prop3:A:h} preservation with the commutator property
\begin{equation*}
	\|\left[\cA_h,\,\cA_h^\star\right] u\|_{L^2} \,\leq \, C\,\|\Phi\|_{W^{2,\infty}} \|u\|_{L^2}\,,
\end{equation*}
where constant $C$ depends only on $R_h$ (see \eqref{hyp:mesh}), it is explicitly given by
\[C\,=\,
2\,+\,
R_h
\,;
\]
\item conservation of the Poincar\'e-Wirtinger inequality: under
  condition \eqref{brac} on $u$ there exists a constant $C_d>0$
  depending only on $\Phi$ and $T_0$ such that
\be
\label{ineg-disc}
\|u\|_{L^2}  \,\leq\,C_d\,\| \cA_h\, u\|_{L^2}\,.
\ee
\end{enumerate}
\el
\begin{remark}
	When the mesh is regular, item \eqref{prop3:A:h} in Lemma \ref{lem:3.1} may be improved into a consistent estimate compared to its continuous analog \eqref{prop3:A}, indeed we easily obtain
	\begin{equation*}
		\|\left[\cA_h,\,\cA_h^\star\right] u\|_{L^2} \,\leq \, \left(\|\Phi\|_{W^{2,\infty}}\,+\,
		\frac{h}{2}\,\|\Phi\|_{W^{3,\infty}}\right)
		 \|u\|_{L^2}\,,
	\end{equation*}
for any $u=\left(u_j\right)_{j\in\J}$, following the same method as in the proof.
\end{remark}
\begin{proof}
To prove item \eqref{prop1:A:h}, we consider any $(u_j)_{j\in\J}$ and $(v_j)_{j\in\J}$, we have after a discrete integration by part and using
periodic boundary conditions
\begin{eqnarray*}
\left\langle \cA_h u,\,v \right\rangle &=&  \sum_{j\in\J} \Delta x_j\,\cA_{j} u\, v_j \\
&= & \sum_{j\in\J} \sqrt{T_0}\left(
\frac{u_{j+1}-u_{j-1}}{2}\, v_j \,-\,  \Delta x_j\,\frac{E_j}{2\,T_0}\, u_j \, v_j\right) 
                                            \\
  &=&  \sum_{j\in\J} -\sqrt{T_0}\left(
\frac{v_{j+1}-v_{j-1}}{2} \,u_j \,+\,  \Delta x_j\,\frac{E_j}{2\,T_0}\,
      v_j \,u_j\right)
       \, =\,  \left\langle u\,,\cA_h^\star\,v\right\rangle.
\end{eqnarray*}
To prove item \eqref{prop2:A:h}, we look for $D=\left(D_k\right)_{k\in\N}$
such that $\cA_h\,D_{0}\,=\,0$, that is,
$$
0\,=\,\cA_i\,D_{0} \,=\, \frac{\sqrt{T_0}}{2\,\Delta x_j}\, \left( \cD_{0,\,j+1} - \cD_{0,\,j-1} \,+\, \frac{\Phi_{j+1}-\Phi_{j-1}}{2\,T_0}\, \cD_{0,\,j}\right)\,.
$$
Hence, from the particular choice of the discrete electric field
\eqref{def:Ei}, we have that
$$
\frac{\cD_{0,\,j+1} -
  \cD_{0,\,j-1}}{\cD_{0,\,j} } \,-\,  \frac{\sqrt\rho_{\infty,\,j+1} -
  \sqrt\rho_{\infty,\,j-1}}{\sqrt\rho_{\infty,\,j}}  \,=\,0 \,,
$$
which yields to definition \eqref{D:infty:h}. 

We turn to the mass conservation property  \eqref{cons:mass0:h}. According to the definition \eqref{def:Ah} of $\cA^{\star}_h$, it holds
\[
\cA^*_j u\, \sqrt\rho_{\infty,j}\,\Delta
x_j
\,=\,
-\sqrt{T_0}\, \left( \sqrt\rho_{\infty,j}\,\frac{u_{j+1} - u_{j-1}}{2}  \,+\,\frac{\sqrt\rho_{\infty,j+1}- \sqrt\rho_{\infty,j-1}}{2}\,u_{j}\right)\,.
\]
Therefore, relation \eqref{cons:mass1:h} is obtained summing the latter over $j \in \cJ$ and performing a discrete integration by part. Relation \eqref{cons:mass:h} is obtained evaluating equation \eqref{discrete} with index $k=0$ and $j\in \cJ$, multiplying by $\sqrt{\rho}_{\infty,j}\,\Delta x_j$, then summing over $j \in \cJ$ and applying relation \eqref{cons:mass1:h} with $u\,=\,D^{n+1}_1$.

We prove item \eqref{prop:A+Astar:h} taking the $L^2$ norm in the following relation
\[
\sqrt{T_0}\,\left(\cA_j+\cA^\star_j\right)u\,=\,-\,\frac{2\,T_0}{\sqrt\rho_{\infty,j}}\,\frac{\sqrt\rho_{\infty,j+1}- \sqrt\rho_{\infty,j-1}}{2\,\Delta x_j}\,u_j\,,
\]
which holds for any $u=\left(u_j\right)_{j\in\J}$.

We turn to item \eqref{prop3:A:h} and compute the commutator for the discrete operator $[\cA_h,\,\cA_h^\star]$ as
\begin{eqnarray*}
  \left[\cA_h,\,\cA^\star_h\right]_j u &=& (\cA_h\,\cA^\star_h -
                                           \cA_h^\star\,\cA_h)_j u
  \\
  &=&  -\frac{E_{j+1}- E_{j-1}}{4\, \Delta x_j}
\,  \left(u_{j+1}+ u_{j-1}\right) \, -    \frac{E_{j+1}-2\, E_j +  E_{j-1}}{4\,\Delta x_j} \left(u_{j+1}- u_{j-1}\right) \,,
\end{eqnarray*}
and therefore, we deduce item \eqref{prop3:A:h} taking the $L^2$ norm in the latter result.\\

Finally, we prove the Poincar\'e inequality \eqref{ineg-disc}. Consider $u=\left(u_j\right)_{j\in\J}$ which meets condition \eqref{brac} and let us denote by $\overline\rho_{\infty}$ the mean of $\rho_{\infty}$
$$
\overline\rho_{\infty} = \sum_{j\in\J} \Delta x_j \, \rho_{\infty,j}\,.
$$
First using the zero weighted average assumption \eqref{brac} on $u$, we remark that the cross term vanishes and
\begin{align*}
	\|u\|_{L^2}^2&=\sum_{j\in\J}\Delta x_j\,
	\left(\frac{u_j}{\sqrt\rho_{\infty,j}}
	\right)^2\,\rho_{\infty,j}\,,
	\\
	&=\frac{1}{2\,\overline\rho_{\infty}}\, \sum_{j\in\J}\sum_{k\in\J}\Delta x_j \,\Delta x_k\,\left(\frac{u_k}{\sqrt\rho_{\infty,k}}
	\,-\,\frac{u_j}{\sqrt\rho_{\infty,j}} \right)^2\,\rho_{\infty,j}
	\,\rho_{\infty,k}\,, \\
	&= \frac{1}{\overline\rho}_{\infty}\, \sum_{k\in\J}\sum_{j<k}\Delta x_j\,\Delta x_k\,\left(\frac{u_k}{\sqrt\rho_{\infty,k}}  - \frac{u_j}{\sqrt\rho_{\infty,j}}\right)^2\,\rho_{\infty,j}\,\rho_{\infty,k}\,.
\end{align*}
For $j<k$, we have
\[
\frac{u_k}{\sqrt\rho_{\infty,k}}  -\frac{u_j}{\sqrt\rho_{\infty,j}}
\,=\,\sum_{l=j}^{k-1}\frac{u_{l+1}}{\sqrt\rho_{\infty,l+1}} \,-\, \frac{u_l}{\sqrt\rho_{\infty,l}}\,,
\]
which yields
\be
\label{q:1}
\|u\|_{L^2}^2 \,\leq\,  \overline{\rho}_\infty\,\left(\sum_{l\in\J}\frac{u_{l+1}}{\sqrt\rho_{\infty,l+1}}-\frac{u_l}{\sqrt\rho_{\infty,l}}\right)^2\,.
\ee
On the other hand, we set for any $j\in\J$
$$
\sqrt{\overline\rho}_{\infty,j} \,=\,
\frac{\sqrt{\rho}_{\infty,j-1}+\sqrt{\rho}_{\infty,j+1}}{2}, \quad{\rm
	and}\quad
\eta_j \,=\,\frac{\sqrt{\rho}_{\infty,j+1}-\sqrt{\rho}_{\infty,j-1}}{2\, \sqrt{\overline\rho}_{\infty,j}},
$$
and  observe that the discrete operator $\cA_hu$ may be written as 
$$
\frac{\Delta x_j}{\sqrt{\overline\rho}_{\infty,j}}\,\cA_j u \,=\,
\frac{\sqrt{T_0}}{2} \,\left[ \left(\frac{u_{j+1}}{\sqrt\rho_{\infty,j+1}} -
\frac{u_{j}}{\sqrt\rho_{\infty,j}} \right) \,\left( 1+ \eta_j\right)  \,+\, \left(\frac{u_{j}}{\sqrt\rho_{\infty,j}} -
\frac{u_{j-1}}{\sqrt\rho_{\infty,j-1}} \right) \,  \,\left( 1-\eta_j\right) \right]\,.
$$
Then we have using periodic boundary conditions
\begin{eqnarray*}
	\sqrt{T_0}\sum_{j\in\J}\left(\frac{u_{j+1}}{\sqrt\rho_{\infty,j+1}}-\frac{u_j}{\sqrt\rho_{\infty,j}}\right)
	&=&
	\frac{\sqrt{T_0}}{2}\sum_{j\in\J}\left(\frac{u_{j+1}}{\sqrt\rho_{\infty,j+1}}-\frac{u_j}{\sqrt\rho_{\infty,j}}\right)
	\,+\,  \left(\frac{u_{j}}{\sqrt\rho_{\infty,j}}-\frac{u_{j-1}}{\sqrt\rho_{\infty,j-1}}\right)
	\\
	&=&\sum_{j\in\J} \frac{\Delta x_j}{\sqrt{\overline\rho}_{\infty,j}}\,\cA_j u 
	\,-\, \sqrt{T_0}
	\left(\frac{u_{j+1}}{\sqrt\rho_{\infty,j+1}}-\frac{u_{j}}{\sqrt\rho_{\infty,j}}\right)\,\frac{\eta_j
		- \eta_{j+1}}{2}
\end{eqnarray*}
Hence using that $\Phi$ is Lipschitzian, we have
$$
|\eta_{j+1} - \eta_j| \leq C_\Phi \,h,
$$
which yields that
\begin{eqnarray*}
	\sqrt{T_0}\sum_{j\in\J}\left|\frac{u_{j+1}}{\sqrt\rho_{\infty,j+1}}-\frac{u_j}{\sqrt\rho_{\infty,j}}\right|
	\,\leq\,
	\sum_{j\in\J} \frac{\Delta x_j}{\sqrt{\overline\rho}_{\infty,j}} \,\left|\,\cA_j\,u\,\right|
	\,+\, C_\Phi\,h \sqrt{T_0}\sum_{j\in\J} \left|\frac{u_{j+1}}{\sqrt\rho_{\infty,j+1}}-\frac{u_{j}}{\sqrt\rho_{\infty,j}}\right|.
\end{eqnarray*}
On the one hand, we consider the case when $h$ is small enough such that
$1-C_\Phi h \geq 1/2$, we  get that
$$
\sum_{j\in\J}\left|\frac{u_{j+1}}{\sqrt\rho_{\infty,j+1}}-\frac{u_j}{\sqrt\rho_{\infty,j}}\right|
\,\leq\, \frac{2}{\sqrt{T_0}}\,\sum_{j\in\J} \frac{\Delta x_j}{\sqrt{\overline\rho}_{\infty,j}} \,\left|\,\cA_j\,u\,\right| 
$$
On the other hand, when $1-C_\Phi\,h \leq 1/2$ (the space step $h$ is large), we use the fact that in
finite dimension, both semi-norms are equivalent. Thus, there exists a
constant $C^\prime_\Phi>0$, independent of $h$, such that
$$
\sum_{j\in\J}\left|\frac{u_{j+1}}{\sqrt\rho_{\infty,j+1}}-\frac{u_j}{\sqrt\rho_{\infty,j}}\right|
\,\leq\, \frac{C_\Phi^\prime}{\sqrt{T_0}}\,\sum_{j\in\J} \frac{\Delta x_j}{\sqrt{\overline\rho}_{\infty,j}} \,\left|\,\cA_j\,u\,\right|. 
$$
Gathering the  latter result with \eqref{q:1}, it yields
$$
\|u\|_{L^2}^2 \,\leq\, \frac{\left( C^\prime_\Phi\right)^2
\,\overline{\rho}_\infty}{T_0}\,\left(\sum_{j\in\J}
\frac{\Delta x_j}{\sqrt{\overline\rho}_{\infty,j}}
\,\left|\,\cA_j\,u\,\right|\right)^2\,.
$$
Using the Cauchy-Schwarz inequality, we obtain  the result
$$
\|u\|_{L^2}^2\, \leq\, C_d^2\, \|\cA_h\,u\|_{L^2}^2\,,
$$
where $C_d^2$ is given by
$$
C_d^2 \,=\,  \frac{\left(C^\prime_\Phi\right)^2\,{\overline\rho}_{\infty}}{T_0}\, \sum_{j\in\J} \frac{\Delta x_j}{|\sqrt{\overline\rho}_{\infty,j}|^2}. 
$$
\end{proof}

From the latter results, we may now get estimates on the solution $u_h$ to \eqref{eq:elliptic:h} as in Lemma \ref{lem:1} in the continuous setting. 
 
\bl
\label{lem:3.3}
 Let us consider the solution $u_h$ to \eqref{eq:elliptic:h} with source term $g=(g_{j})_{j\in\J}$ satisfying the compatibility assumption \eqref{brac}. Then, $u_h$ satisfies the following estimate
 \be
 \label{lem:3.30}
    \|\cA_h\, u_h\|_{L^2} \,\leq\, C_d\,\|g\|_{L^2}\,,
    \ee
    and
    \be
    \label{lem:3.31}
  \|\cA_h^2\, u_h\|_{l^2}  \,\leq\, \left(1\,+\,\frac{C_d}{\sqrt{T_0}}\,\|\partial_x\Phi\|_{L^\infty}\right)\,\|g\|_{L^2}\,.
\ee
Moreover, consider now $(D^n_k)_{k\in\N}$ solution to
\eqref{discrete} and $u_h^n=(u^n_j)_{j\in\J}$ the
 corresponding solution to \eqref{eq:elliptic:h} with the source term
 $D_0^n-\sqrt{\rho}_\infty$. Then  we define $d_{t} u_h^{n+1}$ as
 \be
 \label{def:dtu}
d_{t} u_h^{n+1} \,=\,   \frac{u_h^{n+1}-u_h^n}{\Delta t}\,,
\ee
which satisfies 
\be
\label{lem:3.32}
\eps\,\left\|\cA_h \, d_t u_h^{n+1}\right\|_{L^2}   \,\leq\, \|D_1^{n+1}\|_{L^2}\,.
\ee
\el
\begin{proof}
  We follow the proof of Lemma \ref{lem:1}, we multiply \eqref{eq:elliptic:h} by $\Delta x_i\, u_i$, sum over $i\in\J$ and apply item \eqref{prop1:A:h} of Lemma \ref{lem:3.1}, it yields 
 \[
  \|\cA_h\, u_h\|_{L^2}^2 \,\leq\,\| D -
  \sqrt\rho_{\infty}\|_{L^2}\,\|u_h\|_{L^2}\,,
\]
hence the discrete Wirtinger-Poincar\'e inequality, obtained in Lemma \ref{lem:3.3}, gives,
$$
\|\cA_h \,u_h\|_{L^2} \,\leq\,C_d\,\| D -
  \sqrt\rho_{\infty}\|_{L^2}\,.
  $$
  For the second estimate, we observe  that
  $$
  (\cA_h\,+\,\cA^\star_h)_j\,u_h \,=\,\frac{\sqrt\rho_{\infty,j+1}-\sqrt\rho_{\infty,j-1}}{2\,\Delta x_j\, \sqrt\rho_{\infty,j}}\, u_j
  $$
  hence we  obtain 
  \begin{eqnarray*}
(\cA^2_h)_j \,u_h &=& -\left(\cA_h^\star\,\cA_h\right)_j u_h \,+\,
 \, \frac{\sqrt\rho_{\infty,j+1}-\sqrt\rho_{\infty,j-1}}{2\,\Delta
                x_j\, \sqrt\rho_{\infty,j}}\,\cA_j\,u_h \\
    &=& -\left(D_{0,j}-\sqrt\rho_{\infty,j}\right)\,+\, \frac{\sqrt\rho_{\infty,j+1}-\sqrt\rho_{\infty,j-1}}{2\,\Delta x_j\, \sqrt\rho_{\infty,j}}\,\cA_j\, u_h\,.  
\end{eqnarray*}
Since $\Phi$ is Lipschitzian and applying \eqref{lem:3.30}, we obtain the result
$$
\|\cA^2_h\,u_h \|_{L^2} \,\leq \, C\,\|D(t) -\sqrt\rho_{\infty}\|_{L^2}\,.
$$
For the third estimate we consider now the solution $D^n=(D_{k}^n)_{k\in\N}$ to
\eqref{discrete} and  $u_h^n$ the solution to \eqref{eq:elliptic:h} with source term $D_0^n-\sqrt{\rho}_\infty$. We  get for any $j\in\J$,
$$
  \left(\cA^\star_h\cA_h\right)_j\, d_tu^{n+1}_h\, =\, \frac{\cD_{0,j}^{n+1} -
  \cD_{0,j}^n}{\Delta t} \,=\, -\frac{1}{\eps}\,\cA^\star_j\, D_{1}^{n+1}\,.
  $$ 
Then we multiply by $\Delta x_j\,d_t u^{n+1}_h$, sum over $j\in\J$ and use
\eqref{prop1:A} to get 
\[
 \left\|\cA_h \,d_t u^{n+1}_h\right\|_{L^2}^2\, =\, -\frac{1}{\eps}\left\langle
   D^{n+1}_1,\, \cA_h\,  d_tu_h^{n+1} \right\rangle
 \,\leq\,\frac{1}{\eps}\,\|D^{n+1}_1\|_{L^2}\,\left\| \cA_h \,d_t u^{n+1}_h \right\|_{L^2}\,.
\]
\end{proof}

\subsection{Proof of Theorem \ref{th:main1:h}}
We split the proof of  Theorem \ref{th:main1:h}  into
two steps corresponding to the $L^2$ and $H^1$ convergence
result. Thanks to Lemma \ref{lem:3.3}, the method followed in Section
\ref{sec:2} to prove the continuous analog to this result (Theorem
\ref{th:main1}) directly applies here, excepted for some additional
numerical remainders  for which we give a detailed method in order to
get control over.

We define $\cH_0^n$ as
\begin{equation}
\cH_0^n \,=\,
\frac12\,\|D^n-D_\infty\|_{L^2}^2\,+\, \alpha_0\,\left\langle
  \frac{\tau(\eps)}{\eps}\,\cA_h^\star D_1^n, u_h^n\right\rangle\,,
 \label{eq:H0:discrete}
\end{equation}
where $u^n$ is solution to \eqref{eq:elliptic:h} with
$D^n_0-\sqrt\rho_\infty$ as a source term. First let us point out that $\cH_0^n$ shares the same properties as its continuous analog, indeed it holds 
\bl
Suppose that condition \eqref{cond:tau}  on $\tau(\eps)$ is satisfied.  Then for all
$\alpha_0\in(0,\overline{\alpha}_0)$,  with $\overline{\alpha}_0=1/(4\,\overline{\tau}_0\,C_d)$ 
and $D^n=(\cD^n_{k,j})_{j\in\J\,,\,k\in\N}$, one has
 \begin{equation}
   \label{disc:equiv}
 \frac{1}{4}\,\|D^n -D_\infty\|_{L^2}^2 \leq\,\cH_0^n\,\leq \,\frac{3}{4}\,\|D^n -D_\infty\|_{L^2}^2\,.
 \end{equation}
 \label{lem:2disc}
\el
\begin{proof}
The proof follows the same lines as the one of Lemma \ref{lem:2}.
\end{proof}

We are now able to proceed to the proof of the first item $(i)$ of Theorem \ref{th:main1:h}.
On the one hand, proceeding as the proof of item $(i)$  in Theorem \ref{th:main1}, it yields from Lemma \ref{lem:3.1}
\be
\label{estim:L2:discrete}
\frac{\cH_0^{n+1} -\cH_0^{n}}{\Delta t}  \,=\, \cI_1^{n+1}\,+\, \alpha_0\,\cI_2^{n+1}\,+\, \alpha_0\,\cI_3^{n+1} \, -\, \cR_0^{n+1}\,,
\ee
where
$$
\cI_1^{n+1} \,=\, -\frac{1}{\tau(\eps)}\,\sum_{k\in\N^\star} k\,\left\|D^{n+1}_k\right\|^2_{L^2} 
$$
whereas the other terms correspond to the additional term of the
 modified relative entropy,
 $$
 \left\{\begin{array}{l}
   \ds\cI_2^{n+1}\,:=\,
   -\frac{\tau(\eps)}{\eps^2}\,\left\langle\cA_h^\star\cA_h\,\left(D^{n+1}_0-\sqrt\rho_{\infty}\right) \,-\,
   \sqrt{2}\,(\cA_h^\star)^2\, D_2^{n+1} ,\, u_h^{n+1}\right\rangle  \,
   \textcolor{purple}{-}\, \frac{1}{\eps} \left\langle\cA_h^\star\, D^{n+1}_1
   ,\, u_h^{n+1}\right\rangle\,,
   \\[1.em]
  \ds \cI_3^{n+1}\,:=\, +\frac{\tau(\eps)}{\eps}\,\left\langle \cA_h^\star
  \, D^{n+1}_1,\,d_tu^{n+1}_h\right\rangle\,,
\end{array}\right.
$$
where $d_tu^{n+1}_h$ is given in \eqref{def:dtu}  and $\cR_0$ is a purely numerical remainder given by
\be
\label{7:ju0}
\cR^{n+1}_0 \,=\,  \frac{1}{2\,\Delta t}\| D^{n+1}- D^n\|_{L^2}^2 \,+\,
  \alpha_0\,\frac{\tau(\eps)}{\eps}\,\left\langle
    \cA_h^\star \left(D_1^{n+1}-D_1^n\right),\, d_t u_h^{n+1} \right\rangle\,.
\ee
Both terms  $\cI_2^{n+1}$ and $\cI_3^{n+1}$ can be estimated as in
the proof of item $(i)$  in Theorem \ref{th:main1}, which yields
$$
  \cI_2^{n+1} \,\leq\, -\frac{\tau(\eps)}{\eps^2}\,\left(1\,-\,
       C\,\eta\right) \, \| D^{n+1}_0-D_{\infty,0}
\|_{L^2}^2 \,+\,  \frac{C}{2\,\eta}\, \left( \frac{\tau(\eps)}{\eps^2}\| D^{n+1}_2\|_{L^2} ^2 \,+\,  \frac{1}{\tau(\eps)}\, \| D^{n+1}_1\|_{L^2}^2\right)\,,
$$
for any positive $\eta$ and for some positive constant $C$ depending
only on $T_0$ and $\Phi$ and
$$
\cI^{n+1}_3 \,\leq\, \frac{\tau(\eps)}{\eps^2}\|D^{n+1}_1\|_{L^2}^2\,.
$$
From these latter estimates and taking $\eta=1/(2C)$ and  as long as
\[
\alpha_0
\,<\,
\frac{1}{
C\,(\overline{\tau}_0^2+1)
}\,,
\]
for $C$ great enough and taking $\kappa_0$ such that $3\,\kappa_0/4\,=\,\alpha_0/2$,
we get that
\begin{equation*}
  \frac{\cH_0^{n+1}-\cH_0^{n}}{\Delta t}
  \,+\,
\frac{\tau(\eps)}{\eps^2}\,\kappa_0\,\cH_0^{n+1} \,\leq\,-\cR^{n+1}_0\,.
\end{equation*}
Now we treat the remainder term $\cR^{n+1}_0$, observing that
\begin{equation*}
\left|
 \left\langle \cA_h^\star
  \left(D_1^{n+1}-D_1^n\right),\, d_t u_h^{n+1}
\right\rangle\right| \leq \frac{1}{2\,\Delta t}\,
\left(
\| D^{n+1}_1- D^n_1\|_{L^2}^2
\,+\,
\| \cA_h\left(u_h^{n+1} - u_h^{n}\right) \|_{L^2}^2
\right)\,.
 \end{equation*}
Therefore, applying \eqref{lem:3.30} in Lemma \ref{lem:3.3} with source term $D^{n+1}_0- D^n_0$, we obtain
\begin{equation*}
 \left|
 \left\langle \cA_h^\star
  \left(D_1^{n+1}-D_1^n\right),\, d_t u_h^{n+1}
\right\rangle\right| \leq \frac{1+C_d^2}{2\,\Delta t}\,
\| D^{n+1}- D^n\|_{L^2}^2\,.
 \end{equation*}
Since $\ds\tau(\eps)$ meets assumption \eqref{cond:tau}, the latter estimate ensures that, as long as $\ds\alpha_0\,\leq\,
\left(
\overline{\tau}_0
\,
(1+C_d^2)
\right)^{-1}
$, it holds
\[
0\,\leq\,\cR_0^{n+1}\,,
\]
which yields
\begin{equation*}
  \frac{\cH_0^{n+1}-\cH_0^{n}}{\Delta t}
  \,+\,
\frac{\tau(\eps)}{\eps^2}\,\kappa_0\,\cH_0^{n+1} \,\leq\,0\,.
\end{equation*}
The result follows by applying a discrete Gronwall's lemma and then applying Lemma \ref{lem:2disc} in order to substitute $\cH^n_0$ with the $L^2$ norm of $D^n\,-\,D_{\infty}$ in the latter estimate.

Now we turn to the proof of the second item $(ii)$ of Theorem
\ref{th:main1:h}. Following Section \ref{sec:23}, we   introduce $\cH_1^n$ given by
\begin{equation}
\cH_1^n \,=\,
\frac12\,\|\cB_h D^n \|_{L^2}^2\,+\, \alpha_1\,\left\langle
 \frac{\tau(\eps)}{\eps}\,\cA_h D_0^n, D^n_1\right\rangle\,,
 \label{eq:H1h}
\end{equation}
where $\alpha_1$ has to be determined. Once again, $\cH^n_1$ shares the same properties as its continuous analog 
\bl\label{lem:3:h}
Suppose that condition \eqref{cond:tau}  on $\tau(\eps)$ is satisfied.  Then for all
$\alpha_1\in(0,\overline{\alpha}_1)$,  with $\overline{\alpha}_1=1/(2\,\overline{\tau}_0)$
and $D^n=(D^n_{k})_{k\in\N}$, one has
\begin{equation*}
\label{disc:equiv2}\|\cB_h D^n\|_{L^2}^2
-
\|D^n -D_\infty\|_{L^2}^2
\leq\,4\,\cH_1^n\,\leq \,3\,\|\cB_h D^n\|_{L^2}^2
+\|D^n -D_\infty\|_{L^2}^2\,.
\end{equation*}
\el 
\begin{proof}
The result is obtained applying the same method as in the proof of Lemma \ref{lem:3}.
\end{proof}

We now compute the variation of the modified relative entropy between one
time step from $t^n$ to $t^{n+1}$ and split it into three terms
\[
\frac{\cH_1^{n+1} -\cH_1^{n}}{\Delta t}\, =\, \cJ_1^{n+1} +
\alpha_1\,\cJ_2^{n+1}\, -\, \cR_1^{n+1} ,
\]
where $\cJ_1^{n+1}$ is given by 
$$
  \cJ_1^{n+1} \,:=\, -\frac{1}{\eps}\,\sum_{k\geq 2} \sqrt{k}\,
   \left\langle\, [\cA_h^\star, \cA_h]\,D_{k-1}^{n+1},\, \cA_h^\star
                   D_{k}^{n+1}\right\rangle \,-\,   \frac{1}{\tau(\eps)}\,\sum_{k\in\N^\star} k\,\left\|\cB_{h,k} D^{n+1}_k\right\|^2_{L^2}
$$
 and
 \begin{eqnarray*}
 \cJ_2^{n+1} &:=&
   \frac{\tau(\eps)}{\eps^2}\,\left( \left\langle\cA_h\cA_h^\star D^{n+1}_1 ,\, D^{n+1}_1\right\rangle  \,-\,
   \left\|\cA_h D^{n+1}_0\right\|_{L^2}^2  
   \,+\,\sqrt{2} \left\langle\cA_h D_0^{n+1},\,\cA_h^\star
                  D_2^{n+1}\right\rangle\right)  \\
   &-&\frac{1}{\eps} \left\langle D^{n+1}_1
   ,\,\cA_h D^{n+1}_0\right\rangle\,
 \end{eqnarray*}
 whereas $\cR^n_1$ is given by
\begin{equation}\label{7:ju02}
\cR_1^{n+1} \,=\,  \frac{1}{\Delta t}\left(\frac{1}{2}\| \cB_h \left(D^{n+1}- D^n\right)\|_{L^2}^2 \,+\,
  \alpha_1\,\frac{\tau(\eps)}{\eps}\,\left\langle \cA_h \left(D_0^{n+1}-D_0^n\right),\, D_1^{n+1}-D_1^n \right\rangle\right)\,.
\end{equation}
On the one hand we estimate the terms $\cJ_1^{n+1}$ and $\cJ_2^{n+1}$
following the same method as the one presented to estimate their
continuous analogs $\cJ_1(t)$ and $\cJ_2(t)$ (see the proof item $(ii)$ in Theorem \ref{th:main1}).
On the other hand, the remainder term $\cR_1^{n+1}$ can be treated as $\cR_0^{n+1}$ in
the proof of $(i)$ of Theorem \ref{th:main1:h}. Indeed,
\begin{equation*}
\left|\left\langle \cA_h
  \left(D_0^{n+1}-D_0^n\right), D_1^{n+1}-D_1^n \right\rangle\right| \,\leq\, \frac{1}{2} \left(
          \|D_0^{n+1}-D_0^n\|_{L^2}^2 + \|\cA^*
          \left(D_1^{n+1}-D_1^{n}\right)\|_{L^2}^2
          \right)\,.
\end{equation*}
According to the mass conservation property \eqref{cons:mass:h}, $D_0^{n+1}-D_0^n$ meets condition \eqref{brac}. Therefore we may apply the discrete Poincar\'e inequality \eqref{ineg-disc} to bound $\|D_0^{n+1}-D_0^n\|_{L^2}^2$ in the latter estimate, this yields
\begin{equation*}
\left|\left\langle \cA_h
  \left(D_0^{n+1}-D_0^n\right), D_1^{n+1}-D_1^n \right\rangle\right| \,\leq\, \frac{1+C_d^2}{2}\,
  \|\cB_h
          \left(D^{n+1}-D^{n}\right)\|_{L^2}^2\,.
\end{equation*}
As in the case of $\cR_0^{n+1}$ in the former section, the latter estimate ensures that, as long as $\ds\alpha_0\,\leq\,
\left(
\overline{\tau}_0
\,
(1+C_d^2)
\right)^{-1}
$, it holds
\[
0\,\leq\,\cR_1^{n+1}\,.
\]
Hence, we obtain the result by adapting at the discrete level the proof
of item $(ii)$ in Theorem \ref{th:main1} to bound $\cJ_1^{n+1}$ and
$\cJ_2^{n+1}$ and applying a discrete Gronwall lemma.

\subsection{Proof of Theorem \ref{th:main2:h}}
As in the continuous setting, we prove that the solution $D^n=(D_{k}^n)_{k\in\N}$ to \eqref{discrete} converges to $D^n_{\tau_0}=(D^n_{\tau_0,k})_{k\in\N}$ given by \eqref{disc:eq:lim1}-\eqref{disc:eq:lim0}, whose long time behavior is easily obtained relying on the discrete Poincar\'e inequality \eqref{ineg-disc}
\begin{equation}\label{disc:cv:barD0}
    \|D^n_{\tau_0} -D_{\infty}\|_{L^2}
    \,\leq\,
    \|D^0_{\tau_0} -D_{\infty}\|_{L^2}\,
    \left(
1\,+\,
\frac{2\,\tau_0}{C_d^2}\, \Delta t
\right)^{-\frac{n}{2}}\,,\quad\forall\,t\in\R^+\,.
\end{equation}
We estimate $\ds \left\| D^n_0-D^n_{\tau_0,0}\right\|_{H^{-1}}$ by introducing the intermediate quantity $\cE$, meant to recover coercivity with respect to the first coefficient $D_0^n$
\begin{equation}\label{eq:E2:h}
	\cE^n\,=\,\frac{1}{2}\, \|\cA_h\,v_h^n\|_{L^2}^2\,,
\end{equation}
where $v_h^n$ solves \eqref{eq:elliptic:h} with source term
\[g\,=\, D_{0}^n\,
  +\,\frac{\tau(\eps)}{\eps}\,\cA^\star_h D^n_{1}\,-\, D_{\tau_0,0}^n\,.\]
The following lemma ensures that the quantity $\cE^n$ shares the same properties as its continuous analog. Indeed it holds 
\bl
We consider $\cE^n$ defined by \eqref{eq:E2:h}. It holds uniformly with respect to $\eps$
\begin{equation}
	\label{equiv3:h}
	\cE^n\,\leq\, \|D^n\,-\,D^n_{\tau_0}\|_{H^{-1}}^2
	\,+\,C_d^2\,\frac{\tau(\eps)^2}{\eps^2}\,
	\|\cB_h D^n\|_{L^2}^2\,,
\end{equation}
and
\begin{equation}
	\label{equiv4:h}
	\frac{1}{4}\,
	\|D^n\,-\,D^n_{\tau_0}\|_{H^{-1}}^2
	\,-\,C_d^2\,\frac{\tau(\eps)^2}{2\,\eps^2}\,
	\|\cB_h D^n\|_{L^2}^2
	\,\leq\,
	\cE^n
	\,.
\end{equation}
\label{lem:4:h}
\el 
\begin{proof}
  Defining $w^n_h$ and $u_{\tau_0}$ as the respective solutions to \eqref{eq:elliptic:h} with source term $g\,=\,\cA^\star_h D^n_1$ and $D_{\tau_0,0}\,-\,D_{\infty,0}$, it holds
	\[
	v^n_h\,=\,u^n_h\,-\,u^n_{\tau_0}\,+\,\frac{\tau(\eps)}{\eps}\,w^n_h\,.
	\]
	Applying operator $\cA_h$ to the latter relation, taking the $L^2$ norm, and applying the triangular inequality, it yields 
	\[\sqrt{2\,\cE^n}\,\leq\,
	\left\|
	\cA_h\left(u^n_h-u^n_{\tau_0}\right)
	\right\|_{L^2}
	\,+\,
	\frac{\tau(\eps)}{\eps}\,
	\left\|
	\cA_h\,w^n_h
	\right\|_{L^2}\,,
	\]
	and
	\[
	\left\|
	\cA_h\left(u^n_h-u^n_{\tau_0}\right)
	\right\|_{L^2}
	\,-\,
	\frac{\tau(\eps)}{\eps}\,
	\left\|
	\cA_h\,w^n_h
	\right\|_{L^2}\,\leq\,\sqrt{2\,\cE^n}\,.
	\]
	We estimate 
	$
	\left\|
	\cA_h\,w^n_h
	\right\|_{L^2}
	$ applying \eqref{lem:3.30} in Lemma \ref{lem:3.3}, this yields
	\[
	\sqrt{2\,\cE^n}\,\leq\,
		\|D^n\,-\,D^n_{\tau_0}\|_{H^{-1}}
	\,+\,
	\frac{\tau(\eps)}{\eps}\,
	C_d\,\|\cB_h D^n\|_{L^2}\,,
	\]
	and
	\[
		\|D^n\,-\,D^n_{\tau_0}\|_{H^{-1}}
	\,-\,
	\frac{\tau(\eps)}{\eps}\,
	C_d\,\|\cB_h D^n\|_{L^2}
	\,\leq\,\sqrt{2\,\cE^n}\,.
	\]
	We obtain the result taking the square of the latter inequalities and applying Young's inequality.
\end{proof}

We now treat the asymptotic limit $\eps \rightarrow 0$ corresponding
to the case of $(i)$ in Theorem \ref{th:main2:h} and therefore suppose
that $\tau(\eps)$ fulfills the  assumptions \eqref{cond:tau}, \eqref{cond2a:tau} and \eqref{cond3:tau}. As in the continuous setting, we start by deriving the first result in $(i)$ of Theorem \ref{th:main2:h}. We  already know from the $L^2$  estimate \eqref{estim:L2:discrete} that
	\begin{eqnarray*}
		\ds\frac{\left\|D_{\perp}^{n+1}\right\|^2_{L^2}-\left\|D_{\perp}^n\right\|^2_{L^2}}{2\,\Delta t}\,&+&
		\frac{1}{\tau(\eps)}\, \left\|D_{\perp}^{n+1}\right\|^2_{L^2}\\
		\ds&\leq &-\,
		\left\langle \frac{D_0^{n+1}-D_0^{n}}{\Delta t},\,
		D_0^{n+1}-D_0^{n}\right\rangle
		\,-\,
		\frac{1}{2\,\Delta t}\sum_{k \in \N^*}
		\| D^{n+1}_{k}- D^n_{k}\|_{L^2}^2\\
		\ds&\leq& -\,
		\left\langle \frac{D_0^{n+1}-D_0^{n}}{\Delta t},\,
		D_0^{n+1}-D_{\infty,0}\right\rangle\,.
	\end{eqnarray*}
Therefore, we replace $D_0^{n+1}-D_0^{n}$ according to equation
\eqref{discrete}, and after applying the duality formula of  Lemma \ref{lem:3.1}-\eqref{prop1:A:h}, we obtain
	\begin{equation*}
	\frac{\left\|D_{\perp}^{n+1}\right\|^2_{L^2}-\left\|D_{\perp}^n\right\|^2_{L^2}}{\Delta t}\, +
	\frac{1}{\tau(\eps)}\, \left\|D_{\perp}^{n+1}\right\|^2_{L^2}
	\,\leq\, -\frac{1}{\eps} \left\langle D_1^{n+1},\,
	\cA_h D_0^{n+1}\right\rangle\,, 
\end{equation*}
	Hence, after multiplying by $\Delta t$ and applying the Young inequality to bound the right hand side of the latter inequality, it yields
	$$
	\left(1+\frac{\Delta t}{\tau(\eps)}\right)\left\|D_{\perp}^{n+1}\right\|^2_{L^2}
	\,\leq\,\left\|D_{\perp}^n\right\|^2_{L^2}\,+\,\Delta t\,\frac{\tau(\eps)}{\eps^2}\,\|\cB_{h} D^{n+1}\|_{L^2}^2\,. 
	$$
	To achieve the proof, it remains to bound  $\|\cB_h
        D^{n+1}\|_{L^2}^2$ by applying  Theorem
        \ref{th:main1:h}-$(ii)$ and again following the line of the
        proof of Theorem \ref{th:main2}, 
	we deduce
	\begin{align*}
	&\left\|D_{\perp}^n\right\|^2_{L^2} \,\leq\\ &\left\|D_{\perp}^0\right\|^2_{L^2} \left(1+\frac{\Delta t}{\tau(\eps)}\right)^{-n}\,+\, 6\,\left(C(\overline{\tau}_0^2+1)\,\| D^{0}-D_{\infty}\|_{L^2}^2\,
	+
	\|\cB_h D^{0}\|_{L^2}^2
	\right)
	\,\frac{\tau(\eps)^2}{\eps^2}\,
	\left(1+\frac{\tau(\eps)}{\eps^2}\kappa\Delta t\right)^{-n}\,. 
	\end{align*}
Therefore we obtain the result taking the square root in the latter estimate and substituting $\tau(\eps)$ with $\tau_0\,\eps^2$ according to assumption \eqref{cond3:tau}.\\

To prove the second result of $(i)$ in Theorem \ref{th:main2:h} we evaluate $\cE^n$ as in the proof of
Theorem \ref{th:main2} observing that
$$
\|\cA_h\,v^n_h\|_{L^2}^2 = \left\langle D_0^n\,
+\,\frac{\tau(\eps)}{\eps}\,\cA^\star_h D^n_1
\,-\,D^{n}_{\tau_0,0},\,v^n_h\right\rangle
$$
hence, relying on equations \eqref{discrete} and \eqref{disc:eq:lim0} we deduce 
\begin{equation*}
\frac{\cE^{n+1}-\cE^{n} }{\Delta t} \,=\, -\,\frac{\tau(\eps)}{\eps^2}
\, \| D_0^n
+\frac{\tau(\eps)}{\eps}\,\cA^\star_h D^n_1
-D^n_{\tau_0}\|_{L^2}^2
\,+\,
\cE_{1}^{n+1}
\,+\,
\cE_{2}^{n+1}
\,+\,
\cE_{3}^{n+1}\,-\,\cR^{n+1}_3\,,
\end{equation*} 
where $\ds\cE_{1}^{n+1}$, $
\ds\cE_{2}^{n+1}
$ and $
\ds\cE_{3}^{n+1}$ are the numerical equivalents of the terms $\ds\cE_{1}(t)$, $
\ds\cE_{2}(t)
$ and $
\ds\cE_{3}(t)$ in the proof of Theorem \ref{th:main2}
\begin{equation*}
\left\{
\begin{array}{l}
\ds \cE_{1}^{n+1}
\,=\,
\left(
\tau_0-\frac{\tau(\eps)}{\eps^2}
\right)\,
\left\langle 
\cA^*_h\cA_h\,D_{\tau_0,0}^{n+1}\,,\,
v^{n+1}_h\right\rangle\,,
\\[1.1em]
\ds\cE_{2}^{n+1}
\,=\,\frac{\tau(\eps)^2}{\eps^3}\,
\left\langle 
\cA^*_h\cA_h\, D^{n+1}_1\,,\,
v^{n+1}_h\right\rangle\,,
\\[1.1em]
\ds\cE_{3}^{n+1}
\,=\,\sqrt{2}\,
\frac{\tau(\eps)}{\eps^2}\,
\left\langle 
\left(\cA^*_h\right)^2
D_2^{n+1}\,,\,
v^{n+1}_h\right\rangle\,,
\end{array}\right.
\end{equation*} 
and $\cR^{n+1}_3$ is a numerical dissipation term
\[\cR^{n+1}_3\,=\,\frac{1}{2\Delta t}\,\left\| \cA_h\left( v^{n+1}_h - v^{n}_h \right) \right\|_{L^2}^2\,.\]
Since $\cR^{n+1}_3$ is positive, we apply the same method as the one presented in the proof of Theorem \ref{th:main2}  and therefore we obtain the following estimate for $\cE^{n}$
\begin{align*}
	\left(1+\frac{\tau(\eps)\Delta t}{C_d^2\,\eps^2}\right)
	\cE^{n+1}
	\,\leq
	\cE^{n}&+
	\,C\,\Delta t\,
	\frac{\tau(\eps)}{\,\eps^2}
	\,\left(1+\overline{\tau}_0^2\right)\left\|D_{\perp}^{n+1}\right\|^2_{L^2}\\[0.8em]
	&+C\,\Delta t\,
	\frac{\eps^2}{\tau(\eps)}
	\left|
	\tau_0     
	-
	\frac{\tau(\eps)}{\eps^2}
	\right|^2\,\|D^{n+1}_{\tau_0}-D_{\infty}\|_{L^2}^2\,,
\end{align*} 
	for some constant $C$ depending only on $\Phi$ and $T_0$. 
	In the latter inequality, we bound $\|D^{n+1}_{\tau_0}-D_{\infty}\|_{L^2}^2$ according to \eqref{disc:cv:barD0} and the norm of $\ds D_{\perp}$ according to the first estimate of $(i)$ in Theorem \ref{th:main2:h}. Then we multiply the inequality by $\ds\left(1+\frac{\tau(\eps)\Delta t}{C_d^2\,\eps^2}\right)^{n}$ and sum for $k$ ranging from $0$ to $n-1$, it yields
	\begin{align*}
		\cE^{n} \,&\leq \, 
		\left(
		\cE^0
		\,+\,C\,\frac{\tau(\eps)^2}{\eps^2}\,(\overline{\tau}_0^6+1)\| D^{0}-D_{\infty}\|_{H^1}^2\right)
		\left(1+\frac{\tau(\eps)}{\eps^2}\,\kappa\,\Delta t\right)^{-n}\\
		&+
		C\,\left|
		\frac{\tau_0\,\eps^2}{\tau(\eps)}
		-
		1
		\right|^2\,\|D^{0}_{\tau_0}-D_{\infty}\|_{L^2}^2\,
		\left(
		\frac{2\,\tau_0\,\eps^2}{\tau(\eps)}
		-
		1
		\right)^{-1}
		\left(1+\frac{\tau(\eps)}{\eps^2}\,\kappa\,\Delta t\right)^{-n}
		\,. 
	\end{align*}
To conclude, we substitute $\cE^{n}$ (resp. $\cE^0$) in the latter estimate according to \eqref{equiv4:h} (resp. \eqref{equiv3:h}) in Lemma \ref{lem:4} and then apply assumption \eqref{cond3:tau} on $\tau(\eps)$, which ensures $\ds\left(
\frac{2\,\tau_0\,\eps^2}{\tau(\eps)}
-
1
\right)^{-1}\,\leq\, 3$, this yields 
	\begin{align*}
		\|D^n_0 -D^n_{\tau_0,0}\|_{H^{-1}}^2 \,\leq \, 
		&C\left(
		\|D^0_0 -D^0_{\tau_0,0}\|_{H^{-1}}^2
		\,+\,\frac{\tau(\eps)^2}{\eps^2}\,(\overline{\tau}_0^6+1)\| D^{0}-D_{\infty}\|_{H^1}^2\right)
		\left(1+\frac{\tau(\eps)}{\eps^2}\,\kappa\,\Delta t\right)^{-n}\\
		+\,
		&C\left|
		\frac{\tau_0\,\eps^2}{\tau(\eps)}
		-
		1
		\right|^2\|D^{0}_{\tau_0}-D_{\infty}\|_{L^2}^2\,
		\left(1+\frac{\tau(\eps)}{\eps^2}\,\kappa\,\Delta t\right)^{-n}
		\,. 
	\end{align*}
	We obtain the result taking the square root in the latter estimate and substituting $\tau(\eps)$ with $\tau_0\,\eps^2$according to assumption \eqref{cond3:tau}.
        \\
        Finally the proof of the second item follows the same lines
        replacing $D_{\tau_0}^n$ by $D_\infty$ in the discrete
        functional $\cE^n$. 

%

\section{Numerical simulations}
\label{sec:4}
\setcounter{equation}{0}
\setcounter{figure}{0}
\setcounter{table}{0}

%

We performed several numerical simulations which confirm the accuracy of
the scheme \eqref{discrete}. We do not detail this process here and rather focus on the physical interpretation and the quantitative results obtained in our experiments. We refer to \cite{bessemoulin2022cv} for a precise discussion on that matter.

In this section, we want to illustrate the quantitative estimates of the
solution obtained using  the Hermite Spectral method in velocity and
finite volume scheme in space for the one-dimensional
Vlasov-Fokker-Planck equation.  We choose $\tau(\eps)=\tau_0\,\eps^2$ with $\tau_0=5$ and
consider the Vlasov-Fokker-Planck equation \eqref{vfp0} with $E=-\partial_{\bx}\Phi$ and
$$
\Phi(x) = 0.1\, \cos\left(\frac{2\pi \,x}{L}\right) \,+\, 0.9\, \cos\left(\frac{4\pi \,x}{L}\right), 
$$
The stationary state is given by the Maxwell-Boltzmann distribution
$$
f_\infty(x,v) = \frac{c_0}{\sqrt{2\pi}}\, \exp\left(-\left(\Phi+\frac{|v|^2}{2}\right)\right),
$$
where $c_0$ is given by mass conservation
$$
\int_{\T\times\R} f_\infty \dD v\dD x \,=\,  \int_{\T\times\R} f_0 \dD
v\dD x,\,
$$
where $f_0$ is the initial datum.

In our simulation, we take a time step $\Delta t= 10^{-3}$, a number
of modes $N_H= 200$ and $N_x=64$.  It is worth to mention that all the
numerical simulations presented in this section are not affected by
the numerical parameters, which allows us to focus our discussion on the quantitative results on the diffusive limit $\eps \rightarrow 0 $ and
large time behavior.

\subsection{Test 1 : centred Maxwellian}
\label{sec:4.2}

\begin{figure}
  \centering
  \begin{tabular}{cc}
    \includegraphics[width=3.2in,clip]{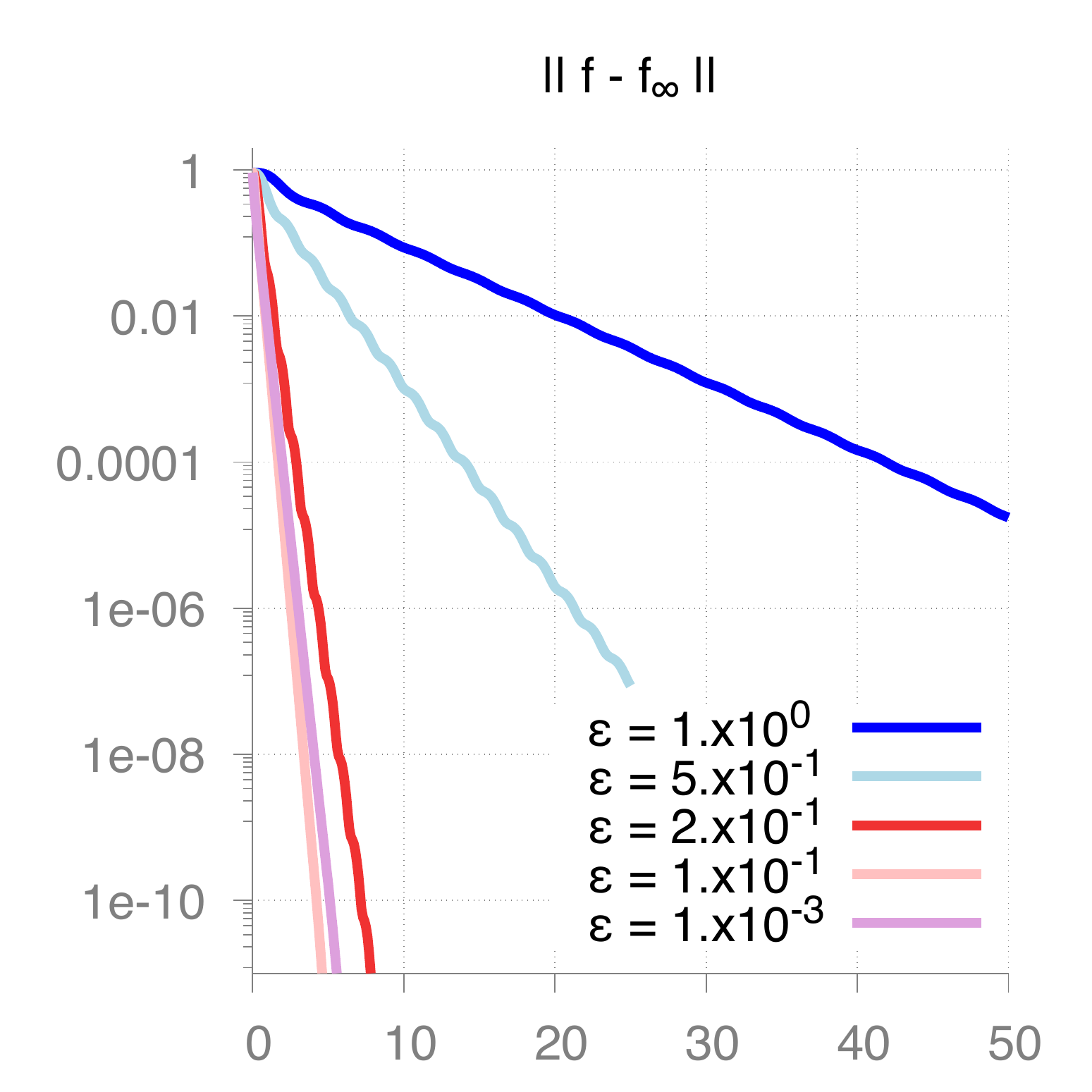}&
                                                     \includegraphics[width=3.2in,clip]{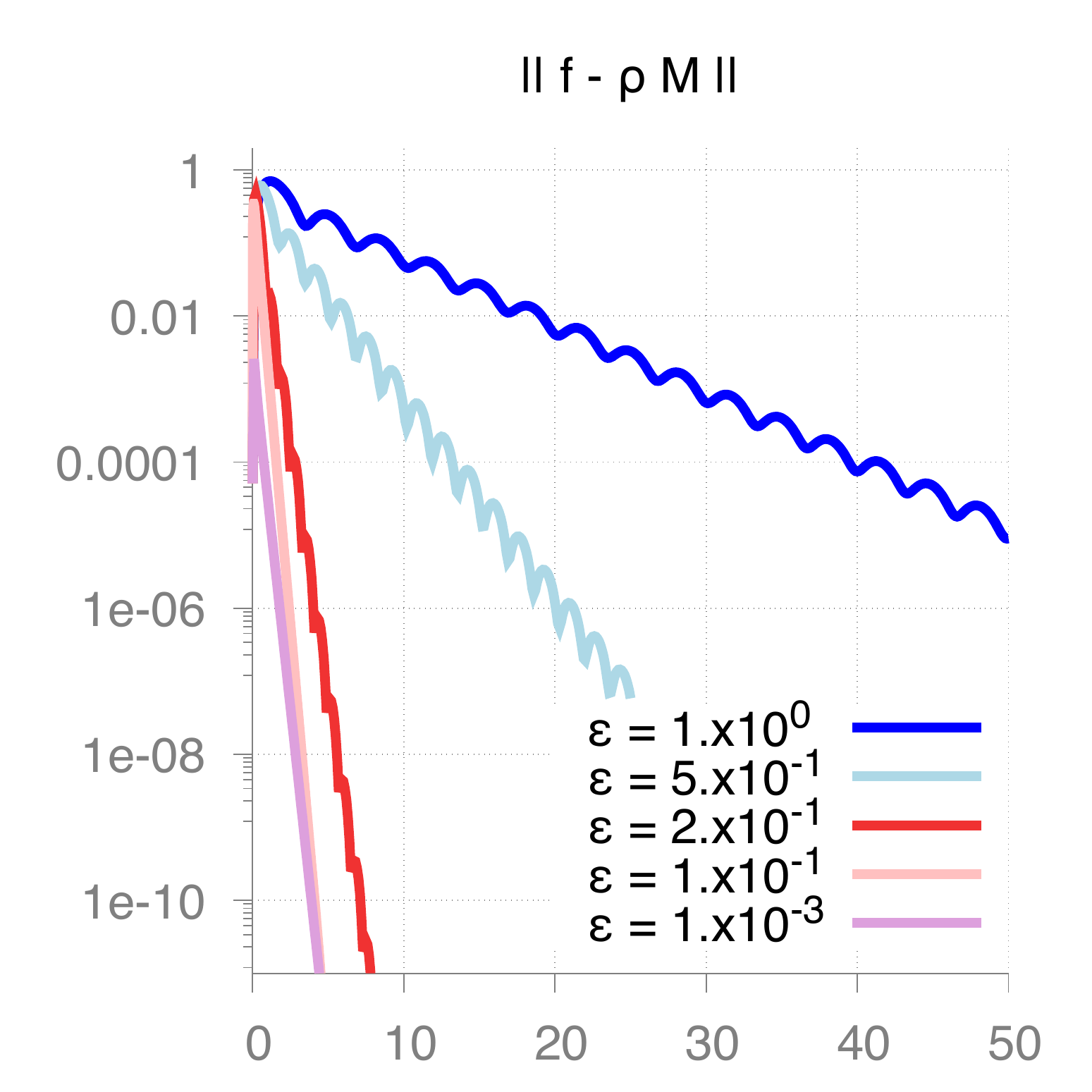}
    \\
    $(a)$ & $(b)$
        \end{tabular}
	\caption{{\bf Test 1 : centred Maxwellian.}  time evolution
          in log scale of $(a)$ $\|f-f_\infty\|_{L^2(f_\infty^{-1})}$,
            $(b)$ $\|f-\rho \, \cM \|_{L^2(f_\infty^{-1})}$. }
	\label{fig:110}
\end{figure}

For the first test, we choose the following initial condition
$$
f_0(x,v) \,=\, \frac{1}{\sqrt{2\pi}}\,\left(1+ \delta \cos\left(\frac{2\pi \,x}{L}\right)\right)\,\exp\left(-\frac{|v|^2}{2}\right)\,,
$$
with $\delta=0.5$ and $L\,=\,10$.\\
On the one hand, we present in Figure \ref{fig:110}  the time
evolution of $\|f-f_\infty\|_{L^2(f_\infty^{-1})} $ and the relative entropy on
$f$,
$$
\| f - \rho \,\cM \|_{L^2(f_\infty^{-1})}  \,=\,   \| D_\perp(t) \|_{L^2}.
$$

The most striking feature in this test consists in the oscillatory
behavior of the relative entropy which unfolds in the relaxation of
$f$ towards its equilibrium. These oscillations may be observed in
Figure \ref{fig:110}-$(b)$  and occur for
various values of $\eps$ ranging from $1$ represented by blue curves
to $2.10^{-1}$ represented by red curves.

\begin{figure}
  \centering
  \begin{tabular}{cc}
      \includegraphics[width=3.2in,clip]{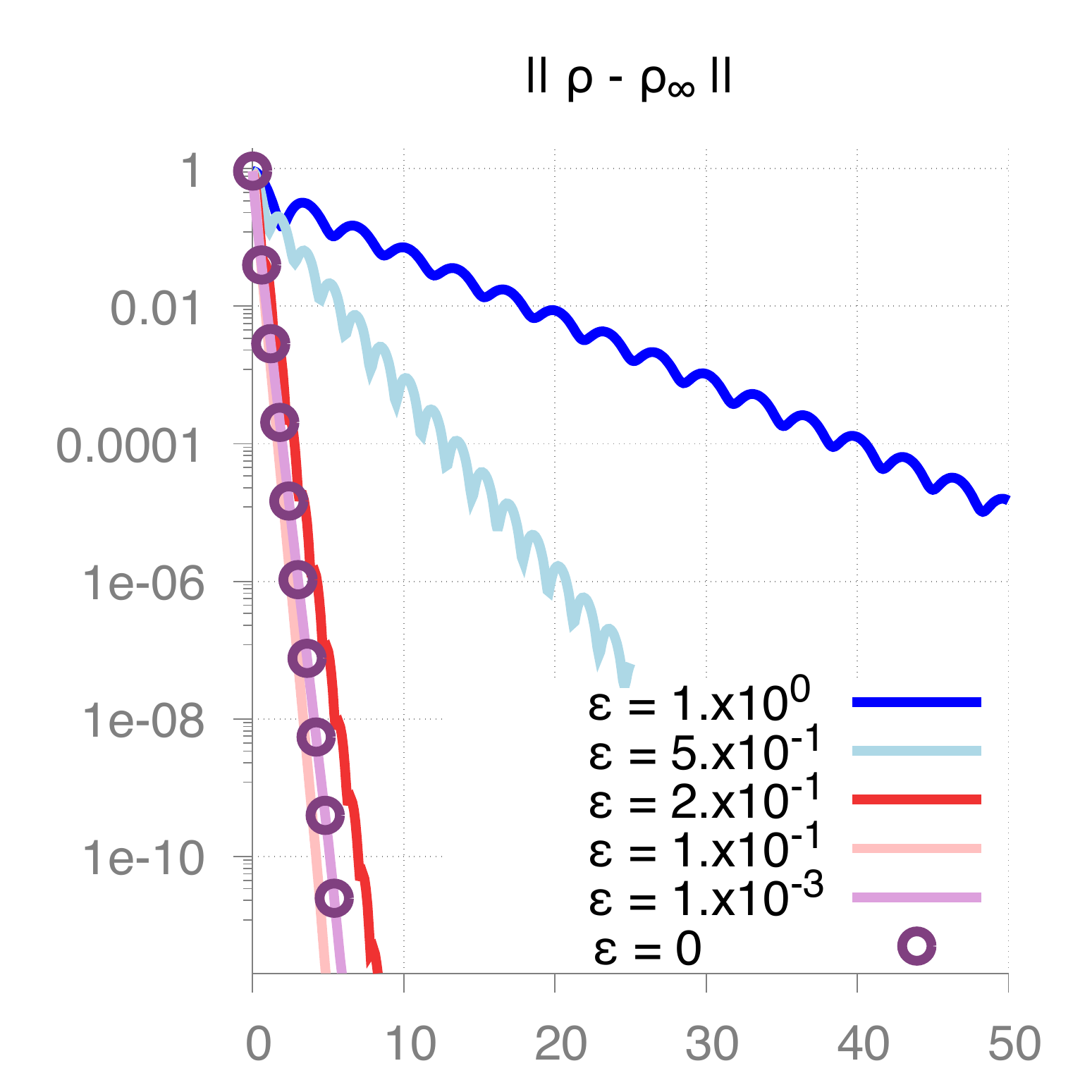}&
                                                     \includegraphics[width=3.2in,clip]{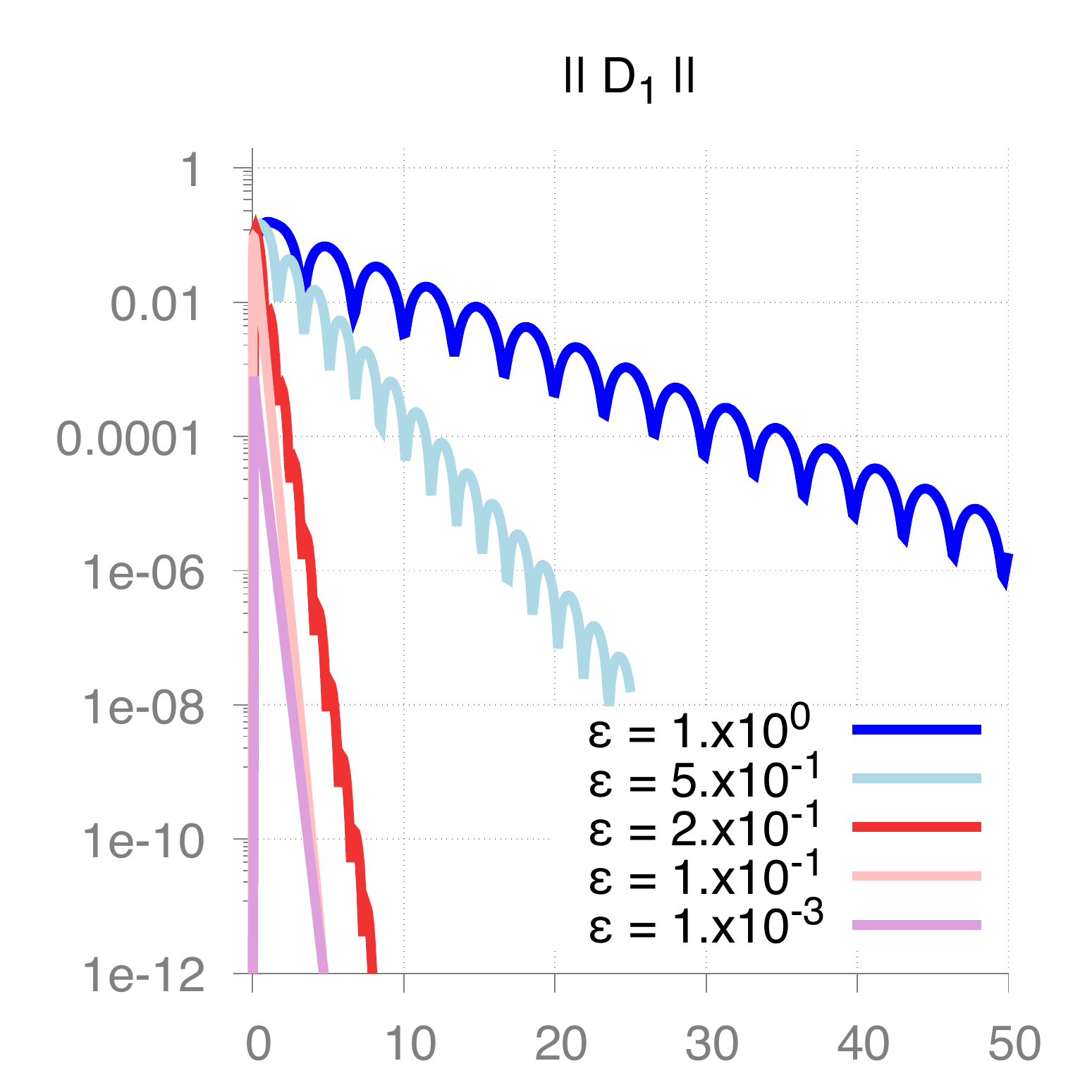}
    \\
    (a)&(b) \\
        \end{tabular}
	\caption{{\bf Test 1 : centred Maxwellian.}  time evolution
          in log scale of $(a)$
              $\|\rho-\rho_\infty  \|_{L^2(\rho_\infty^{-1})}$ and  $(b)$ $\|D_1  \|_{L^2}$. }
	\label{fig:111}
\end{figure}

We also present in Figure \ref{fig:111} the relaxation to equilibrium
of  macroscopic quantities 
$$
\|D_0 - D_{\infty,0}\|_{L^2} = \| \rho- \rho_{\infty}\|_{L^2(f_\infty^{-1})}
$$
and the norm of the first moment $D_1$. Time oscillations, observed on the
distribution function,  seem to affect  macroscopic quantities
associated to the solution as moments $D_0$ and $D_1$.

\begin{figure}
  \centering
  \begin{tabular}{cc}
    \includegraphics[width=3.2in,clip]{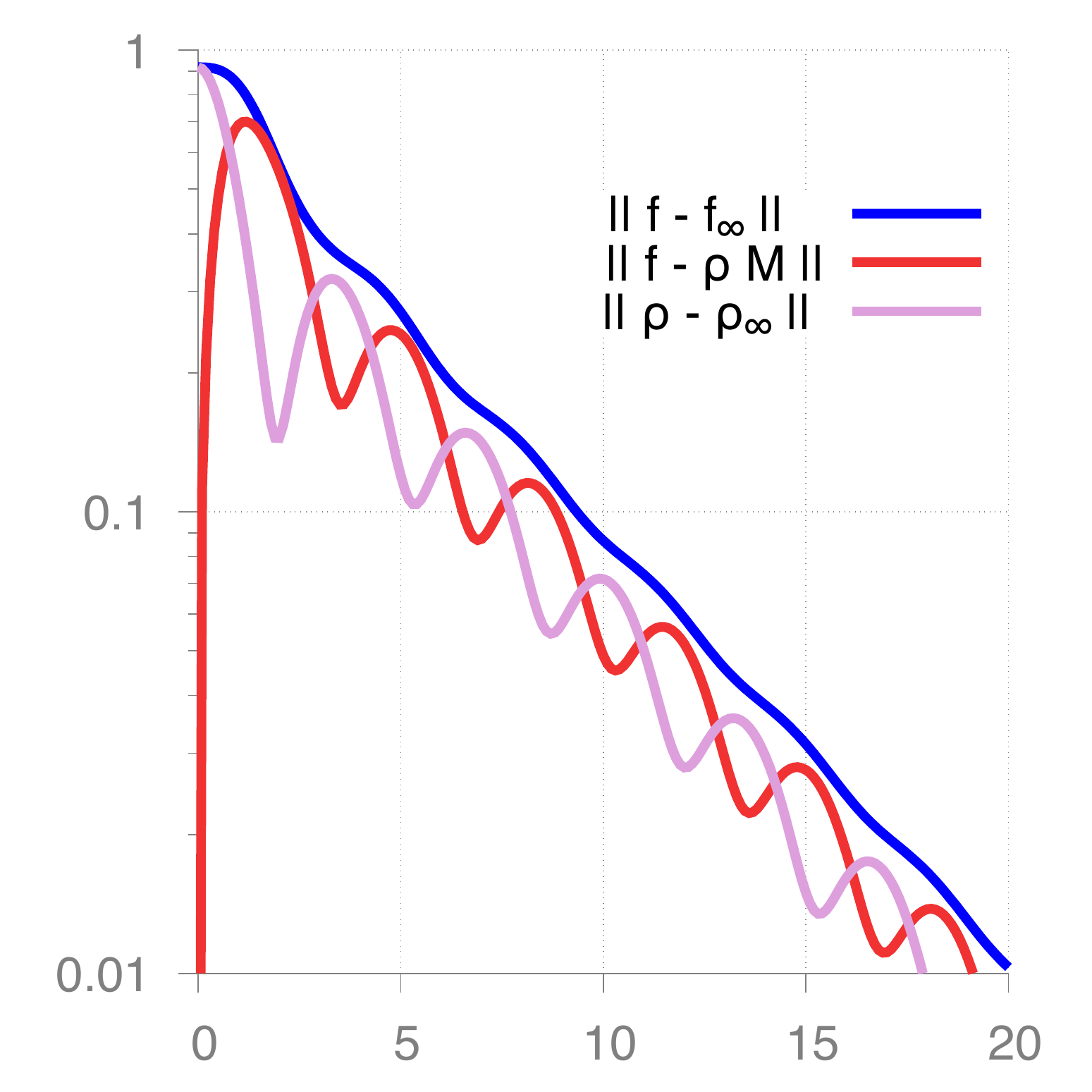}&
                                                      \includegraphics[width=3.2in,clip]{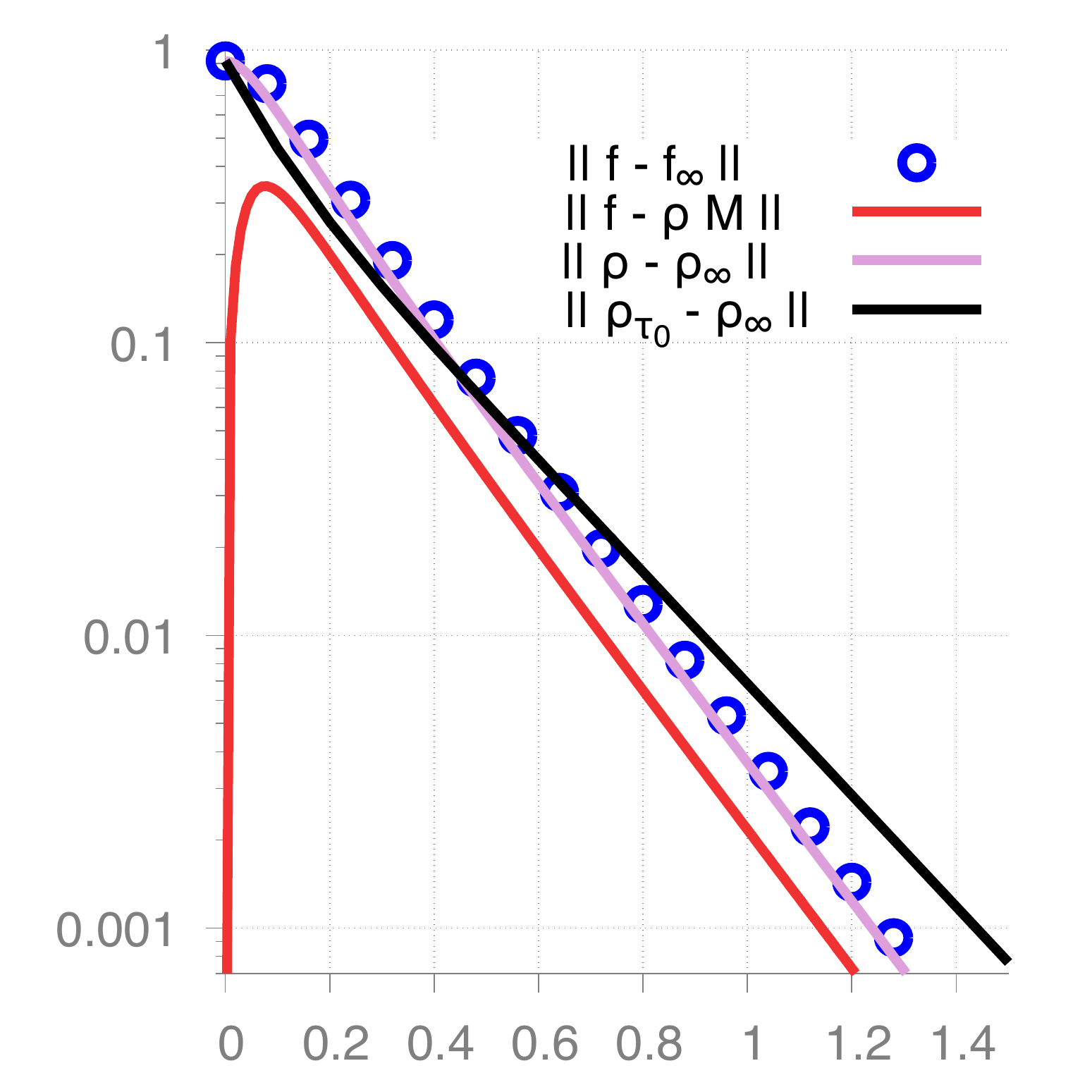}\\
    $(a)$ $\eps=1$ & $(b)$  $\eps=10^{-1}$
  \end{tabular}
  \caption{{\bf Test 1 : centred Maxwellian.}  time evolution
          in log scale of  $\|f-f_\infty\|_{L^2(f_\infty^{-1})}$ (blue),
             $\|f-\rho \, \cM \|_{L^2(f_\infty^{-1})}$ (red),
              $\|\rho-\rho_\infty  \|_{L^2(\rho_\infty^{-1})}$ (pink)
              and $\|\rho_{\tau_0}-\rho_\infty  \|_{L^2(\rho_\infty^{-1})}$
              (black)  for
              $(a)$ $\eps=1$ and $(b)$ $\eps=10^{-1}$. }
            	\label{fig:12}
              \end{figure}

On the other hand, we provide In Figure \ref{fig:12},  a detailed description in
the case $\eps=1$, where we see that the oscillations of the spatial
density and the ones of the higher modes in velocity are asynchronous,
this may be interpretated as a transfer of information between these
two quantities. This phenomenon has already been investigated for
non-linear kinetic models (see  \cite{filbet2006}) but we show through these experiments that
even the simple model at play here captures this phenomena.

These oscillations stay visible for surprisingly small values of
$\eps$, up to $10^{-1}$. It showcases the robustness of our scheme,
which is still able to capture them at low computational cost. To be noted that our numerical experiments indicate that a non zero
external force field seems to be mandatory  to observe this
oscillatory behavior. We also emphasize that these oscillations seem
to be quite sensitive to the choice of the initial data and the external
field (see the second numerical test with a different initial
data, where such oscillations disappear for large time).

This leads us to the second feature of this test, which is the
asymptotic preserving property of the scheme for various values of $\eps$. The method is
accurate on large time intervals in the situation where $\eps\,=\,1$
(see Figure \ref{fig:12}-$(a)$), which corresponds to the long time
behavior of the model but it is also accurate when  $\eps\ll
1$. Indeed, as it is shown in Figure \ref{fig:111}-$(a)$,  the purple
error curve of the density $\rho$ corresponds exactly to the circled
error curve of the macroscopic model $\rho_{\tau_0}$ when
$\eps=10^{-3}$ and even smaller (not shown since the curves coincide).

Finally we focus on the intermediate value $\eps=10^{-1}$, for which
we observe in Figures \ref{fig:110}-$(a)$, \ref{fig:111}-$(a)$ and
\ref{fig:12}-$(b)$, a somehow surprising phenomenon: the kinetic model
relaxes faster towards equilibrium than the macroscopic one. This
appears to be a consequence of our choice of initial data which is
already at local equilibrium at time $t=0$. This aspect of the
experiment justifies our efforts to cover a wide range of values for
the scaling parameter $\eps$: it enables to capture intermediate
regimes which may display peculiar phenomena. As we will see in the
next section, the reverse situation is possible as well, when the
initial condition is far from equilibrium.

\begin{figure}
  \centering
  \begin{tabular}{cc}
    \includegraphics[width=3.2in,clip]{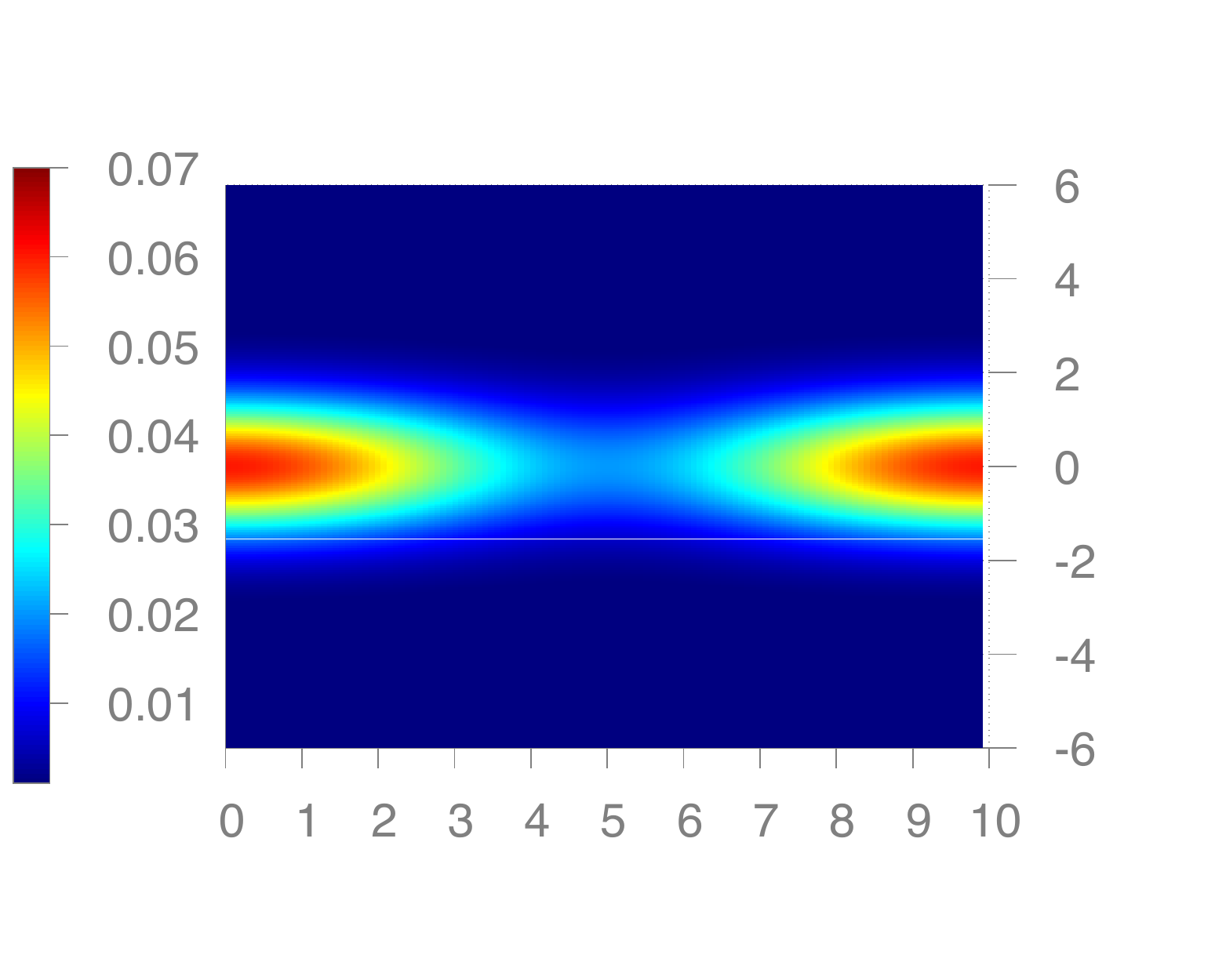}&
                                                      \includegraphics[width=3.2in,clip]{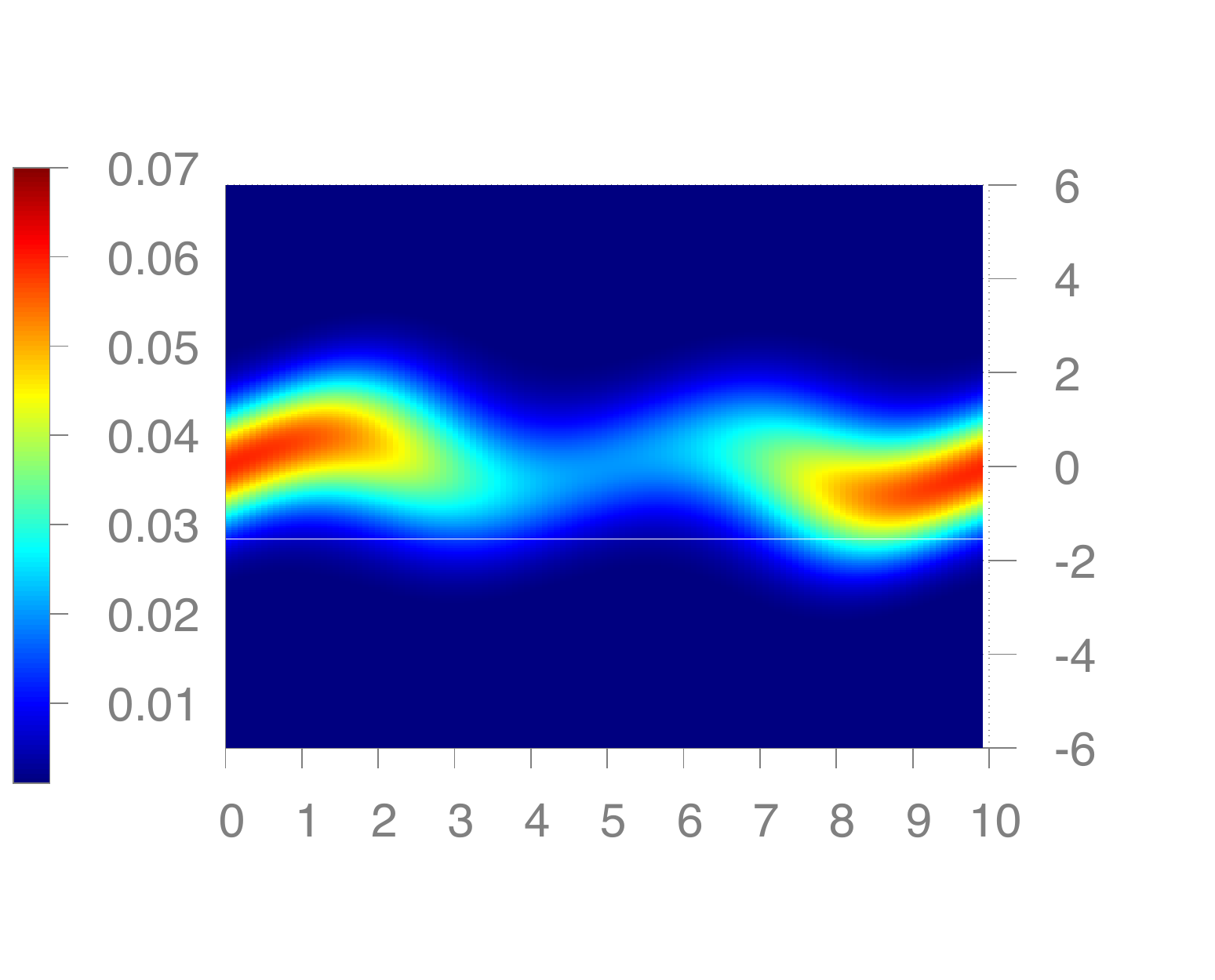}\\
    $(a)$ $t=0$& $(b)$  $t=0.5$ \\
    \includegraphics[width=3.2in,clip]{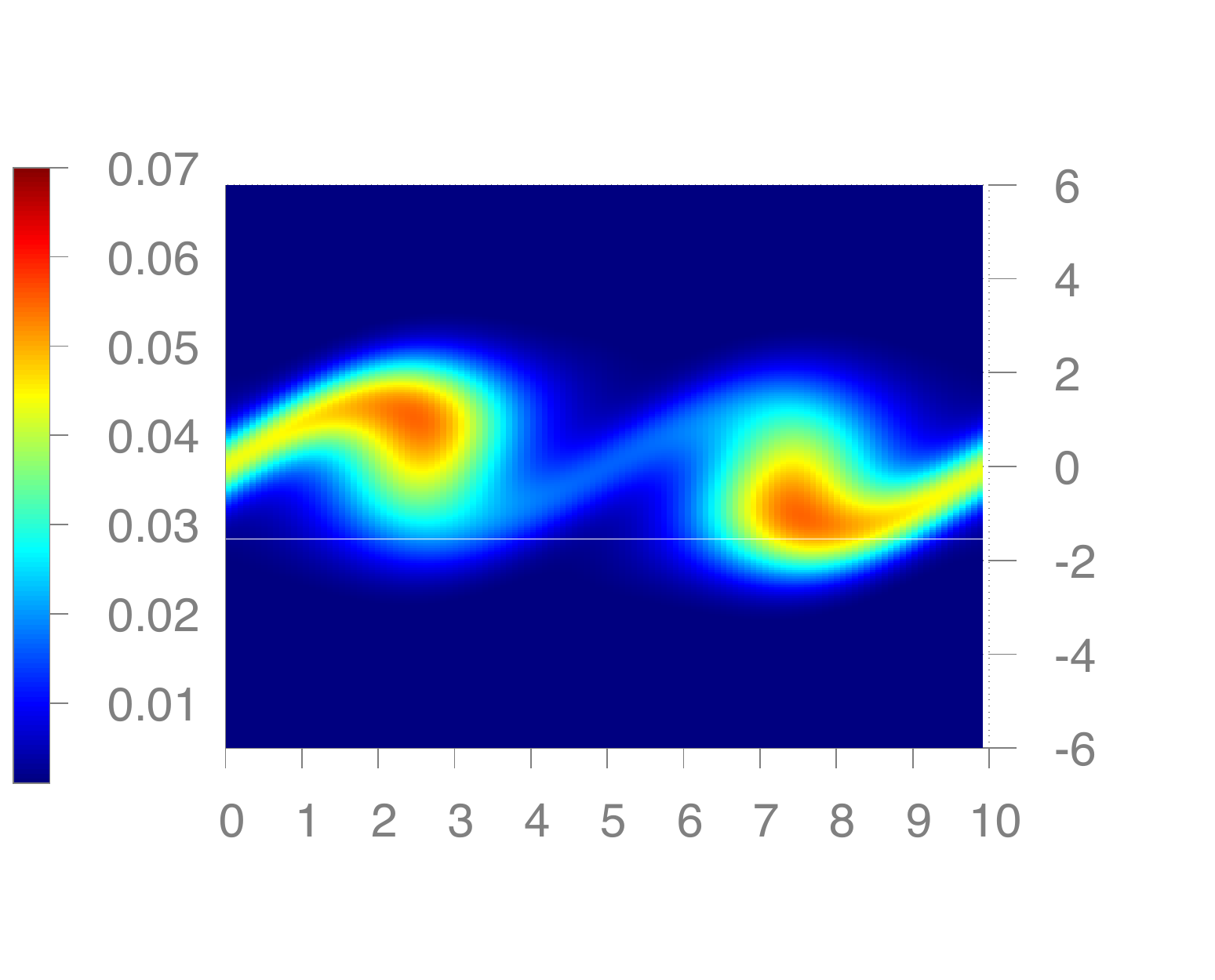}&
                                                      \includegraphics[width=3.2in,clip]{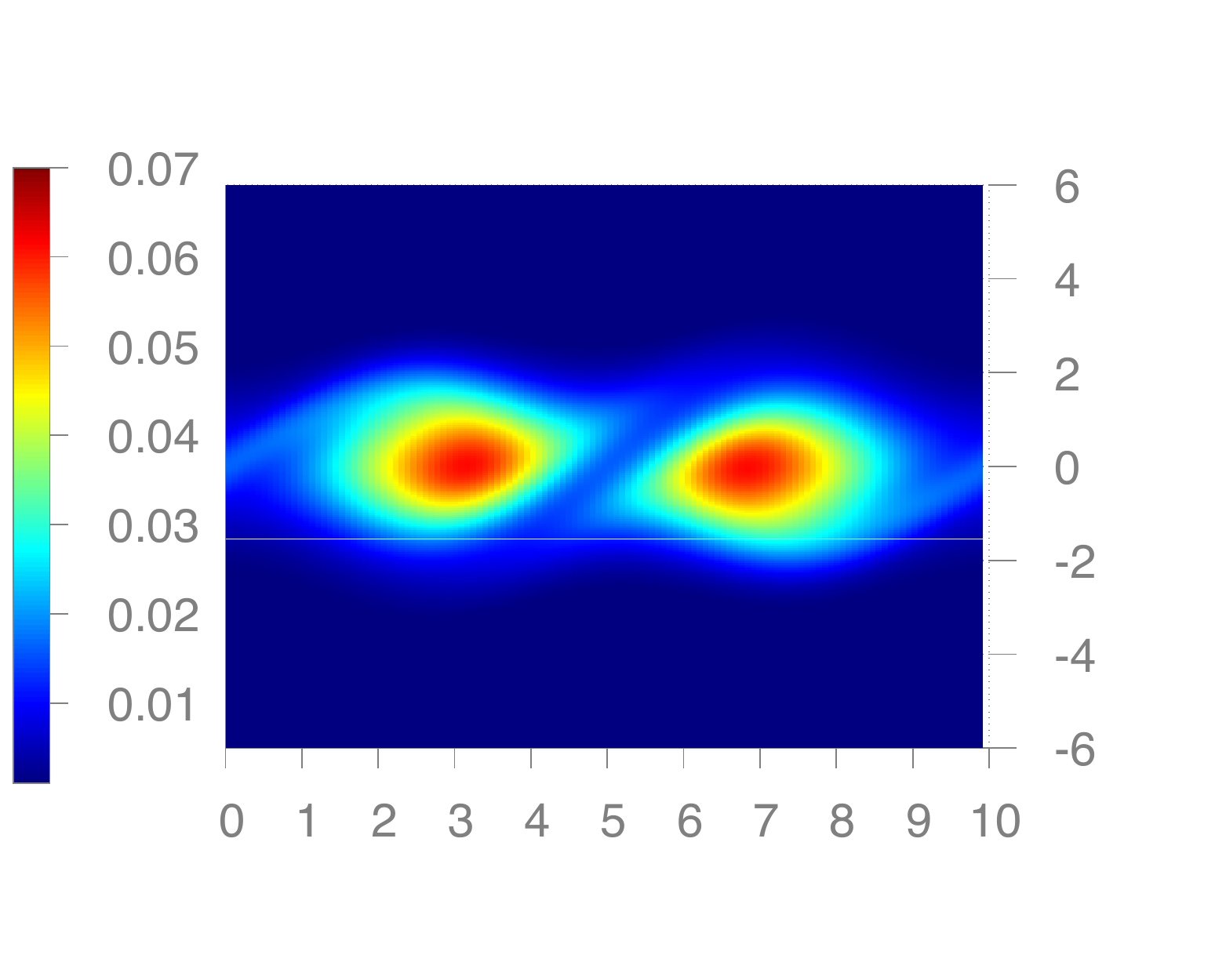}\\
    $(c)$ $t=1.5$ & $(d)$ $t=3$\\
         \includegraphics[width=3.2in,clip]{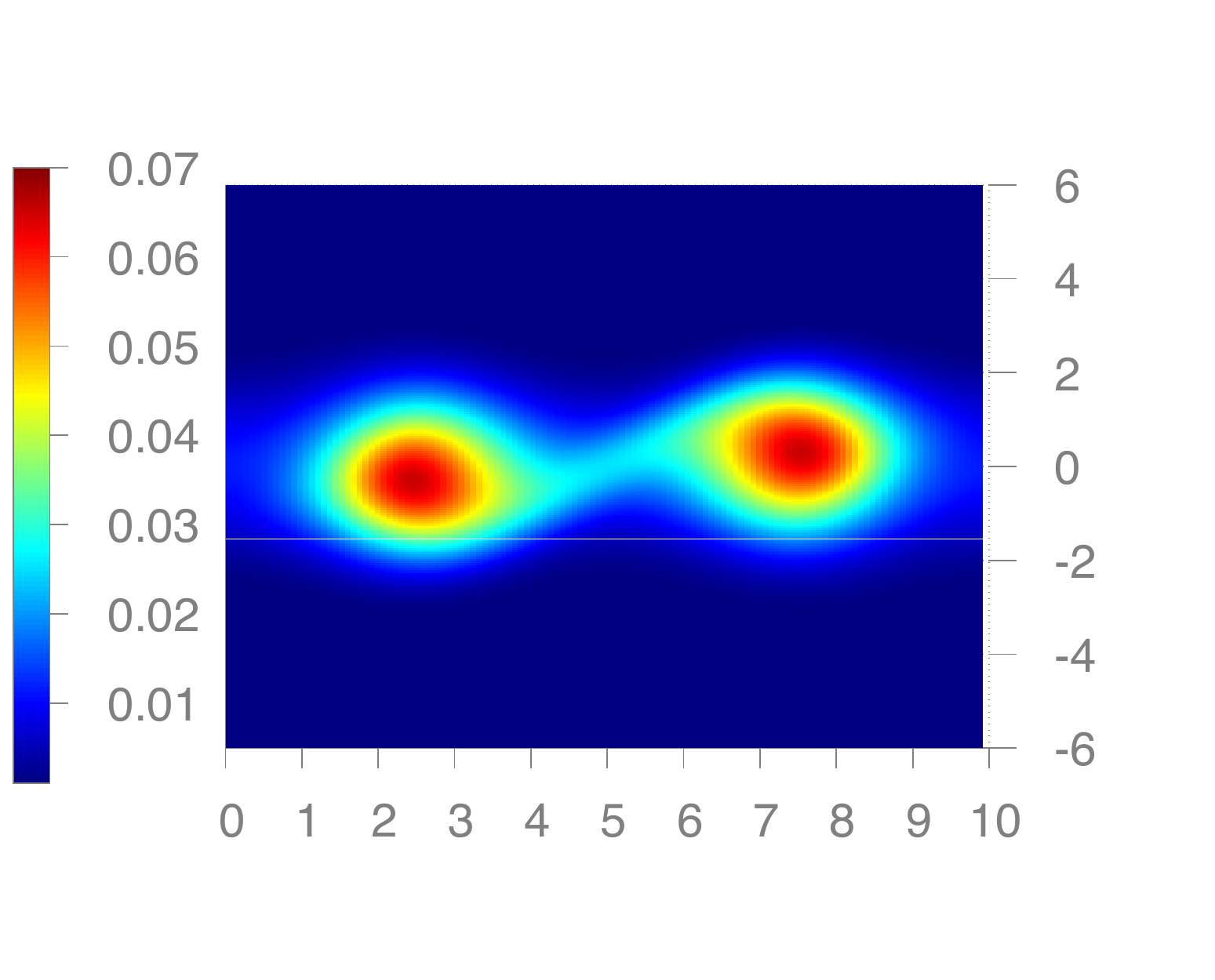}&
                                                      \includegraphics[width=3.2in,clip]{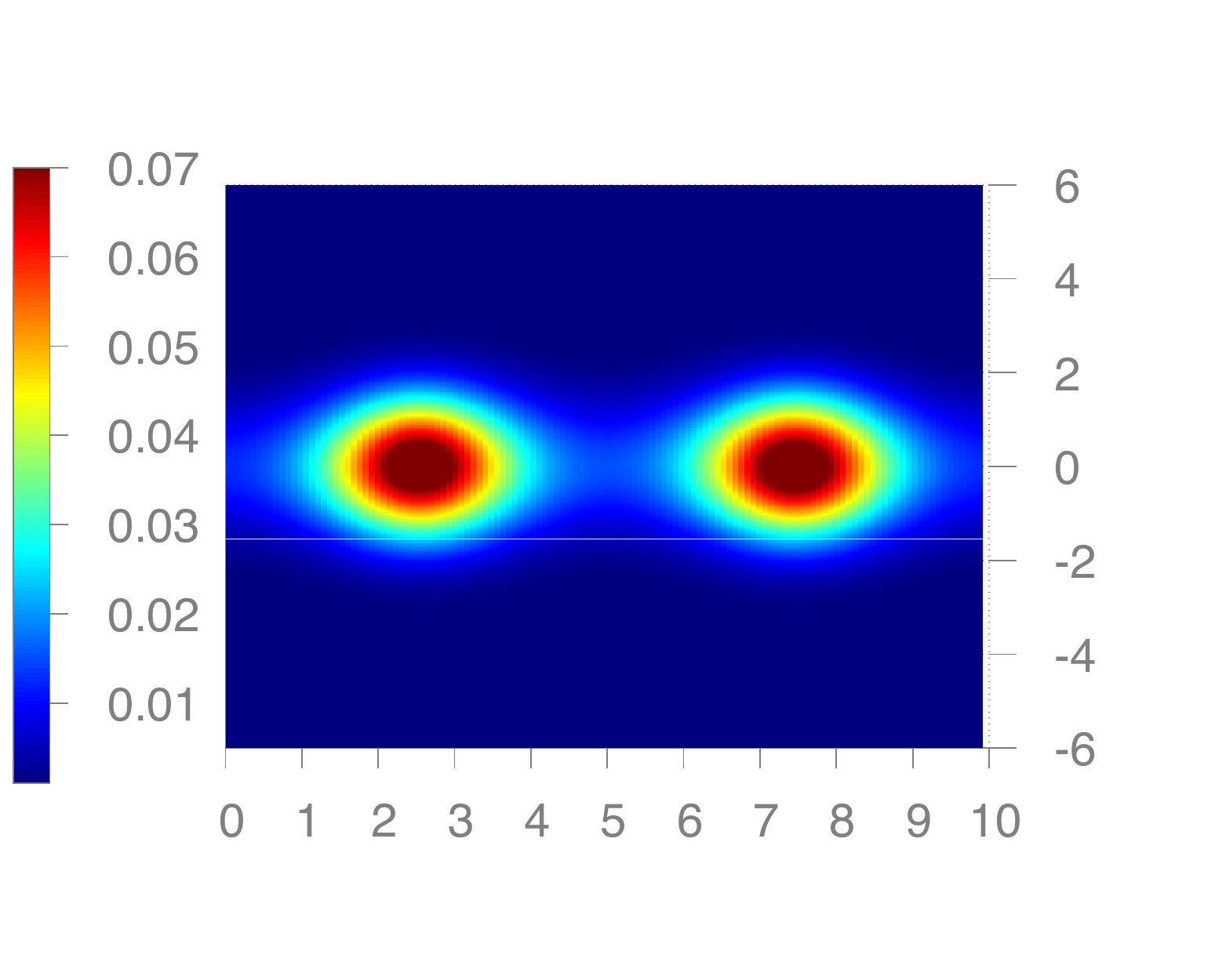}\\
    $(e)$ $t=5$ & $(f)$ $t=20$
  \end{tabular}
  \caption{{\bf Test 1 : centred Maxwellian.}  snapshots of the
    distribution function for $\eps=1$ at time $t=0, \, 0.5, \, 1.5,
    \, 3,\, 5$ and $20$. }
	\label{fig:13}
\end{figure}

We conclude this section by drawing the readers attention towards Figure \ref{fig:13}, which features the graph of the solution $f$ at different times, in the case $\eps=1$ and on which we witness its intricate relaxation towards equilibrium.

\subsection{Test 2 : shifted Maxwellian}
\label{sec:4.1}
We now choose the same parameter as before excepted that the initial
condition is a shifted Maxwellian 
$$
f_0(x,v) \,=\, \frac{1}{\sqrt{2\pi}}\,\left(1+ \delta \cos\left(\frac{2\pi \,x}{L}\right)\right)\,\exp\left(-\frac{|v-u_0|^2}{2}\right),
$$
with $u_0=1$, which is far from equilibrium.

\begin{figure}
  \centering
  \begin{tabular}{cc}
	\includegraphics[width=3.2in,clip]{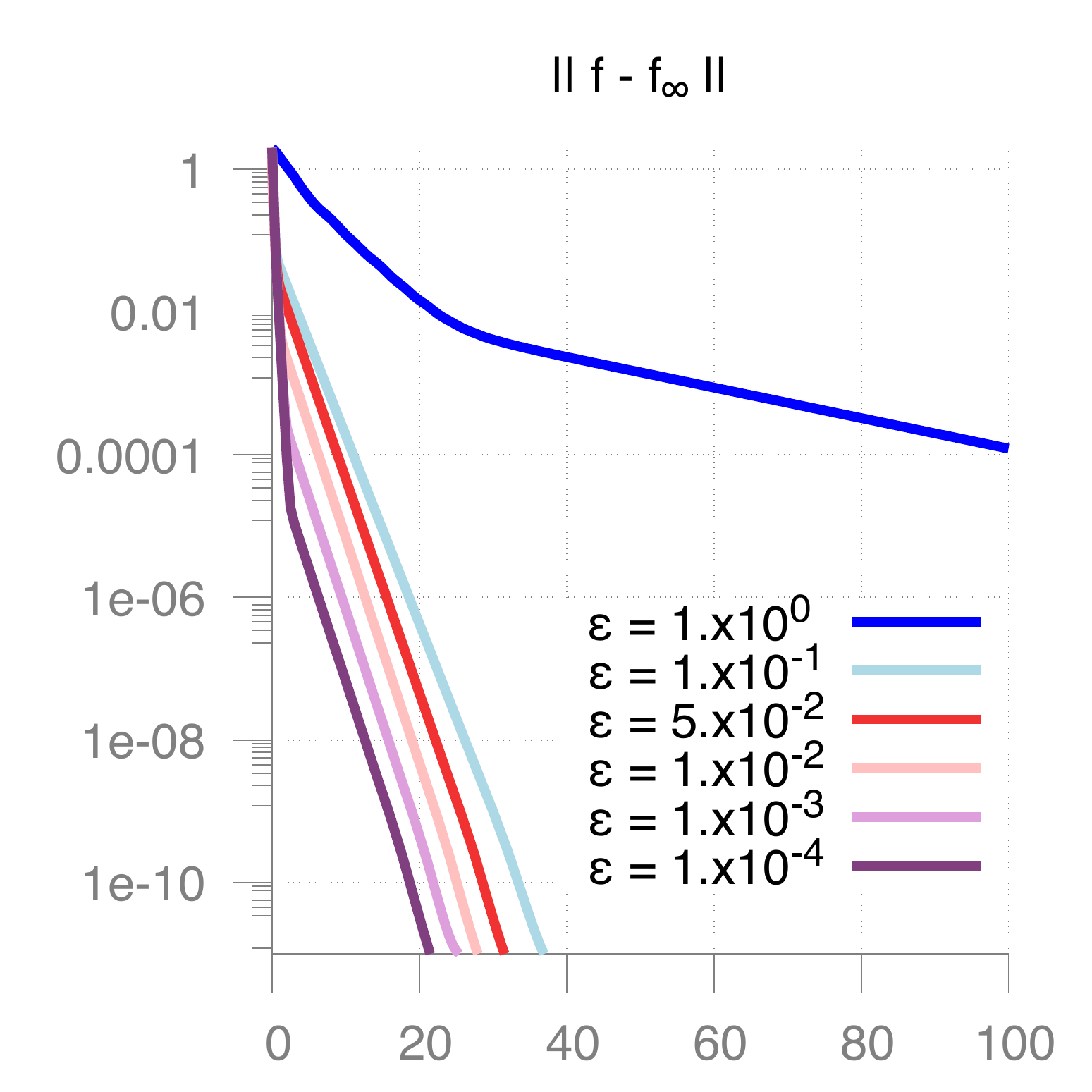}&
                                                         \includegraphics[width=3.2in,clip]{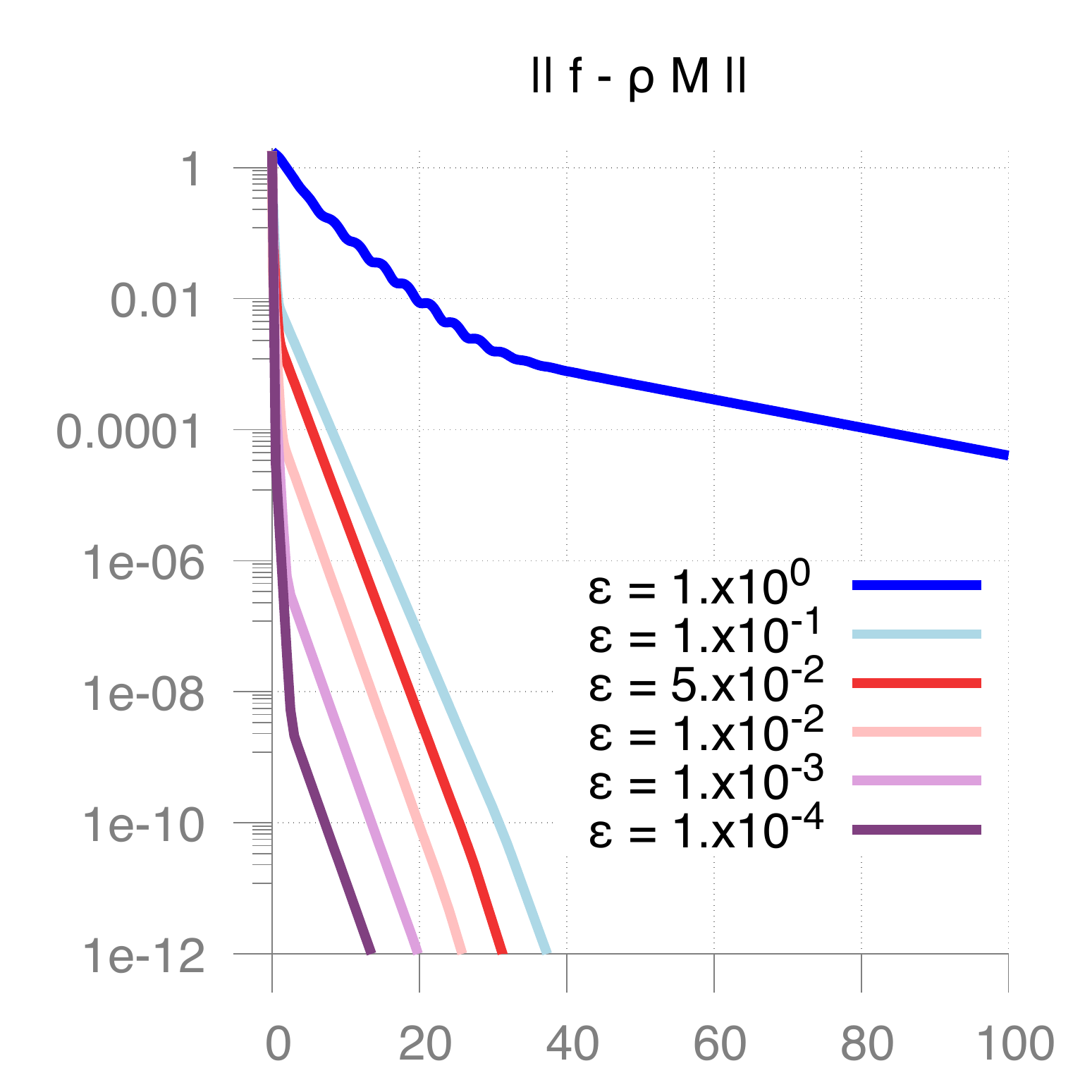}\\
    (a)&(b) \\
    \includegraphics[width=3.2in,clip]{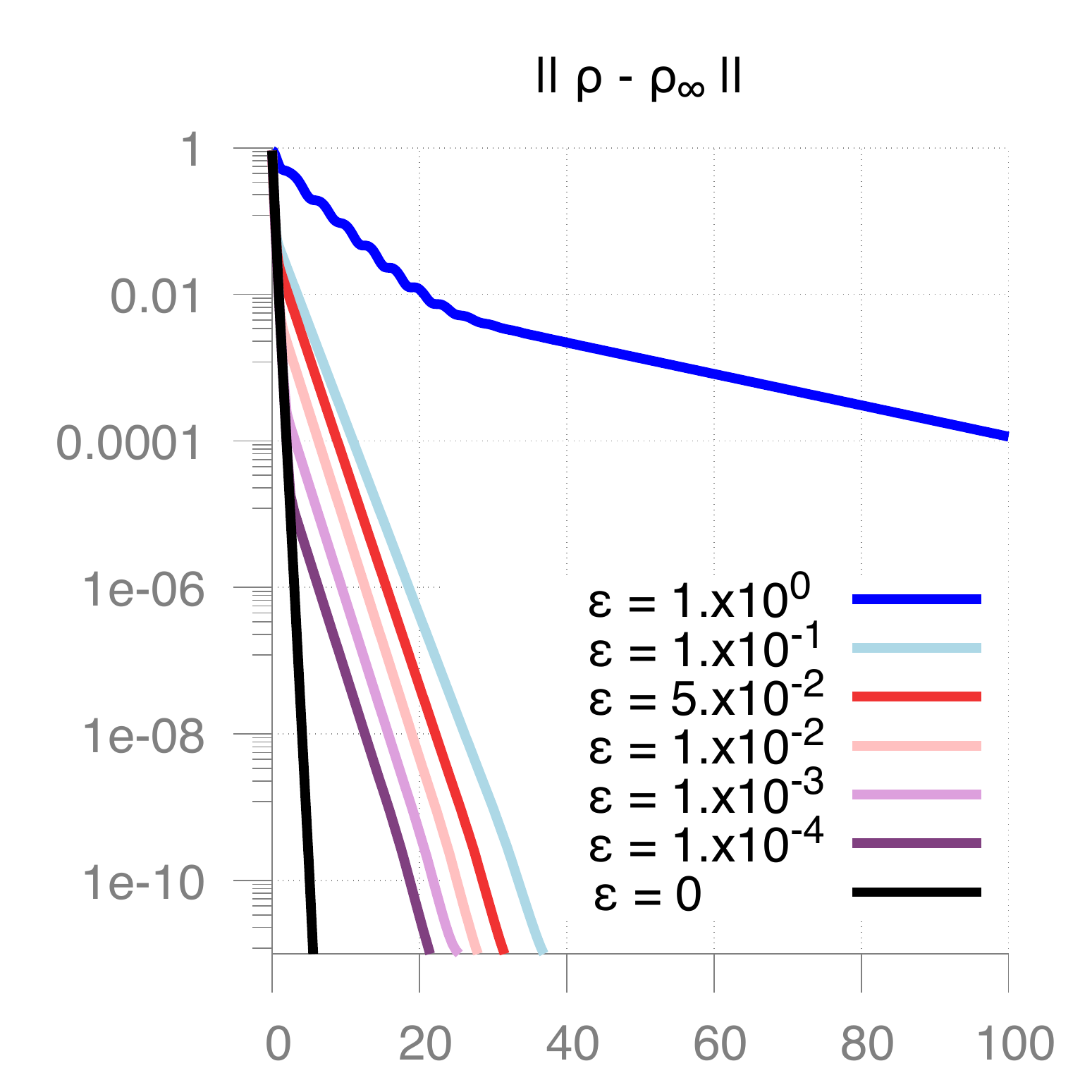}&
                                                      \includegraphics[width=3.2in,clip]{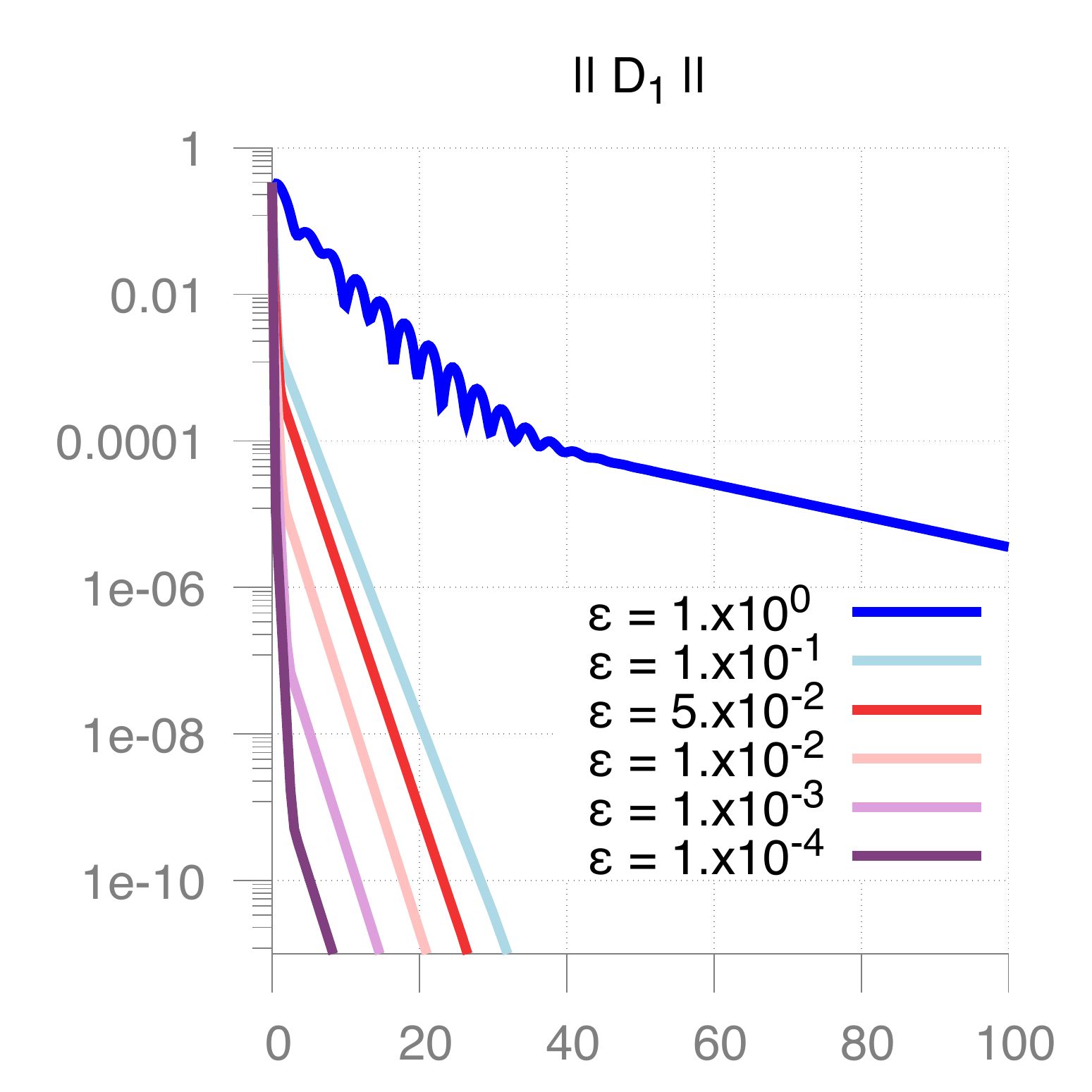}
                                            
    \\
    (c) & (d)
            \end{tabular}
\caption{{\bf Test 2 : shifted Maxwellian.}  time evolution
          in log scale of $(a)$ $\|f-f_\infty\|_{L^2(f_\infty^{-1})}$,
            $(b)$ $\|f-\rho \, \cM \|_{L^2(f_\infty^{-1})}$,  $(c)$
              $\|\rho-\rho_\infty  \|_{L^2(\rho_\infty^{-1})}$ and  $(d)$ $\|D_1  \|_{L^2}$. }
            
	\label{fig:21}
\end{figure}

\begin{figure}
  \centering
  \begin{tabular}{cc}
	\includegraphics[width=3.2in,clip]{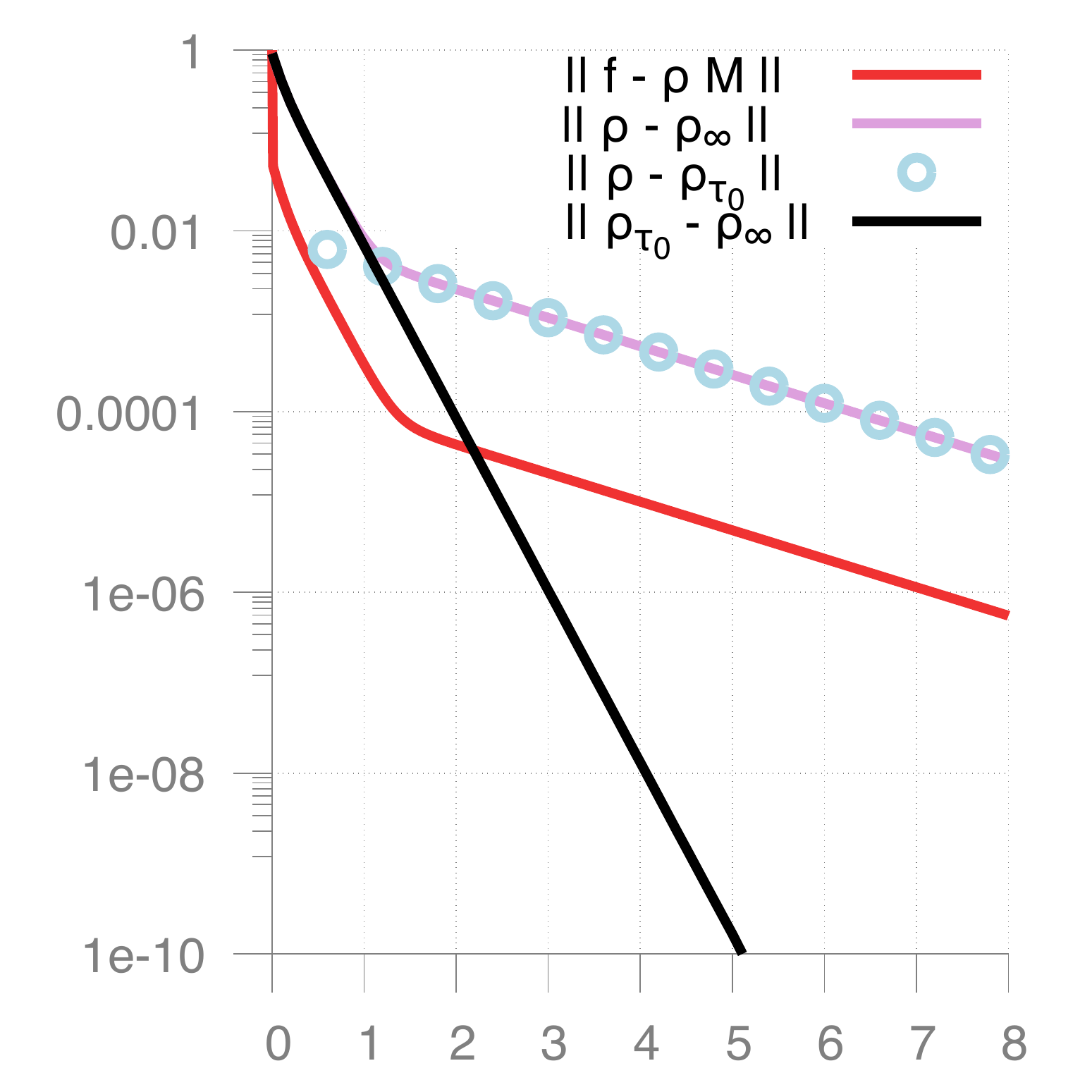}&
	\includegraphics[width=3.2in,clip]{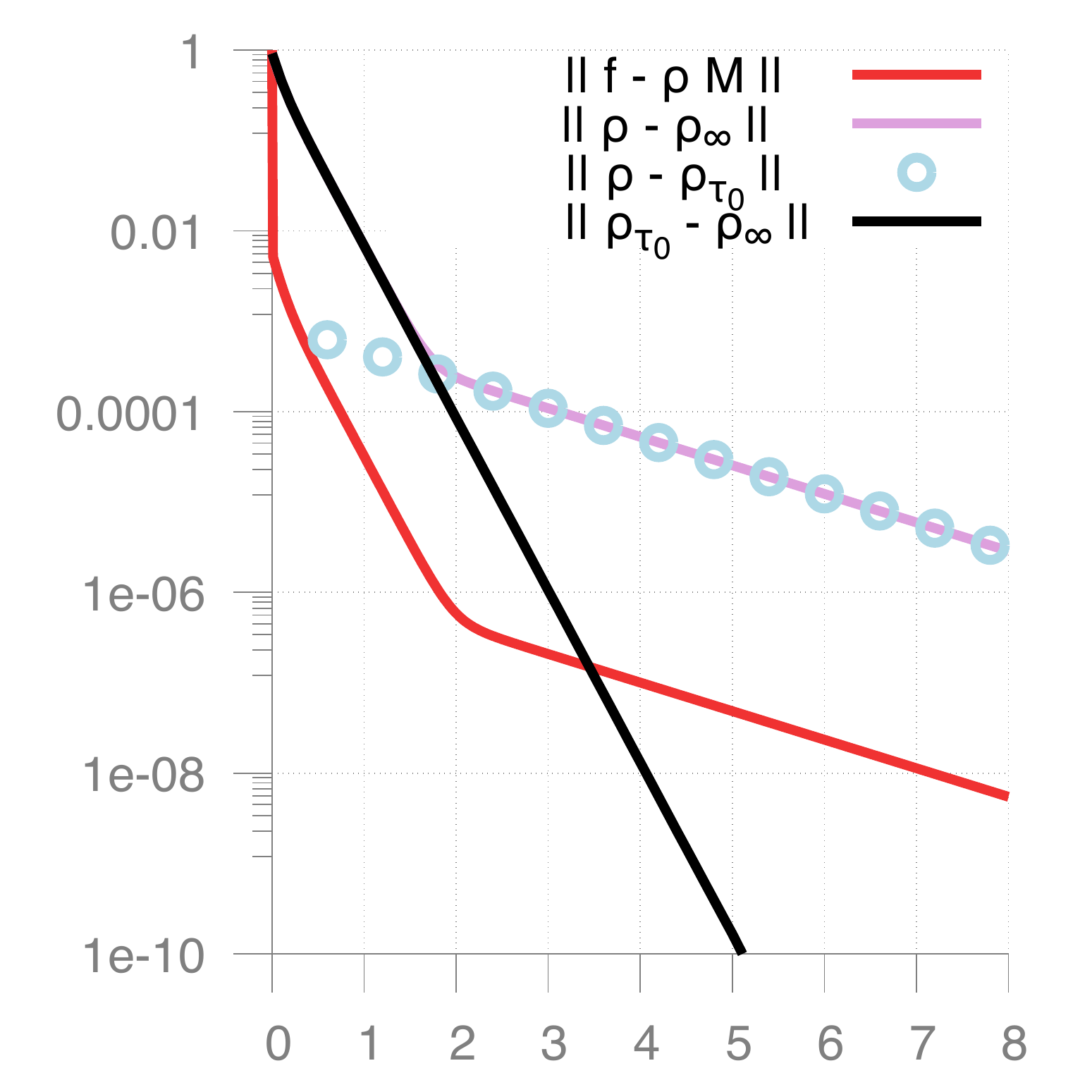}\\
    (a) $\eps=10^{-2}$  &(b) $\eps=10^{-3}$  \\
    \includegraphics[width=3.2in,clip]{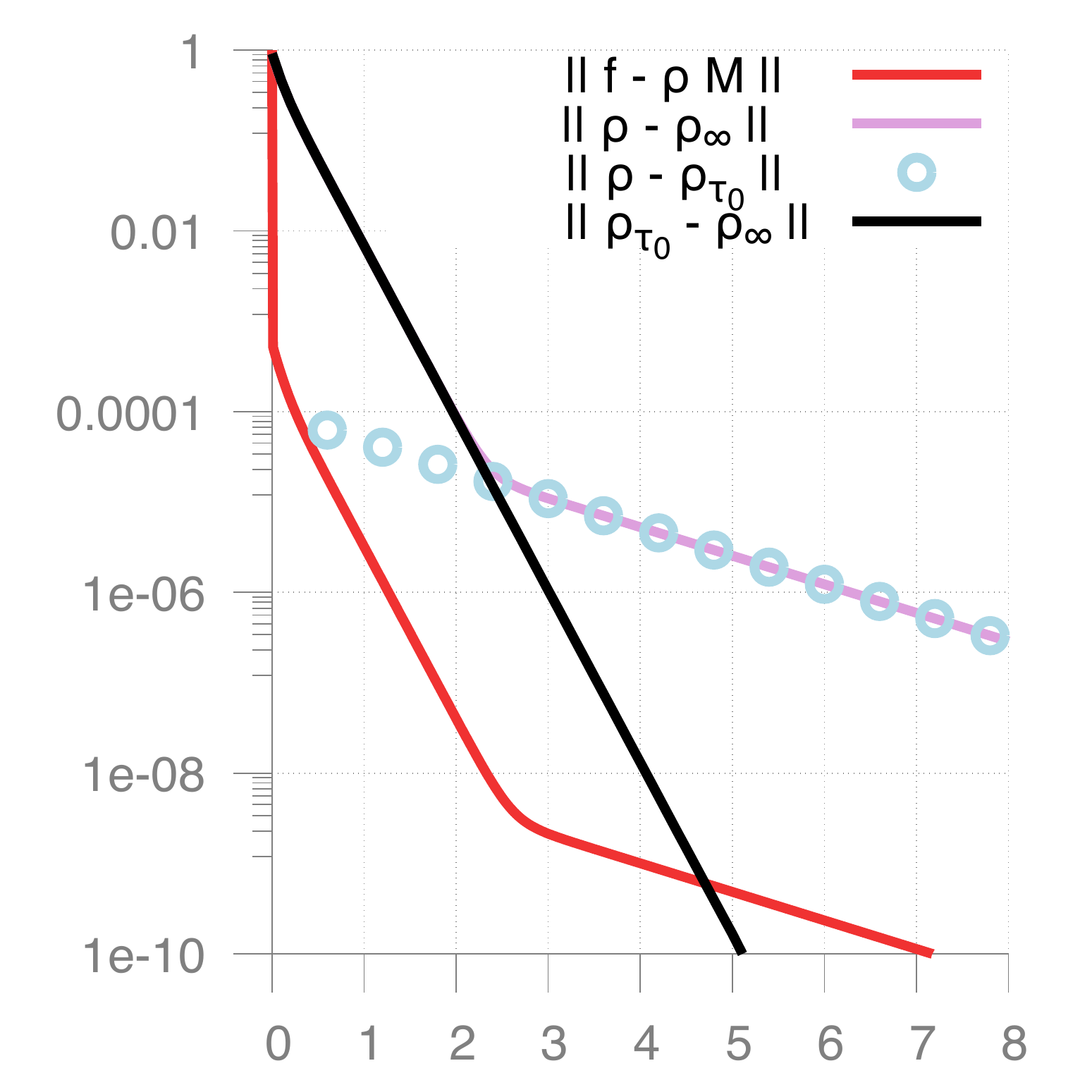}&
	\includegraphics[width=3.2in,clip]{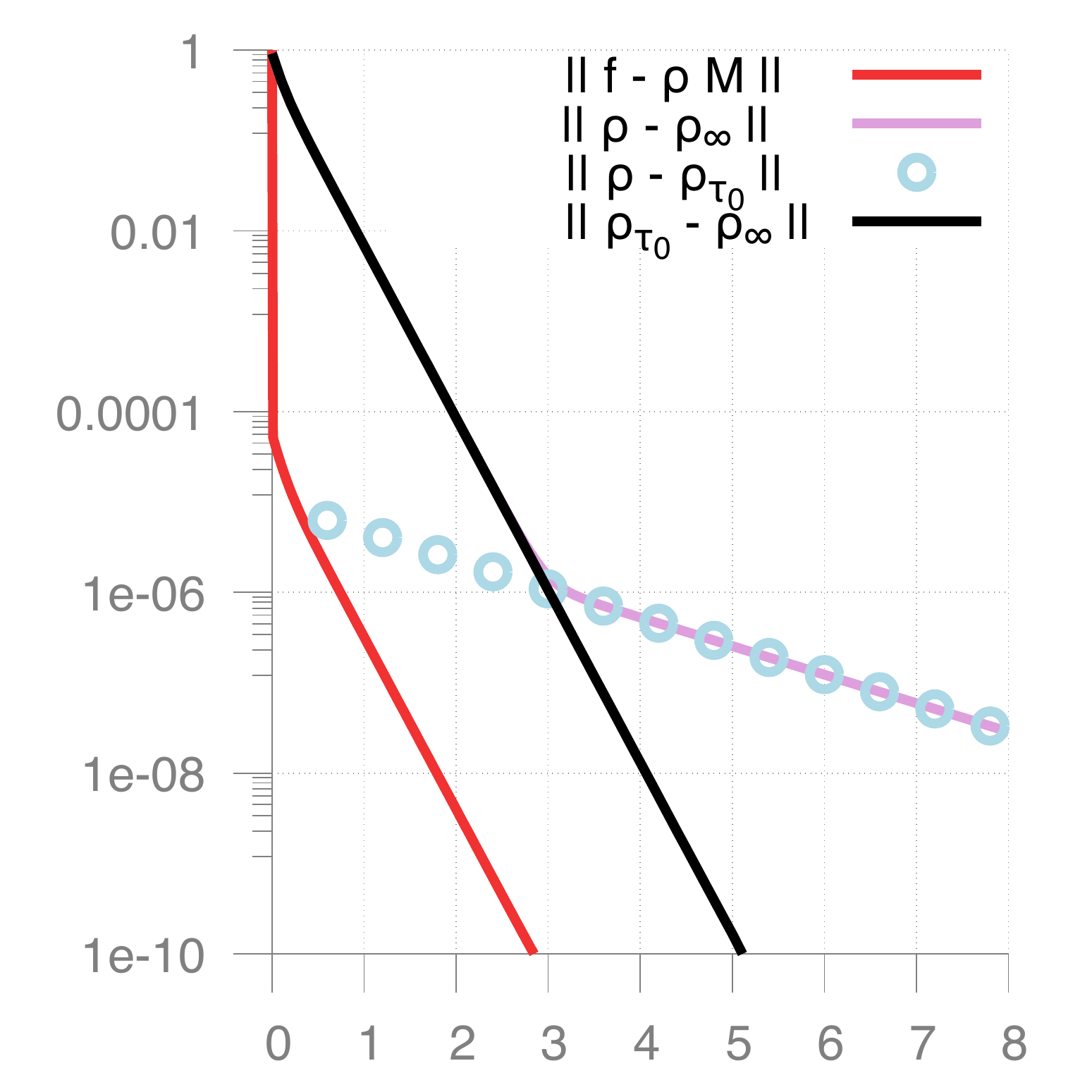}\\
    (c) $\eps=10^{-4}$  & (d) $\eps=10^{-5}$ \\
  \end{tabular}
  \caption{{\bf Test 2 : shifted Maxwellian.}  time evolution
          in log scale of 
             $\|f-\rho \, \cM \|_{L^2(f_\infty^{-1})}$ (red), 
              $\|\rho-\rho_\infty  \|_{L^2(\rho_\infty^{-1})}$
              (pink), $\|\rho-\rho_{\tau_0}
              \|_{L^2(\rho_\infty^{-1})}$ (blue points)
              and 
              $\|\rho_{\tau_0}-\rho_\infty
              \|_{L^2(\rho_\infty^{-1})}$ (black) for  $\eps=10^{-2}, \, 10^{-3}, \, 10^{-4}$ and $10^{-5}$. }
            
	\label{fig:22}
\end{figure}

First, we focus on the case $\eps=1$ displayed in Figure \ref{fig:21},
where we observe that unlike in the previous test, the oscillatory
relaxation stops after a short time and is replaced by a slower but
straight relaxation towards equilibrium. Another interesting comment on Figure \ref{fig:21} is that all the curves associated to value of $\eps$ below $5.10^{-2}$ (red, beige, pink and purple) are parallel. These two features might be explained by a fine spectral analysis of the model at play.\\
We now zoom in to focus on smaller time intervals and  propose a
detailed description of these dynamics in Figure \ref{fig:22}, where
we distinguish three phases constituting a great illustration for
the result presented in item $(i)$ of Theorem \ref{th:main2:h}:
\begin{enumerate}
\item the first phase is the initial time layer, it occurs on
negligible time intervals compared to the time scale chosen in Figure
\ref{fig:22} but it is still visible if we focus on the red curves, reprensenting the norm of $D_{\perp}$, in
plots $(a)$ to $(d)$. As predicted by the first result in $(i)$ of
Theorem \ref{th:main2:h}, higher Hermite modes gathered in the
quantity $D_{\perp}$  undergo a
steep exponential descent with theoretical rate of order
$(\eps^2\,\tau_0)^{-1}$, until they reach a critical level of order
$\eps$;
\item the second phase corresponds to the diffusive
regime where $f$ is close to $\rho_{\tau_0}\,\cM$. Indeed we see that for times ranging from $\sim 0$ up to
$t=1$ in the case $\eps=10^{-2}$ and increasing up to $t=3$ in the case $\eps=10^{-5}$, the red curve, which represents the norm of $D_{\perp}$, is parallel to the pink line corresponding to the norm of $\rho-\rho_{\tau_0}$ which itself coincides with the black curve reprensenting the norm of $\rho_{\tau_0}-\rho_{\infty}$. It indicates that, for a finite amount of time which increases as $\eps$ goes to zero, the kinetic model behaves like the macroscopic one;
\item the last phase is the long time behavior, it starts as the error
  between $\rho_{\tau_0}$ and $\rho$ is of the same order as the error
  between $\rho$ and $\rho_{\infty}$. In Figure \ref{fig:22}
  $(a)$-$(d)$, it corresponds to the intersection between circled blue
  and black lines. As predicted by the second result in $(i)$ of
  Theorem \ref{th:main2:h},  this circled curve, representing the
  error $\|\rho-\rho_{\tau_0}\|$, starts with an ordinate of order
  $\eps$  at time $t=0$, then it decays with a rate proportional to
  $\tau_0$ but smaller than the relaxation rate of the macroscopic
  model. This constitutes a striking illustration of "hypocoercivity"
  phenomenon induced by the transport term proper to kinetic
  equations. During this final phase, the
solution $f$ to \eqref{vlasov} slowly relaxes towards equilibrium. A surprising and unexpected fact is that the transition from diffusive
regime to long time behavior occurs synchronisingly for the spatial density and higher modes in velocity. Indeed,  as it can be observed in plots $(a)$ to $(c)$ of Figure \ref{fig:22},  the inflexions points of the red and the pink curves are almost aligned. 
\end{enumerate}

\section{Conclusion and perspectives}
In the present article, we design a numerical method capable to
capture a rich variety of regimes for a Vlasov-Fokker-Planck equation
with external force field. We prove quantitative estimates for all the
regimes of interest, and do this uniformly with respect to all
parameter at play. We illustrate the robustness of our scheme by
proposing several numerical tests in which we capture a wide variety
of situations (exponential decay with oscillations, transition phase
between diffusive regime an long time behavior, initial time layer, etc ...). Furthermore, we built the method such that it should be easily adaptable in any dimension, at least for cartesian mesh.\\

Two questions arise naturally from this work. The first one is to build on the groundworks layed in this article in order to design a scheme which takes into account non-linear coupling with Poisson for the electric force field. This challenging perspective would be a great improvement since even for the continuous model, there exists to our knowledge very few results which treat the longtime behavior and the diffusive regime with the accuracy proposed in this article. Up to our knowledge, all the works on this subject have restrictions on the dimension of the phase-space and therefore, it would naturally be interesting to propose a method which applies in the physical case $d=3$.\\
Another interesting question arose from our numerical tests, in which we witnessed oscillating behaviors in the solution's relaxation towards equilibrium as well as transition phase between diffusive regime and longtime behavior. It would be of great interest to carry out a fine spectral analysis of the model both at the continuous and the discrete level in order to provide a quantitative description of these phenomena: we may hope for precise and enlightening results due to the simplicity of our model.

%
\section*{Acknowledgement}

Both authors are partially funded by the ANR Project Muffin (ANR-19-CE46-0004).

\bibliographystyle{abbrv}
\bibliography{refer}
\end{document}